\input amstex
\documentstyle{amsppt}
\magnification 1100
\font\tenscr=rsfs10
\font\sevenscr=rsfs7
\font\fivescr=rsfs5
\skewchar\tenscr='177 \skewchar\sevenscr='177 \skewchar\fivescr='177
\newfam\scrfam \textfont\scrfam=\tenscr \scriptfont\scrfam=\sevenscr
\scriptscriptfont\scrfam=\fivescr
\define\scr#1{{\fam\scrfam#1}}
\input xy
\xyoption{all} 
\pagewidth{420pt}
\define\a{\alpha}\define\be{\beta}\define\Gam{\varGamma}
\define\Ph{\varphi}\define\z{\zeta}\define\zn{\zeta_{p^n}}
\define\G{\text{\rm Gal}}\define\Dc{\bold D_{\text{\rm cris}}}\define\La{\varLambda}
\define\Dd{\bold D_{\text{\rm dR}}}

\define\Drig{\bold D_{\text{\rm rig}}}
\define\RG{\bold R\Gamma}
\define\res{\text{\rm res}}
\define \R{\bold R}
\define \boB{\bold B}
\define \Brigplus{\widetilde {\bold B}_{\text{\rm rig}}^+}
\define \A{\bold A}
\define \E{\bold E}
\define \bD{\bold D}

\define \co{\text{\rm c}}
\define \cl{\text{\rm cl}}
\define \Tr{\text{\rm Tr}}
\define \Tot{\text{\rm Tot}}
\define \Iw{\text{\rm Iw}}
\define \Zp{\Bbb Z_p}

\define\ep{\varepsilon}
\define\Bd{\bold B_{\text{\rm dR}}}
\define\Bc{\bold B_{\text{\rm cris}}}\define\Exp{\text{\rm Exp}}

\define\Hom{\text {\rm Hom}}
\define\Ext{\text {\rm Ext}}\define\M {\Cal M}

\define\F{\text{\rm Fil}}
\define \g{\gamma}
\define \gn{\gamma_n}

\define\Ind{\text{\rm Ind}}

\define\Tam{\text{\rm Tam}}
\define\Hi{H_{\text{\rm Iw}}}

\define\Ddagrig{\bold D^{\dagger}_{\text{\rm rig}}}

%\define\Brigdag{\bold B^{\dagger}_{{\text{\rm rig},K}}}
\define\Brigdag{\CR (K)}

\define\iso{\overset{\sim}\to{\rightarrow}}

\define\Rep{\text{\bf Rep}}
\define\Gal{\text{\rm Gal}}
\define\dop{\partial}
\define \HG{\Cal H(\Gamma)}
\define \bExp{\text{\bf E}\text{\rm xp}}
\define \Vc{\bold V_{\text{\rm cris}}}
\define \CR{\Cal R}
\define \CDcris{\scr D_{\text{\rm cris}}}
\define\loc{\text{\rm loc}_p}
\define \sha{\hbox{\text{\bf III}\hskip -11pt\vrule width10pt
depth0pt height0.4pt\hskip 1pt}}
\topmatter \nologo \TagsOnRight
\NoBlackBoxes
\title
{\smc On extra zeros of $p$-adic $L$-functions: the crystalline case}
\endtitle
\author
{ Denis Benois}
\endauthor
\date April 2013
\enddate
\abstract We formulate a conjecture about  extra zeros of $p$-adic $L$-functions at near central points 
which generalizes the conjecture formulated in \cite{Ben2}. We prove that
this conjecture is compatible with Perrin-Riou's theory of $p$-adic $L$-functions.
Namely, using   Nekov\'a\v r's  machinery of Selmer complexes we prove that
our $\scr L$-invariant appears as an additional factor in the Bloch-Kato type formula 
for special values of Perrin-Riou's module of  $L$-functions.
$$
$$
\tenpoint

$\qquad\qquad\qquad\qquad\qquad\qquad\text{\it Nous avons toutefois suppos\'e pour simplifier que}$

$\qquad\qquad\qquad\qquad\qquad\qquad\text{\it les op\'erateurs $1-\Ph$ et $1-p^{-1}\Ph^{-1}$ sont
inversibles}$

$\qquad\qquad\qquad\qquad\qquad\qquad \text{\it laissant les autres cas, pourtant extr\^emement }$

$\qquad\qquad\qquad\qquad\qquad\qquad\text{\it int\'eressant pour plus tard.}$
\newline
$\,$

$\qquad\qquad\qquad\qquad\qquad\qquad\text{\rm Introduction to Chapter III of \cite{PR2}}$
\endabstract
\address Denis Benois, Institut de Math\'ematiques,
Universit\'e Bordeaux I, 351, cours de la Lib\'eration 33405
Talence, France
\endaddress
\email denis.benois\@math.u-bordeaux1.fr
\endemail
\leftheadtext{Denis Benois} \rightheadtext{Extra zeros}
\toc
\nofrills {\bf Table of contents} \widestnumber \head{\S 5.}
\widestnumber\subhead{5.2} \specialhead{}Introduction
\endspecialhead 
\head \S 1. Preliminaries
\endhead
\subhead 1.1. $(\Ph,\Gamma)$-modules
%The rings $\boB$, $\boB^{\dagger}$ and $\Brigdag$
\endsubhead
 %The rings $\Bc$ and $\Bd$
%\endsubhead
\subhead 1.2. Crystalline representations
\endsubhead
\head \S 2. The exponential map 
\endhead
\subhead 2.1. The Bloch-Kato exponential map
\endsubhead
\subhead 2.2. The large exponential map
\endsubhead
\head \S 3. The $\scr L$-invariant
\endhead
\subhead 3.1. Definition of the $\Cal L$-invariant
\endsubhead
\subhead 3.2. $\scr L$-invariant and the large exponential map
\endsubhead
\head \S 4. Special values of $p$-adic $L$-functions 
\endhead
\subhead 4.1. The Bloch-Kato conjecture
\endsubhead
\subhead 4.2. $p$-adic $L$-functions
\endsubhead
\head \S 5. The module of $p$-adic $L$-functions
\endhead
\subhead 5.1. The Selmer complex
\endsubhead
\subhead 5.2. The module of $p$-adic $L$-functions 
\endsubhead
\head  {Appendix. Galois cohomology of $p$-adic representations}
\endhead
\endtoc
\endtopmatter
\document

\head{\bf Introduction}
\endhead

\flushpar 
{\bf 0.1. Extra zeros.} Let $M$ be a pure motive over $\Bbb Q$.  Assume that the complex
$L$-function $L(M,s)$ of $M$ extends to a meromorphic function on the whole
complex plane $\Bbb C$. Fix an odd prime $p$. It is expected that one can construct
$p$-adic analogs of $L(M,s)$ interpolating $p$-adically   algebraic
parts of its special values. This program was realised and the corresponding 
$p$-adic $L$-functions were constructed in many cases, but the general theory  remains conjectural. 
 In \cite{PR2}, Perrin-Riou formulated 
precise conjectures about the existence and arithmetic properties of $p$-adic
$L$-functions in the case then the $p$-adic realisation $V$
of $M$ is crystalline at $p$. Let $\Dc (V)$ denote the filtered Dieudonn\'e module
associated to $V$ by the theory of Fontaine. Let $D$ be a subspace of $\Dc (V)$
of dimension $d_{+}(V)=\dim_{\Bbb Q_p} V^{c=1}$ stable under the action of $\Ph$.
One says that $D$ is regular if one can associate to $D$ a $p$-adic analog of the
six-term exact sequence of Fontaine and Perrin-Riou (see \cite{PR2} for exact definition).
Fix a lattice $T$ of $V$ stable under the action of the Galois group 
and a lattice $N$ of a regular module $D$.
Perrin-Riou conjectured that one can associate to this data a $p$-adic 
$L$-function $L_p(T,N,s)$ satisfying some explicit interpolation property.
Let $r$ denote the order of vanishing of $L(M,s)$ at $s=0$ and let $L^*(M,0)=\lim_{s\to 0} s^{-r}L(M,s).$
Then at $s=0$ the interpolation property 
 writes
$$
\lim_{s\to 0}\frac{L_p(T,N,s)}{s^r} \,=\,
\Cal E(V,D)R_{V,D}(\omega_{V,N})\,\frac{L^*(M,0)}{R_{M,\infty}(\omega_M)}.
$$
Here  $R_{M,\infty} (\omega_M)$ (resp.  $R_{V,D}(\omega_{V,N})$) is the determinant 
of the Beilinson (resp. the $p$-adic regulator) computed in some compatible bases $\omega_M$ and $\omega_{V,V}$ and $\Cal E(V,D)$ is an Euler-like factor given by
$$
\Cal E(V,D)=
\det (1-p^{-1}\Ph^{-1} \mid D)\,\det (1-\Ph \mid \Dc (V)/D).
$$
If either $D^{\Ph=p^{-1}}\ne 0$ or $(\Dc (V)/D)^{\Ph=1}\ne 0$ we have
$\Cal E(V,D)=0$ and the order of vanishing of $L_p(N,T,s)$ should be $>r.$
In this case we say that $L(T,N,s)$ has an extra zero at $s=0.$
The same phenomenon occurs in the  case then $V$ is semistable and non-crystalline
at $p$. An architypical example is provided by elliptic curves
having split multiplicative reduction \cite{MTT}.
Assume that $0$ is a  {\it critical} point for  $L(M,s)$  and that  $H^0(M)=H^0(M^*(1))=0.$
In \cite{Ben2} using the theory of $(\Ph,\Gamma)$-modules we associated to each regular $D$
an invariant $\scr L(V,D)\in \Bbb Q_p$ generalising  both Greenberg's
$\scr L$-invariant \cite{G}  and Fontaine-Mazur's $\scr L$-invariant \cite{M}.
This allows  to formulate a quite general conjecture about the behavior
of $p$-adic $L$-functions at  extra zeros in the spirit of \cite{G}.
To the best of our knowledge this conjecture is actually proved in
the following cases:

1) Kubota-Leopoldt $p$-adic $L$-functions \cite{FG}, \cite{GK}. Here the $\scr L$-invariant
can be interpreted in terms of Gross  $p$-adic regulator\cite{Gs}.

2) Modular forms of even weight \cite{K}, \cite{GS}, \cite{S}. Here the $\scr L$-invariant
coincides with Fontaine-Mazur's $\scr L(f)$.

3) Modular forms of odd weight \cite{Ben3}. The associated $p$-adic representation $V$ is
either crystalline or potentially crystalline at $p$ and we do need the theory
of $(\Ph,\Gamma)$-modules to define the $\scr L$-invariant.

4) Symmetric square of an elliptic curve having either split multiplicative reduction \cite{Ro}
or a good ordinary reduction (Dasgupta, work in progress).
Here $V$ is ordinary  and the $\scr L$-invariant reduces to  Greenberg's construction
\cite{G}.

5) Symmetric powers of CM-modular forms \cite{HL}.
\newline
\newline
{\bf 0.2. Extra zero conjecture.} In this paper we generalise the  conjecture from \cite{Ben2}
to the {\it non critical} point case. Assume that $V$ is crystalline at $p$. Then the weight argument
shows that $\Cal E(V,D)$ can vanish only if $wt(M)=0$ or $-2$.  In particular, we expect that
the interpolation factor does not vanish at $s=0$ if $wt (M)=-1$ i.e. that the $p$-adic $L$-function can not
have an extra zero at the central point in the good reduction case. To fix ideas assume that
$wt (M)\leqslant -2$ and that $M$ has no subquotients isomorphic to $\Bbb Q(1)$
\footnote{The last condition is not really essential and can be suppressed}.
Then $D$ is regular if and only if the associated $p$-adic regultor map
$$
r_{V,D}\,\,:\,\, H_f^1(V)@>>>\Dc (V)/(\F^0\Dc (V)+D)
$$
is an isomorphism. The semisimplicity of $\Ph \,:\,\Dc (V)@>>>\Dc (V)$ (which conjecturally
always holds) allows to decompose $D$ into a direct sum
$$
D=D_{-1}\oplus D^{\Ph=p^{-1}}.
$$ 
Under some mild assumptions (see 3.1.2 and 4.1.2 below) we associate to $D$ an $\scr L$-invariant
$\scr L(V,D)$ which is a direct generalization of the main construction of \cite{Ben2}.
The Beilinson-Deligne conjecture predicts that $L(M,s)$ does not vanish at $s=0$ and that
$L(M^*(1),s)$ has a zero of order $r=\dim_{\Bbb Q_p}H^1_f(V)$ at $s=0.$ 
%Denote by $R_{M,\infty}(\omega_M)$ and   $R_{V,D}(\omega_{V,N})$ the determinants of
%the Beilinson regulator $r_{M,\infty}$  and the $p$-adic regulator $r_{V,D}$  computed in some %compatibles bases $\omega_M$ and $\omega_{V,D}$ (see section 4.2).
We propose 
the following conjecture:
\newline
{\,}

\proclaim{ Extra zero conjecture} Let $D$ be a regular subspace of $\Dc (V)$ and
let $e=\dim_{\Bbb Q_p}(D^{\Ph=p^{-1}}).$ 
Then

1) The $p$-adic $L$-function $L_p(T,N,s)$ has a zero of order $e$ at $s=0$ and
$$
\frac{L^*_p(T,N,0)}{R_{V,D}(\omega_{V,N})} \,=\,
-\scr L(V,D)\,\Cal E^+(V,D)\,\frac{L(M,0)}{R_{M,\infty}(\omega_M)}.
$$
 
2) Let $D^{\perp}$ denote the orthogonal complement to $D$ under 
the canonical duality 
\linebreak
$\Dc (V)\times \Dc (V^*(1))@>>>\Bbb Q_p.$
The $p$-adic $L$-function $L_p(T^*(1),N^{\perp},s)$ has a zero of order
$e+r$ where $r=\dim_{\Bbb Q}H^1_f(V)$ at $s=0$ and 
$$
\frac{L_p^*(T^*(1),N^{\perp},0)}{R_{V^*(1),D^{\perp}}(\omega_{V^*(1),N^{\perp}})} \,=\,
\scr L(V,D)\,
\Cal E^+(V^*(1),D^{\perp})\,\frac{L^*(M^*(1),0)}{R_{M^*(1),\infty}(\omega_{M^*(1)})}.
$$
In the both cases
$$
\Cal E^+(V,D)=\Cal E^+(V^*(1),D^{\perp})=\det (1-p^{-1}\Ph^{-1} \mid D_{-1})\,
\det (1-p^{-1}\Ph^{-1} \mid \Dc (V^*(1))).
$$
\endproclaim

\flushpar
{\bf Remarks.} 1) $\Cal E^+(V,D)$ is obtained from $\Cal E(V,D)$ by excluding
zero factors. It can be also written in the form
$$
\Cal E^+(V,D)= ,E_p^*(V,1)\,{\det}_{\Bbb Q_p} \left (\frac{1-p^{-1}\Ph^{-1}}{1-\Ph} \,\vert D_{-1} \right )
$$
where $E_p(V,t)={\det} (1-\Ph t\mid \Dc (V))$ is the Euler factor at $p$ and   
$E_p^*(V,t)\,=\,E_p(V,t)\,\left (1-\dsize\frac{t}{p}\right )^{-e}.$

2) Assume that $H^1_f(V)=0$. Since $H^1_f(V^*(1))$ should also vanish
by the weight reason,  our conjecture in this cases reduces to the conjecture 2.3.2 from \cite{Ben2}.

3) The regularity of $D$ supposes that the localisation  $H_f^1(V) @>>>H^1_f(\Bbb Q_p,V)$
is injective.  Jannsen's conjecture (precised by Bloch and Kato) says that the $p$-adic
realisation map $H^1_f(M)\otimes \Bbb Q_p @>>> H^1_f(V)$ is an isomorphism. 
The composition $H^1_f(M)@>>>H^1_f(\Bbb Q_p,V)$ of these two maps  is essentially the syntomic regulator.  Its injectivity seems to be  a difficult open  problem. 
\newline
\newline
{\bf 0.3. Selmer complexes and Perrin-Riou's theory.} In the last part of the paper
we show that our extra zero conjecture is compatible with the Main Conjecture
of Iwasawa theory as formulated in \cite{PR2}. The main technical tool here is the descent
theory for Selmer complexes \cite{Ne2}. We hope that the approach to Perrin-Riou's theory
based on the formalism of Selmer complexes can be of independent interest.

For a profinite group $G$ and a continuous $G$-module $X$ we denote by $C_c^\bullet (G,X)$
the standard complex of continuous cochains. Let $S$ be a finite set of primes
containing $p$. Denote by  $G_S$ the Galois group of the maximal algebraic extension of $\Bbb Q$
unramified outside $S\cup \{\infty\}$.  Set
$\RG_S(X)= C^{\bullet}_c (G_S,X)$ and $\RG (\Bbb Q_v,X)=C_c^\bullet (G_v,X)$, where 
$G_v$ is the absolute Galois group  of  $\Bbb Q_v$. Let $\Gamma$ be the Galois group
of the cyclotomic $p$-extension $\Bbb Q(\zeta_{p^{\infty}})/\Bbb Q),$ 
$\Gamma_1=\Gal (\Bbb Q(\zeta_{p^{\infty}})/\Bbb Q (\zeta_p)$ and $\Delta =\Gal (\Bbb Q(\zeta_p)/\Bbb Q)$.
Let $\Lambda (\Gamma) =\Bbb Z_p[[\Gamma]]$ denote the Iwasawa algebra of $\Gamma .$
Each $\Lambda (\Gamma)$-module  $X$ decomposes into the direct sum of its isotypical components
$X=\dsize\underset{\eta \in \hat\Delta}\to \oplus X^{(\eta)}$ and we denote by 
$X^{(\eta_0)}$ the component which corresponds to the trivial character $\eta_0.$
Set $\Lambda =\Lambda (\Gamma)^{(\eta_0)}.$
Let $\Cal H$ denote the algebra of power series with coefficients in $\Bbb Q_p$
which converge on the open unit disk. We will denote again by $\Cal H$ the
associated large Iwasawa algebra $\Cal H(\Gamma_1).$
In this paper we consider only the trivial character component of the module of
$p$-adic $L$-functions because it is sufficient for applications to trivial zeros,
but in the general case the construction is exactly the same.
We keep notation and assumptions of section 0.2
 Assume that the weak Leopoldt conjecture holds for  $(V, \eta_0)$ and $(V^*(1),\eta_0).$ 
We consider global and local   Iwasawa cohomology 
$\RG_{\text{\rm Iw},S} (T)\,=\RG_S(\Lambda (\Gamma)\otimes_{\Bbb Z_p} T)^\iota)$
and $\RG_{\text{\rm Iw}}(\Bbb Q_v,T)\,=\,\RG (\Bbb Q_v,(\Lambda (\Gamma)\otimes_{\Bbb Z_p} T)^\iota)$
where $\iota$ is the canonical involution on $\Lambda (\Gamma).$ Let $D$ be a regular submodule
of $\Dc (V).$  For each non archimedian
place $v$ we define a local condition at $v$ in the sense of \cite{Ne2} as follows.
If $v\ne p$ we use the unramified local condition which is defined by
$$
\RG_{\text{\rm Iw},f}^{(\eta_0)}(\Bbb Q_v,N,T)= \RG_{\text{\rm Iw},f}^{(\eta_0)}(\Bbb Q_v,T)\,=\,
\left [T^{I_v}\otimes \Lambda^{\iota} @>f_v>> T^{I_v}\otimes \Lambda^{\iota} \right ]
$$
where $I_v$ is the inertia subgroup at $v$ and $f_v$ is the geometric Frobenius.
If $v=p$ we define 
$$ 
\RG_{\text{\rm Iw},f}^{(\eta_0)}(\Bbb Q_v,N,T) = (N\otimes \Lambda) [-1].
$$
The derived version of  the large exponential map $\text{\rm Exp}_{V,h},$ $h\gg 0$ (see \cite{PR1}) 
gives a morphism
$$
\RG_{\text{\rm Iw},f}^{(\eta_0)}(\Bbb Q_v,N,T) @>>>\RG_{\text{\rm Iw}}^{(\eta_0)}(\Bbb Q_p,T)
\otimes \Cal H.
$$
Therefore we have a diagram
$$
\xymatrix{
\RG_{\text{\rm Iw},S}^{(\eta_0)}(T)\otimes_{\Lambda}\Cal H
\ar[r] & \underset{v\in S} \to \oplus \RG_{\text{\rm Iw}}^{(\eta_0)}(\Bbb Q_v,T)\otimes_{\Lambda} \Cal H\\
 & \left (\underset{v\in S}\to \oplus \RG_{\text{\rm Iw},f}^{(\eta_0)}(\Bbb Q_v,N,T)
\right ) \otimes_{\Lambda} \Cal H\,\,.
 \ar[u]
}
$$
Let $\RG_{\text{\rm Iw},h}^{(\eta_0)}(D,V)$ denote the Selmer complex associated to this data.
By definition it sits in the distinguished triangle
$$
\bold R\Gamma_{\text{\rm Iw},S}^{(\eta_0)}(D,V)@>>>
\left (\RG_{\text{\rm Iw},S}^{(\eta_0)}(V) \oplus \left (\underset{v\in S}\to \oplus
\RG_{\text{\rm Iw},f}^{(\eta_0)}(\Bbb Q_v,D,V)\right )\right )\otimes{\Cal H} @>>>
\left (\underset{v\in S}\to \oplus
\RG_{\text{\rm Iw}}^{(\eta_0)}(\Bbb Q_v,V)\right )\otimes \Cal H . \tag{0.1}
$$
Define
$$
\Delta_{\text{\rm Iw},h} (N,T)\,=\,
{\det}^{-1}_{\Lambda}\left (\RG_{\text{\rm Iw},S}^{(\eta_0)}(T) \oplus \left (\underset{v\in S}\to \oplus
\RG_{\text{\rm Iw},f}^{(\eta_0)}(\Bbb Q_v,N,T)\right )\right ) \otimes 
{\det}_{\Lambda} \left (\underset{v\in S}\to \oplus
\RG_{\text{\rm Iw}}^{(\eta_0)}(\Bbb Q_v,T)\right ).
$$
Our results can be summarized as follows (see Theorems 5.1.3, 5.2.5 and Corollary 5.2.7).

\proclaim{Theorem 1} Assume that $\scr L(V,D)\ne 0.$  
Then

i) The cohomology  $\bold R^i\Gamma_{\text{\rm Iw},h}^{(\eta_0)}(D,V)$ are $\Cal H$-torsion modules for all $i$.

ii)  $\bold R^i\Gamma_{\text{\rm Iw},h}^{(\eta_0)}(D,V)=0$ for  $i\ne 2,3$ and
$$
\bold R^3\Gamma_{\text{\rm Iw},h}^{(\eta_0)}(D,V)\simeq \left(H^0(\Bbb Q(\zeta_{p^{\infty}}),V^*(1))^*
\right )^{(\eta_0)}\otimes_{\Lambda}\Cal H.
$$

iii) The complex $ \bold R\Gamma_{\text{\rm Iw},h}^{(\eta_0)}(D,V)$ is semisimple i.e.
for each $i$ the natural map 
$$
\bold R^i\Gamma_{\text{\rm Iw},h}^{(\eta_0)}(D,V)^{\Gamma} @>>>
\bold R^i\Gamma_{\text{\rm Iw},h}^{(\eta_0)}(D,V)_{\Gamma}
$$
is an isomorphism.
\endproclaim 

Assume that $\scr L(V,D)\ne 0.$ Let $\Cal K$ be the field of fractions of $\Cal H.$
Then Theorem 1 together with (0.1) define an injective map
$$
i_{V,\text{\rm Iw},h}\,:\,\Delta_{\text{\rm Iw},h}(N,T) @>>> \Cal K 
$$  
and the module of $p$-adic $L$-functions is defined as 
$$
\bold L^{(\eta_0)}_{\text{\rm Iw},h}(N,T) =
i_{V,\text{\rm Iw},h} (\Delta_{\text{\rm Iw},h}(N,T)) \subset \Cal K.
$$
Let $\gamma_1$ be a fixed generator of $\Gamma_1$.
Choose a generator $f(\gamma_1-1)$ of the free $\Lambda$-module  $\bold L^{(\eta_0)}_{\text{\rm Iw},h}(N,T)$ and define a meromorphic
$p$-adic function
$$
L_{\text{\rm Iw},h}(T,N,s)=f(\chi (\gamma_1)^s-1),
$$
where $\chi \,\,:\,\,\Gamma @>>>\Bbb Z_p^*$ is the cyclotomic character.

\proclaim{Theorem 2} Assume that $\scr L(V,D)\ne 0.$ Then

1) The $p$-adic $L$-function $L_{\text{\rm Iw},h}(T,N,s)$ has a zero of order $e=\dim_{\Bbb Q_p}(D^{\Ph=p^{-1}})$ at $s=0.$

2) One has 
$$
\frac{L_{\text{\rm Iw},h}^*(T,N,0)}{R_{V,D}(\omega_{T,N})}
\sim_p \Gamma (h)^{d_{+}(V)} \,\scr L (V,D)\, 
\Cal E^+(V,D)
\,\,\frac{\#\sha(T^*(1))\,\Tam_{\omega_M}^0 (T)}{\#H^0_S(V/T)\,\#H^0_S(V^*(1)/T^*(1))},
$$
where $\sha(T^*(1))$ is the Tate-Shafarevich group of Bloch-Kato \cite{BK} and   $\Tam_{\omega_M}^0 (T)$ is the product of local Tamagawa numbers of $T$.

\endproclaim

\flushpar
{\bf Remarks.} 1) Using the compatibility of Perrin-Riou's theory with the functional equation
we obtain analogous results for the $L_p(T^*(1),N^{\perp},s)$ (see section 5.2.9).

2) It $\Dc (V)^{\Ph=1}=\Dc (V)^{\Ph=p^{-1}}=0$  the phenomenon
of extra zeros does not appear, $\scr L(V,D)=1$  and Theorem 2 was proved in
\cite{PR2}, Theorem 3.6.5. We remark that even in this case our proof is different. 
We compare the leading term of $L_{\text{\rm Iw},h}^*(T,N,s)$
with the trivialisation $i_{\omega_M,p} \,:\,\Delta_{EP}(T) @>>> \Bbb Q_p$ of the 
Euler-Poincar\'e line $\Delta_{EP}(T)$ (see \cite{F3}) and show
that in compatible bases one has
$$
\frac{L_{\text{\rm Iw},h}^*(T,N,0)}{R_{V,D}(\omega_{V,N})}
\sim_p \Gamma (h)^{d_{+}(V)} \,\scr L (V,D)\, 
\Cal E^+(V,D)
\,\,i_{\omega_M, p}\,(\Delta_{\text{\rm EP}} (T)) \tag{0.2}
$$
(see Theorem 5.2.5). Now Theorem 2 follows from the well known computation
of $i_{\omega_M, p}\,(\Delta_{\text{\rm EP}} (T))$ in terms of the Tate-Shafarevich group
and Tamagawa numbers (\cite{FP}, Chapitre II).

3) Let $E/\Bbb Q$ be an elliptic curve having  good reduction at $p.$ 
Consider the $p$-adic representation  $V=\text{\rm Sym}^2 (T_p(E))\otimes \Bbb Q_p$, where $T_p(E)$ is the $p$-adic Tate module 
of $E.$ It is easy to see that  $D=\Dc (V)^{\Ph=p^{-1}}$ is  one dimensional. In  this case some versions of Theorem 2 
were proved in \cite{PR3} and \cite{D} with  an ad hoc definition  of the $\scr L$-invariant.
Remark that  $p$-adic $L$-functions attached to the symmetric square
of a newform were constructed by Dabrowski and Delbourgo \cite{DD}.

4) Another approach to Iwasawa theory in the non-ordinary case was developped by Pottharst in 
\cite{Pt1}, \cite{Pt2}.  Pottharst  uses the formalism of Selmer complexes
but works with local conditions coming from submodules of 
the $(\Ph,\Gamma)$-module associated to $V$ rather then with the large exponential map.  
This approach has many advantages, in particular it allows to develop
an interesting theory for representations which are not necessarily crystalline.
Nevertheless it seems that the large exponential map is crucial   for the study of extra zeros
at least in the good reduction case.

5) The Main conjecture of Iwasawa theory \cite{PR2}, \cite{C2} says that the analytic 
$p$-adic $L$-function $L_p(N,T,s)$ multiplied by a simple explicit $\Gamma$-factor depending on $h$ can be written in the form  $L_p(N,T,s)=f(\chi (\gamma_1)^s-1)$ for an appropriate generator  $f(\gamma_1-1)$ of 
$\bold L^{(\eta_0)}_{\text{\rm Iw},h}(N,T).$ Therefore the main conjecture implies 
Bloch-Kato style formulas for special values of $L_p(N,T,s).$ We remark that the Bloch-Kato conjecture
predicts that 
$$
\frac{L^*(M,0)}{R_{M,\infty}(\omega_M)} \sim_p \frac{\#\sha(T^*(1))\,\Tam_{\omega_M}^0 (T)}{\#H^0_S(V/T)\,\#H^0_S(V^*(1)/T^*(1))}
$$
and therefore Theorem 2 implies the compatibility of our extra zero conjecture with 
the Main conjecture. Note that this also follows directly from (0.2) if we use
the formalism of Fontaine and Perrin-Riou \cite{F3} to formulate Bloch-Kato conjectures.
\newline
\newline
{\bf 0.4.} The organisation of the paper is as follows. In \S1 we review the theory of $(\Ph,\Gamma)$-modules which is the main technical tool  in our definition of the $\scr L$-invariant.
We also give the derived version of  computation of Galois cohomology in terms of $(\Ph,\Gamma)$-modules.
This  follows easily from the results of Herr \cite{H1} and Liu \cite{Li} and the proofs are 
placed in Appendix. Similar results can be found in \cite{Pt1}, \cite{Pt2}.
In \S2 we recall preliminaries on the Bloch-Kato exponential map and review the construction
of the large exponential map of Perrin-Riou given by Berger \cite{Ber3} using again the basic
language of derived categories. 
The $\Cal L$-invariant is constructed is section 3.1. In section 3.2 we relate this construction
to the derivative of the large exponential map. This result  plays a key role in the proof
of Theorem 2. The extra zero conjecture is formulated in \S4. In \S5 we interpret
Perrin-Riou's theory in terms of Selmer complexes and prove Theorems 1 and 2.  
\newline
\newline  
{\bf Acknowledgements.} I am very grateful to  Jan Nekov\'a\v r and Daniel Delbourgo for several interesting discussions and comments concerning this work. 
\newpage

\head {\bf \S1. Preliminaries}
\endhead

\flushpar
{\bf 1.1. $(\Ph,\Gamma)$-modules.}
\newline
{\bf 1.1.1. The Robba ring} (see \cite{Ber1},\cite{C3}).
In this section  $K$ is a finite unramified extension of $\Bbb Q_p$ with residue field $k_K$, $O_K$ its ring of integers,
and $\sigma$ the absolute Frobenius of $K$. Let 
$\overline K$ an algebraic closure of $K$, $G_K=\text{Gal}(\bar K/K)$ and  $C$ the completion of $\overline K .$ 
Let  $v_p\,\,:\,\,C@>>>\Bbb R\cup\{\infty\}$ denote the $p$-adic valuation normalized so that $v_p(p)=1$ and 
set $\vert x\vert_p=\left (\frac{1}{p}\right )^{v_p(x)}.$ Write $B(r,1)$ for the $p$-adic annulus 
$B(r,1)=\{ x\in C \,\mid \, r\leqslant \vert x\vert <1 \}.$
As usually,  $\mu_{p^n}$  denotes the group of $p^n$-th roots of  unity.
Fix a system of primitive roots of unity   $\ep=(\zeta_{p^n})_{n\geqslant 0}$,
$\,\zeta_{p^n} \in \mu_{p^n} $ such that $\zeta_{p^n}^p=\zeta_{p^{n-1}}$ for all $n$.
Set $K_n=K(\zeta_{p^n})$, $K_{\infty}= \bigcup_{n=0}^{\infty}K_n$, $H_K=\Gal (\bar K/K_\infty)$, $\Gam =\G(K_{\infty}/K)$
and denote by $\chi \,:\,\Gam @>>>\Bbb Z_p^*$ the cyclotomic character.

Set 
$$
\widetilde {\bold E}^+=\varprojlim_{x\mapsto x^p} O_C/\,p\,O_C\,=\,\{x=(x_0,x_1,\ldots ,x_n,\ldots )\,\mid \,
x_i^p=x_i \,\,\forall i\in \Bbb N\}.
$$
Let $\hat x_n\in O_C$ be a lifting of $x_n$.
Then for all $m\geqslant 0$ the sequence $\hat x_{m+n}^{p^n}$ converges to 
$x^{(m)}=\lim_{n\to \infty} \hat x_{m+n}^{p^n}\in O_C$
which does not depend on the choice of  liftings.
 The ring $\widetilde {\bold E}^+$ equipped with the valuation $v_{\bold E}(x)=v_p(x^{(0)})$ is a
complete local ring of characteristic $p$ with  residue field
 $\bar k_K$. Moreover it is integrally closed in his field
of fractions $\widetilde {\bold E}=\text {\rm Fr}(\widetilde {\bold E}^+)$. 

Let $\widetilde \A=W(\widetilde \E)$ be the ring of Witt vectors with coefficients
in $\widetilde \E$. Denote by $[\,\,]\,:\,\widetilde \E@>>>W(\widetilde \E)$  the Teichmuller lift.
Any $u=(u_0,u_1,\ldots )\in \widetilde \A$ can be written in the form
$$
u=\underset{n=0}\to{\overset{\infty}\to \sum} [u^{p^{-n}}]p^n.
$$

Set $\pi=[\ep]-1$, $\A_{K_0}^+=O_{K_0} [[\pi]]$ and denote by $\A_{K}$  the $p$-adic completion
of $\A_{K}^+\left [1/{\pi}\right ]$. 
%This is
%a discrete valuation ring with residue field $\E_{K_0}=k_K((\ep-1))$. 
Let $\widetilde\boB= \widetilde \A\left [{1}/{p}\right ]$,  $\boB_{K}=\A_{K}\left [{1}/{p}\right ]$ and
let $\boB$ denote  the completion of the maximal unramified extension of $\boB_{K}$ in $\widetilde \boB$.
Set $\A=\boB\cap \widetilde \A$, $\widetilde \A^+=W(\E^+)$, $\A^+= \widetilde \A^+\cap \A$ and $\boB^+=\A^+\left [{1}/{p}\right ].$
All these rings are endowed with natural  actions of the Galois group $G_K$ and   Frobenius $\Ph$.

Set $\A_K=\A^{H_K}$ and  $\boB_K=\A_K\left [{1}/{p}\right ].$
Remark that $\Gamma$ and $\Ph$ act on $\boB_{K}$ by
$$
\aligned
& \tau (\pi)=(1+\pi)^{\chi (\tau)}-1,\qquad \tau \in \Gamma\\
&\Ph (\pi)=(1+\pi)^p-1.
\endaligned
$$
For any $r>0$ define
$$
\widetilde {\bold B}^{\dagger,r}\,=\,\left \{ x\in \widetilde
{\bold B}\,\,|\,\, \lim_{k\to +\infty} \left (
v_{\E}(x_k)\,+\,\dsize \frac{pr}{p-1}\,k\right )\,=\,+\infty
\right \}.
$$
Set ${\bold B}^{\dagger,r}=\boB \cap \widetilde\boB^{\dagger,r}$, 
$\boB_{K}^{\dagger,r}=\boB_{K} \cap \boB^{\dagger,r}$, 
${\bold B}^{\dagger}=\underset{r>0}\to \cup \boB^{\dagger,r}$
${\bold A}^{\dagger}=\bold A \cap {\bold B}^{\dagger}$
and $\bold B^\dag_K=\underset{r>0}\to \cup \boB_{K}^{\dagger,r}$.

It can be shown that for any $r\geqslant  p-1$

$$
\boB_{K}^{\dagger,r}=\left \{ f(\pi)=\sum_{k\in \Bbb Z}
a_k\pi^k\,\mid \, \text{\rm $a_k\in K$ and $f$ is holomorphic and bounded on $B(r,1)$} \right \}.
$$
Define
$$
\boB^{\dag,r}_{\text{rig},K}\,=\,\left \{ f(\pi)=\sum_{k\in \Bbb Z}
a_k\pi^k\,\mid \, \text{\rm $a_k\in K$ and $f$ is holomorphic  on $B(r,1)$} \right \}.
$$ 
%Let  $\boB_K^{\dagger}= \underset{r\geqslant p-1}\to \cup \boB^{\dag,r}_{K}$,
%$\bold A_K^{\dag}=\{f(\pi) = \sum_{k\in \Bbb Z}
%a_k\pi^k\,\mid \,f(\pi)\in \boB_K^{\dagger}\,\, \text{\rm and} \,\, a_k\in O_K\}$
Set    
$\CR (K) =\underset{r\geqslant  p-1}\to \cup \boB^{\dag,r}_{\text{rig},K}$ and $\CR^+(K)=\CR (K) \cap K[[\pi]].$ 
It is not difficult to check that these rings  are stable under $\Gamma$ and  $\Ph .$
To simplify notations we will write  $\Cal R=\CR (\Bbb Q_p)$ and $\Cal R^+= \CR^+(\Bbb Q_p).$
As usual, we set
$$
t=\log (1+\pi)=\underset{n=1}\to{\overset\infty\to\sum} (-1)^{n+1} \frac{\pi^n}{n} \in \Cal R
$$
Note that $\Ph (t)=pt$ and $\tau (t)=\chi (\gamma) t$, $\tau \in \Gamma.$
\newline
{\,}
\flushpar
{\bf 1.1.2. $(\Ph,\Gamma)$-modules} (see \cite{F2}, \cite{CC1}). 
Let $A$ be either $\boB_K^{\dag}$ or $\CR (K).$
A $(\Ph,\Gamma)$-module over A  is a finitely generated free $A$-module $\bD$
equipped with semilinear actions of $\Ph$ and  $\Gamma$ commuting to each other and such that
the induced linear map $\Ph \,:\,A\otimes_{\Ph} \bD @>>>\bD$ is an isomorphism.
Such a module    is said to be etale   if it admits a $\bold A_K^{\dag}$-lattice $N$ stable under $\Ph$ and $\Gamma$ and such that 
$\Ph \,:\,\bold A_K^{\dag}\otimes_{\Ph} N @>>>N$ is an isomorphism. The functor $\bD\mapsto \CR(K)\otimes_{\boB_K^{\dag}} \bD$
induces an equivalence between the category of etale $(\Ph,\Gamma)$-modules over $\boB_K^{\dag}$ and the category
of $(\Ph,\Gamma)$-modules over $\CR (K)$ which are of slope $0$ in the sense of Kedlaya's theory  (\cite{Ke} and
\cite{C5}, Corollary 1.5). Then Fontaine's classification of $p$-adic representations \cite{F2} together
with the main result of \cite{CC1} lead to the following statement.   

\proclaim{ Proposition 1.1.3} i) The functor 
$$
\bD^{\dagger}\,\,:\,\,V \mapsto \bD^{\dagger}(V)=(\boB^{\dagger}\otimes_{\Bbb Q_p}V)^{H_K}
$$
establishes an equivalence between the category of $p$-adic representations
of $G_K$ and the category of etale $(\Ph,\Gamma)$-modules over $\boB^{\dagger}_K .$

ii) The functor $\Ddagrig (V)=\Brigdag \otimes
_{\boB^{\dagger}_K} \bD^{\dagger}(V)$ 
gives  an equivalence between the
category of $p$-adic representations of $G_K$ and the category of
$(\Ph,\Gamma)$-modules over $\Brigdag$ of slope $0$.
\endproclaim
\demo{Proof} see \cite{C4}, Proposition 1.7.
\enddemo
{\,}

\flushpar
{\bf 1.1.4.  Cohomology of $(\Ph,\Gamma)$-modules} (see \cite{H1}, \cite{H2}, \cite{Li}).
 Fix a generator $\gamma$ of $\Gamma$. If $\bD$ is a $(\Ph,\Gamma)$-module over $A$, we denote by
$ C_{\Ph,\gamma}(\bD)$ the complex
$$
C_{\Ph,\gamma} (\bD)\,\,:\,\,0@>f>>\bD @>>> \bD\oplus \bD@>g>> \bD@>>>0
$$
where $f(x)=((\Ph-1)\,x,(\gamma -1)\,x)$ and
$g(y,z)=(\gamma-1)\,y-(\Ph-1)\,z.$ Set $H^i(D)=H^i(C_{\Ph,\gamma}(D)).$ A short
exact sequence of $(\Ph,\Gamma)$-modules
$$
0@>>>\bD'@>>>\bD@>>>\bD''@>>>0
$$
gives rise to an exact cohomology sequence:
$$
0@>>>H^0(\bD')@>>>H^0(\bD)@>>>H^0(\bD'')@>>>H^1(\bD')@>>>\cdots @>>>
H^2(\bD'')@>>>0.
$$

\proclaim{Proposition 1.1.5} Let $V$ be a $p$-adic representation of $G_K.$
Then  the complexes $\RG (K,V)$,  $C_{\Ph,\gamma}(\bD^{\dag}(V))$ and $C_{\Ph,\gamma}(\Ddagrig (V))$ are isomorphic in the derived
category of $\Bbb Q_p$-vector spaces $\Cal D(\Bbb Q_p).$

%ii) The natural map $C_{\Ph,\gamma}^{\dag}(V)@>>>C_{\Ph,\gamma}(V)$ is an isomorphism 
%in $\Cal D(\Bbb Q_p).$
\endproclaim
\demo{Proof} This is a derived version of Herr's  computation of Galois cohomology \cite{H1}. 
The  proof is given in the Appendix, Propositions A.3 and  Corollary A.4.
\enddemo 

{\,}
\flushpar
{\bf 1.1.6. Iwasawa cohomology.}  Recall that $\Lambda$ denotes the Iwasawa algebra of $\Gamma_1$, $\Delta=\Gal (K_1/K)$
and   $\Lambda (\Gamma)=\Bbb Z_p[\Delta]\otimes_{\Bbb Z_p}\Lambda$. Let $\iota \,\,:\,\,
\Lambda (\Gamma) @>>>\Lambda (\Gamma)$ denote the involution defined by $\iota (g)=g^{-1},$
$g\in \Gamma .$ 
If $T$ is a $\Bbb Z_p$-adic representation of $G_K$,  then the induced module
$\text{\rm Ind}_{K_\infty/K}(T)$ is isomorphic to $(\Lambda (\Gamma)\otimes_{\Bbb Z_p}T)^\iota$
and we set
$$
\RG_{\Iw}(K,T)\,=\,\RG (K, \text{\rm Ind}_{K_\infty/K}(T)).
$$
Write  $H^i_{\Iw}(K,T)$  for  the Iwasawa cohomology
$$
H^i_{\Iw}(K,T)=\varprojlim_{\text{cor}_{K_n/K_{n-1}}} H^i(K_n,T).
$$

Recall that there are canonical and functorial isomorphisms
 $$
 \aligned
 &\bold R^i \Gamma_{\Iw}(K,T)\,\simeq H^i_{\Iw}(K,T),\qquad i\geqslant 0,\\
 &\RG_{\Iw}(K,T)\otimes^{\bold L}_{\Lambda (\Gamma)} \Bbb Z_p[G_n] \simeq \RG (K_n,T)
 \endaligned
 $$
 (see  \cite{Ne2}, Proposition 8.4.22). The interpretation of the Iwasawa cohomology 
in terms of $(\Ph,\Gamma)$-modules was found by Fontaine (unpublished but see \cite{CC2}). 
We give here the derived version of this result. Let 
$\psi \,\,:\,\,\boB@>>>\boB$ be the operator
defined by the formula
$
\psi (x)\,=\,\frac{1}{p} \Ph^{-1} \left (\text{\rm Tr}_{\boB/\Ph (\boB)} (x)\right ).
$
We see immediately that  $\psi \circ  \Ph={\text{id}}.$ Moreover $\psi$ commutes with the action of $G_K$ 
and $\psi ({\bold A}^{\dag})={\bold A}^{\dag}$. Consider the complexes
$$
\aligned
&C_{\Iw,\psi}(T)\,\,:\,\, \bD (T)@>\psi-1>>\bD (T),\\
&C^{\dag}_{\Iw,\psi}(T)\,\,:\,\, \bD^{\dag} (T)@>\psi-1>>\bD^{\dag} (T).
\endaligned
$$

\proclaim{Proposition 1.1.7} i) The complexes $\RG_{\Iw}(K,T)$, 
$C_{\Iw,\psi}(T)$ and $C_{\Iw,\psi}^{\dag}(T)$ are naturally isomorphic in the derived category  $\Cal D(\Lambda (\Gamma))$
of $\Lambda (\Gamma)$-modules.
\endproclaim
\demo{Proof} See Proposition A.7 and Corollary A.8.
\enddemo 

{\,}
\flushpar
{\bf 1.1.8.$(\Ph,\Gamma)$-modules of rank $1$.} Recall the computation of the cohomology of
$(\Ph,\Gamma)$-modules of rank $1$ following Colmez \cite{C4}. As in {\it op. cit.},  we
consider the case  $K=\Bbb Q_p$ and 
put  $\Cal R=\boB^{\dag}_{\text{\rm rig},\Bbb Q_p}$ and $\Cal R^+=\boB^+_{\text{\rm rig},\Bbb Q_p}.$
The differential operator
$\partial =(1+\pi)\dsize\frac{d}{d\pi}$ acts on $\CR$ and $\CR^+.$
If $\delta \,:\,\Bbb Q_p^*@>>> \Bbb Q_p^*$ is a continuous character, we write  $\Cal R (\delta)$
for  the $(\Ph,\Gamma)$-module $\Cal R e_{\delta}$ defined by
$\Ph (e_{\delta})=\delta (p) e_{\delta}$ and 
$\gamma (e_{\delta})=\delta (\chi(\tau))\,e_{\delta}.$ Let $ x$ denote
the character induced by the natural inclusion  of $\Bbb Q_p$ in $L$
and $\vert x\vert$ the character defined by $\vert
x\vert=p^{-v_p(x)}.$ 
%To simplify notations, set $D_m\,=\,\Cal R(\vert x\vert x^m)=\Cal Rw_m.$

\proclaim{Proposition 1.1.9} Let $\delta \,:\,\Bbb Q_p^*@>>>\Bbb Q_p^*$
be a continuous character. Then:

i) 

 $$H^0(\Cal R(\delta))=\cases \Bbb Q_p t^m &\text{\rm if
$\delta=x^{-m}$, $k\in \Bbb N$}\\
0 & \text{\rm otherwise.}
\endcases
$$

ii) 

$$\text{\rm dim}_{\Bbb Q_p} (H^1(\Cal R(\delta)))\,=\,\cases
2 &\text{\rm if either $\delta (x)= x^{-m}$, $m \geqslant 0$ or
$\delta (x)=\vert x\vert x^{m},$ $k\geqslant 1,$}\\
 1 &\text{\rm otherwise}.
\endcases
$$

iii) Assume that $\delta (x)= x^{-m},$ $m\geqslant 0.$ The classes 
$\cl (t^m,0)\,e_{\delta}$ and $\cl (0,t^m)\,e_{\delta}$ form a basis of  $H^1(\Cal R(x^{-m}))$.

iv) Assume that $\delta (x)= \vert x\vert x^{m},$ $m\geqslant 1.$ Then 
$H^1(\CR (\vert x\vert x^m)),$ $m \geqslant 1$  is generated by $\cl (\alpha_m)$ and $\cl (\beta_m)$
where
$$
\align
&\alpha_m=\frac{(-1)^{m-1}}{(m-1)!}\, \partial^{m-1} \left (\frac{1}{\pi}+\frac{1}{2},a \right )\,e_{\delta} ,\qquad
(1-\Ph)\,a=(1-\chi (\gamma) \gamma)\,\left (\frac{1}{\pi}+\frac{1}{2} \right ),\\
&\beta_m=\frac{(-1)^{m-1}}{(m-1)!}\,\partial^{m-1} \left (b,\frac{1}{\pi} \right )\,e_{\delta}, \qquad
(1-\Ph)\,\left (\frac{1}{\pi}\right )\,=\,(1-\chi (\gamma)\,\gamma)\,b
\endalign
$$
\endproclaim
\demo{Proof} See  \cite{C4}, sections 2.3-2.5.
\enddemo
\flushpar
\newline
\newline
{\bf 1.2. Crystalline representations.}
\newline
{\bf 1.2.1. The rings $\Bc$ and $\Bd$} (see \cite{F1}, \cite{F4}).
Let $\theta_0\,:\,\A^+@>>>O_C$ be the map given by the formula
$$
\theta_0 \left (\underset{n=0}\to{\overset{\infty}\to \sum} [u_n]p^n \right )\,=\,
\underset{n=0}\to{\overset{\infty}\to \sum} u_n^{(0)} p^n.
$$
It can be shown that $\theta_0$ is a surjective ring homomorphism and that $\ker (\theta_0)$ 
is the principal ideal generated by $\omega \,=\,\underset{i=0}\to {\overset{p-1} \to \sum} [\epsilon]^{i/p}.$
By linearity, $\theta_0$ can be extended to a map $\theta\,:\, \tilde {\boB}^+@>>>C$. The ring $\Bd^+$
is defined to be the completion of $\tilde {\boB}^+$ for the $\ker (\theta)$-adic topology:
$$
\Bd^+\,=\,\varprojlim_n \tilde {\boB}^+/ker (\theta)^n.
$$
This is a complete discrete valuation ring with residue field $C$ equipped with a natural
action of $G_K.$ Moreover, there exists a canonical embedding $\bar K  \subset \Bd^+.$
The series $t=\underset{n=0}\to{\overset \infty \to \sum}(-1)^{n-1}{\pi}^n/n$ converges
in the topology of $\Bd^+$ and it is easy to see that $t$ generates the maximal ideal of $\Bd^+.$ 
The Galois group acts on $t$ by the formula $g (t)=\chi (g)t.$  Let $\Bd=\Bd^+[t^{-1}]$ be the
field of fractions of $\Bd^+.$ This is a complete discrete valuation field equipped with 
a $G_K$-action and an exhaustive separated decreasing filtration $\F^i \Bd=t^i\Bd^+.$ 
As $G_K$-module, $\F^i\Bd/\F^{i+1}\Bd \simeq C(i)$ and  $\Bd^{G_K}=K.$

Consider the $PD$-envelope of $\A^+$ with a respect to the map $\theta_0$
$$
\A^{\text{\rm PD}}=\A^+\left [\frac{\omega^2}{2!},\frac{\omega^3}{3!},\ldots ,\frac{\omega^n}{n!},\ldots \right ]
$$
and denote by ${\A}^+_{\text{\rm cris}}$ its $p$-adic completion. Let $\Bc^+={\A}^+_{\text{\rm cris}}\otimes_{\Bbb Z_p}\Bbb Q_p$
and $\Bc=\Bc^+[t^{-1}].$ Then $\Bc$ is a subring of $\Bd$ endowed with  the induced filtration and  Galois action.
Moreover, it is equipped with a continuous Frobenius $\Ph$, extending the map $\Ph \,:\,\A^+@>>>\A^+.$
One has $\Ph (t)=p\,t.$ 
\newline
\, 
\flushpar
{\bf 1.2.2. Crystalline representations} (see \cite{F5}, \cite{Ber1}, \cite{Ber2}).
\newline
Let $L$ be a finite extension of $\Bbb Q_p$. Denote by $K$ its maximal
unramified subextension.  A filtered Dieudonn\'e module over $L$
is a finite dimensional $K$- vector space $M$ equipped with the
following structures:

$\bullet$ a $\sigma$-semilinear bijective map $\Ph \,:\,M@>>>M;$

$\bullet$ an exhaustive decreasing filtration $(\F^i M_L)$ on
the $L$-vector space $M_L=L\otimes_{K} M.$

A $K$-linear map $f\,:\,M@>>>M'$ is said to be a morphism of
filtered modules if

$\bullet$ $f(\Ph (d))=\Ph (f(d)),\qquad \text{\rm for all $d\in M;$}$

$\bullet$ $f(\F^i M_L)\subset \F^i M'_L, \qquad \text{\rm for
all $i\in \Bbb Z.$}$

The category  $\bold M\bold F^{\Ph}_L$ of filtered Dieudonn\'e modules is additive,
has kernels and cokernels but is not abelian. 
%It can be equipped
%with a tensor product defined by
%$$
%\aligned
%& \F ^i(D'\otimes D'')_K=\underset{j+k=i}\to\sum \F^j
%D'_K\otimes \F^k D''_K;\\
%& \Ph (d'\otimes d'')=\Ph (d')\otimes \Ph (d'').
%\endaligned
%$$
Denote by $\bold 1$  the vector space
$K_0$ with the natural action of $\sigma$ and  the
filtration given by
$$
\F^i\bold 1=\cases K,&\text{\rm if $i\leqslant 0,$}\\
0, &\text{\rm if $i>0.$}
\endcases
$$
Then $\bold 1$ is a unit object of $\bold M\bold F^{\Ph}_L$ i.e.
$
M\otimes \bold 1\simeq \bold 1\otimes M\simeq M
$
 for
any $M$.

If $M$ is a one dimensional Dieudonn\'e module and $d$ is a basis vector
of $M$, then $\Ph (d)=\alpha v$ for some $\alpha \in K$. Set
$t_N(M)= v_p(\alpha)$ and denote by $t_H(M)$ the unique filtration jump  
of $M.$ If $M$ is of an arbitrary finite dimension $d$, set 
$t_N(M)= t_N(\overset d\to \wedge M)$ and  $t_H(M)= t_H(\overset d\to \wedge M).$   
A Dieudonn\'e module $M$ is said to be weakly admissible if $t_H(M)=t_N(M)$
and if $t_H(M')\leqslant t_N(M')$  for any $\Ph$-submodule $M'\subset M$
equipped with the induced filtration. Weakly admissible modules 
form a subcategory of $\bold M\bold F_L$ which we denote by $\bold M\bold F_L^{\Ph,f}.$

If $V$ is a $p$-adic representation of $G_L$, define
$
\Dd (V)=(\Bd \otimes V)^{G_L}.
$
Then $\Dd (V)$ is a $L$-vector space  equipped with the  decreasing filtration
$\F^i\Dd (V)\,=\,(\F^i\Bd\otimes V)^{G_L}$.  One has 
$\dim_L \Dd (V)\leqslant \dim_{\Bbb Q_p}(V)$ and 
$V$ is said to be de
Rham if
$
\dim_L \Dd (V)= \dim_{\Bbb Q_p}(V).
$
Analogously one defines
$
\Dc (V)=(\Bc \otimes V)^{G_L}.
$
Then $\Dc (V)$ is a filtered Dieudonn\'e module over $L$ of
dimension $\dim_{K} \Dc(V)\leqslant \dim_{\Bbb Q_p}(V)$ and $V$ is
said to be crystalline if  the equality holds here.
In particular, for crystalline representations one has 
$
\Dd (V)= \Dc(V)\otimes_{K}L.
$
By the theorem of Colmez-Fontaine \cite{CF},  the functor $\Dc $
establishes an equivalence between the category of crystalline
representations of $G_L$ and $\bold M\bold F_L^{\Ph,f}$.
Its quasi-inverse $\bold V_{\text{\rm cris}}$ is given by
$
\bold V_{\text{\rm cris}}(D)=\F^0(D\otimes_{K_0} \Bc)^{\Ph=1}.
$

An important result of Berger (\cite{Ber 1}, Theorem 0.2)
says that
$\Dc (V)$ can be recovered from the $(\Ph,\Gamma)$-module $\Ddagrig (V).$
The situation is particularly simple  if 
If $L/\Bbb Q_p$ is unramified. In this case 
set 
${\bold D}^+(V)=(V\otimes_{\Bbb Q_p}\boB^+)^{H_K}$ and 
${\bold D}_{\text{\rm rig}}^+(V)= \CR^+(K)\otimes_{\boB_K^+}{\bold D}^+(V).$ Then
$$
\Dc (V)\,=\,\left ( {\bold D}_{\text{\rm rig}}^+(V) \left [\frac{1}{t}\right ]\right )^{\Gamma}
$$
(see \cite{Ber2}, Proposition 3.4).

\head {\bf \S 2. The  exponential map}
\endhead

\flushpar
{\bf 2.1. The Bloch-Kato exponential map} (\cite{BK}, \cite{Ne1}, \cite{FP}).
\newline
{\bf 2.1.1. Cohomology of Dieudonn\'e modules.}  Let $L$ be a finite extension of $\Bbb Q_p.$ Recall that we denote by $\bold M\bold F^{\Ph}_L$ 
the category of filtered Dieudonn\'e modules over $L.$ If $M$ is an object of $\bold M\bold F^{\Ph}_L$, define
$$
H^i (L,M)=\Ext^i_{\bold M\bold F^{\Ph}_L}(\bold 1, M), \qquad i=0,1.
$$
Remark that   $H^*(L,M)$ can be computed explicitly  as
the cohomology of the complex
$$
C^{\bullet}(M)\,\,:\,\,M@>f>>({M_L}/{\F^0 M_L}) \oplus M
$$
where the modules are placed in degrees $0$ and $1$ and $f(d)=
(d\pmod {\F^0 M_L},(1-\Ph)\,(d))$ (\cite{Ne1},\cite{FP}). Remark that if $M$ is weakly admissible
then each extension $0@>>>M@>>>M'@>>>\bold 1@>>>0$ is weakly admissible too and we can write
$ H^i (L,M)=\Ext^i_{\bold M\bold F^{\Ph,f}_L}(\bold 1, M).$ 
\newline
\,
\flushpar
{\bf 2.1.2. The exponential map.}  Let $\Rep_{\text{\rm cris}}(G_K)$ denote the category of crystalline representations
of $G_K.$ For any object $V$ of $\Rep_{\text{\rm cris}}(G_K)$ define
$$
H_f^i(K,V)=\Ext_{\Rep_{\text{\rm cris}}(G_K)}^i(\Bbb Q_p(0),V).
$$
An easy computation shows that
$$
H^i_f(K,V)=\cases H^0(K,V), &\text{if $i=0$,}\\
\ker\,(H^1(K,V)@>>>H^1(K,V\otimes \Bc)), &\text{if $i=1$,}\\
0, &\text{if $i\geqslant 2.$}
\endcases
$$
Let $t_V(K)=\Dd (V)/\F^0 \Dd (V)$ denote the tangent space of $V$. 
The rings $\Bd$ and $\Bc$ are related to each other via the fundamental
exact sequence
$$
0@>>>\Bbb Q_p@>>>\Bc @>f>>\Bd/\F^0\Bd \oplus \Bc @>>>0
$$
where $f(x)=( x\pmod{\F^0\Bd}, (1-\Ph)\,x)$ (see \cite{BK}, \S4).
Tensoring this sequence with $V$ and taking cohomology one obtains an exact sequence
$$
0@>>>H^0(K,V)@>>>\Dc (V)@>>>t_V(K)\oplus \Dc (V) @>>>H^1_f(K,V)@>>>0.
$$
The last map of this sequence gives rise to the  Bloch-Kato exponential map
$$
\exp_{V,K}\,\,:\,\, t_V(K)\oplus \Dc (V)@>>>H^1(K,V).
$$
Following \cite{F3}  set
$$
\RG_f(K,V)\,=\,C^\bullet (\Dc (V))\,=\,\left [ \Dc (V)@>f>> t_V(K)\oplus
\Dc (V) \right ].
$$

From the classification of crystalline representations in terms
of Dieudonn\'e modules  it follows that  the functor $\Vc$ induces natural isomorphisms
$$
r_{V,p}^i\,\,:\,\,\bold R^i\Gamma_f (K,V)@>>> H^i_f(K,V),\qquad i=0,1.
$$
The composite homomorphism
$$
t_K(V)\oplus \Dc (V)@>>> \bold R^1\Gamma_f (K,V)@>r_{V,p}^1>> H^1(K,V)
$$
 coincides with the Bloch-Kato exponential map
$\exp_{V,K}$ (\cite{Ne1}, Proposition 1.21).
\newline
\,
\flushpar
{\bf 2.1.3. The map $\RG_f (K,V)@>>>\RG (K,V)$.} 
Let $g\,:\,B^{\bullet}@>>>C^{\bullet}$ be a morphism of complexes.
We denote by $\Tot^{\bullet} (g)$ the complex
$\Tot^n(g)=C^{n-1}\oplus B^{n}$ with differentials
$d^n\,:\,\Tot^n(g)@>>>\Tot^{n+1}(g)$ defined by the formula
$
d^n(c,b)=((-1)^n g^n(b)+d^{n-1}(c),d^n(b)).
$
It is well known that if
$
0@>>>A^{\bullet}@>f>>B^{\bullet}@>g>>C^{\bullet}@>>>0
$
is an exact sequence of complexes, then $f$
induces a quasi isomorphism
$
A^\bullet \overset{\sim}\to \rightarrow \Tot^{\bullet}(g).
$
In particular, tensoring the fundamental exact sequence with $V$, we obtain an
exact sequence of complexes
$$
0@>>>\RG (K,V)@>>>C_c^{\bullet} (G_K,V\otimes
\Bc)@>f>>C_c^{\bullet} (G_K,(V\otimes (\Bd /\F^0 \Bd)) \oplus
(V\otimes \Bc))@>>>0 
$$
which gives a quasi isomorphism  $ \RG (K,V) \overset{\sim}\to \rightarrow \Tot^{\bullet}
(f). $
Since $\RG_f(K,V)$ coincides tautologically with the complex
$$ C_c^0 (G_K,V\otimes \Bc)@>f>>C_c^0 (G_K,(V\otimes (\Bd /F^0 \Bd)) \oplus
(V\otimes \Bc))
$$
we obtain a diagram
$$
\xymatrix{ \RG (K,V) \ar[r]^{\sim} &\Tot^{\bullet} (f)   \\
&  \RG_f(K,V)
\ar[u]
\ar
%@<0ex> 
@{.>}[ul]}
$$
 which defines a morphism 
$
\RG_f (K,V)@>>>\RG (K,V)
$
in $\Cal D(\Bbb Q_p)$ (see \cite{BF}, Proposition 1.17). Remark that the induced homomorphisms $\R^i\Gamma_f(K,V)@>>>
H^i(K,V)$ ($i=0,1$) coincide with the composition of $r_{V,p}^i$ 
with natural embeddings $H^i_f(K,V)@>>> H^i(K,V).$
\newline
\,
\flushpar
{\bf 2.1.4. Exponential map for $(\Ph,\Gamma)$-modules.} In this subsection we define an analog of the exponential
map for crystalline $(\Ph,\Gamma)$-modules. See  \cite{Na} for a more general
setting.
Let  $K/\Bbb Q_p$ be an unramified extension. 
If $\bD$ is a $(\Ph,\Gamma)$-module over $\CR (K)$
define
$$
\Cal D_{\text{\rm cris}}(\bD)\,=\,\left (\bD\left [{1}/{t} \right ] \right )^{\Gamma}.
%\qquad \Cal D_{\text{\rm st}}(\bD)\,=\,\left (\bD\otimes_{\CR (K)} \CR_{\log}(K)[1/t] \right )^{\Gamma}
$$
It can be shown that   $\Cal D_{\text{\rm cris}}(\bD)$ is a finite dimensional $K$-vector space
equipped with a natural decreasing filtration $\F^i\Cal D_{\text{\rm cris}}(\bD)$ and
a semilinear action of $\Ph$. One says  that $D$ is crystalline if
$$
\dim_K(\Cal D_{\text{\rm cris}}(\bD))\,=\,\text{\rm rg}  (\bD).
$$
From \cite{Ber4}, Th\'eor\`eme A it follows that the functor  $\bD\mapsto \Cal D_{\text{\rm cris}}(\bD)$ is an equivalence
between the category of crystalline $(\Ph,\Gamma)$-modules and $\bold M\bold F^{\Ph}_K.$
Remark that if $V$ is a $p$-adic representation of $G_K$ then 
$
\Dc (V)=\Cal D_{\text{\rm cris}}(\Ddagrig (V))
$ 
and $V$ is crystalline if and only if $\Ddagrig (V)$ is.

Let $\bD$ be a $(\Ph,\Gamma)$-module. To any cocycle $\alpha =(a,b)\in Z^1(C_{\Ph,\gamma}(\bD))$ 
one can associate the extension
$$
0@>>>\bD@>>>\bD_{\alpha}@>>> \CR(K)@>>>0
$$
defined by 
$$
\bD_{\alpha}=D\oplus \CR (K)\,e,\qquad (\Ph-1)\,e=a, \quad (\gamma-1)\,e=b.
$$
As usual, this gives rise to an isomorphism  $H^1(\bD)\simeq \text{\rm Ext}^1_{\CR} (\CR (K),\bD).$ 
We say that $\cl (\alpha)$ is crystalline if
$\dim_{K} \left (\Cal D_{\text{\rm cris}}(\bD_{\alpha})\right )=
\dim_{K} \left (\Cal D_{\text{\rm cris}}(\bD)\right )+1
$
and define 
$$
H^1_f(\bD)\,=\,\{ \cl (\alpha)\in H^1(\bD)\,\,\vert \,\, \text{\rm $\cl (\alpha)$ is crystalline}\,\}
$$
(see \cite{Ben2}, section 1.4.1). If $\bD$ is crystalline (or more generally potentially
semistable ) one has a natural isomorphism
$$
H^1(K, \Cal D_{\text{\rm cris}}(\bD))@>>>H^1_f(\bD).
$$
Set $t_{\bD}= \Cal D_{\text{\rm cris}}(\bD)/\F^0 \Cal D_{\text{\rm cris}}(\bD)$ and denote by
$
\exp_{\bD}\,\,:\,\,t_{\bD}\oplus \Cal D_{\text{\rm cris}}(\bD)@>>>H^1(\bD)
$
the composition of this isomorphism with the projection 
$t_{\bD}\oplus \Cal D_{\text{\rm cris}}(\bD) @>>>H^1(K, \Cal D_{\text{\rm cris}}(\bD))$ and the embedding
$H^1_f(\bD)\hookrightarrow H^1(\bD).$
\newline
\,

\flushpar 
Assume  that $K=\Bbb Q_p.$ To simplify notation we will write $\bD_m$ for $\CR (\vert x\vert x^m)$ and $e_m$ for its canonical basis.
Then  $\Cal D_{\text{\rm cris}}(\bD_m)$ is the one dimensional $\Bbb Q_p$-vector space 
generated by $t^{-m} e_m$.   
As in \cite{Ben2}, we normalise the basis $(\cl (\alpha_m), \cl (\beta_m))$ of $H^1(\bD_m)$
putting $\alpha_m^*=\left (1-1/p\right )\,\cl (\alpha_m)$ and 
$\beta_m^*=\left (1-1/p \right ) \log (\chi (\gamma))\, \cl (\beta_m).$

\proclaim{Proposition 2.1.5} i) $H^1_f(\bD_m)$ is the one-dimensional $\Bbb Q_p$-vector space generated by 
$\alpha^*_m$. 
 
ii) The exponential map
$$
\exp_{\bD_m}\,\,:\,\,t_{\bD_m} @>>>H^1(\bD_m)
$$
sends $t^{-m}w_m$ to 
$- \,\alpha^*_m.
$
\endproclaim
\demo{Proof} This is a reformulation of \cite{Ben2}, Proposition 1.5.8 ii). 
\enddemo
$\,$

\flushpar
{\bf 2.2. The large exponential map.}
\newline
{\bf 2.2.1. Notation.} In this section $p$ is an odd prime number, $K$ is a finite unramified extension
of $\Bbb Q_p$ and $\sigma$ the absolute Frobenius acting on $K.$ 
Recall that $K_n=K(\zn)$ and $\dsize
K_\infty =\cup_{n=1}^{\infty} K_n.$ We set  $\Gamma =\Gal (K_\infty
/K),$  $\Gamma_n=\Gal (K_\infty/K_n)$ and $\Delta=\Gal (K_1/K)$ . 
Let  $\La=\Bbb Z_p[[\Gamma_1]]$ and $\Lambda (\Gamma)=\Bbb Z_p[\Delta]\otimes_{\Bbb Z_p}\Lambda.$
We will consider the following operators acting on the ring $K[[X]]$ of formal power series with coefficients in $K$:  

$\bullet$ The ring homomorphism $\sigma \,:\, K[[X]]@>>>K[[X]]$ defined
by
$
\sigma \left (\dsize \sum_{i=0}^\infty  a_iX^i\right )=\dsize\sum_{i=0}^\infty
\sigma (a_i)X^i;
$

$\bullet$  The ring homomorphism $\Ph \,:\,K[[X]]@>>>K[[X]]$ defined by
$$
\Ph \left (\sum_{i=0}^\infty  a_iX^i\right )=\sum_{i=0}^\infty
\sigma (a_i)\Ph(X)^i, \qquad \Ph (X)=(1+X)^p-1.
$$

$\bullet$ The differential operator $\dop=(1+X)\dsize\frac{d}{dX}.$ One
has
$
\dop\circ  \Ph=p\Ph \circ \dop.
$

$\bullet$ The operator $\psi \,: \,K[[X]]@>>> K[[X]]$ defined by
$
\psi (f(X))= \dsize\frac{1}{p} \Ph^{-1} \left
(\dsize\sum_{\z^p=1}f((1+X)\z-1)\right ).
$
It is easy to see that $\psi $ is a left inverse to $\Ph,$ i.e.
that
$
\psi \circ \Ph=\text{\rm id}.
$
%Let $K((X))$ be the field of fractions of $K[[X]].$ Then
%$K((X))/K((X))^{\Ph}$ is an algebraic extension of degree $p$ and
%one has a well defined  trace map $\Tr_{\Ph}\,:\, K((X))@>>>
%K((X))^{\Ph}.$ It it easy to check that
%$$
%\psi=\frac{1}{p}\Ph^{-1}\circ \Tr_{\Ph}.
%$$

$\bullet$  An action of $\Gamma$ given by
$
\g \left (\dsize\sum_{i=0}^\infty  a_iX^i\right )=\dsize \sum_{i=0}^\infty a_i
\g( X)^i, \qquad \g (X)=(1+X)^{\chi (\g)}-1.
$

Remark that these formulas are compatible with the definitions from sections 1.1.1 and 1.1.6.
Fix a generator $\g_1\in \Gamma_1$ and define 
$$
\Cal H\,=\,\{f(\g_1-1)\,|\,f \in \Bbb Q_p[[X]] \,\,
\text{{\rm is holomorphic on $B(0,1)$}}\},\qquad
\Cal H(\Gamma)\,=\,\Bbb Z_p[\Delta]\otimes_{\Bbb Z_p}\Cal H.
$$
{\,}

\flushpar
{\bf 2.2.2. The map $\boldsymbol \Xi_{V,n}^{\ep}$. } It is well known that  $ \Bbb Z_p[[X]]^{\psi=0}$ is a free $\Lambda$-module generated by
$(1+X)$  and the operator $\dop$ is bijective on $\Bbb
Z_p[[X]]^{\psi=0}.$ If $V$ is  a crystalline representation of $G_K$ 
put
$
\Cal D(V)=\Dc (V)\otimes _{\Bbb Z_p} \Bbb Z_p[[X]]^{\psi=0} .
$
Let 
$\boldsymbol \Xi_{V,n}^{\ep}\,:\,\Cal D(V)_{\Gamma_n} [-1]@>>>\RG_f
(K_n,V)
$
be the map defined by
$$
 \boldsymbol{\Xi}_{V,n}^\ep (\alpha)\,=\,\cases p^{-n} (\sum_{k=1}^n
(\sigma \otimes \Ph)^{-k} \alpha (\zeta_{p^k}-1),\,-\alpha (0))
&\text{\rm if
$n\geqslant 1$,} \\
\text{\rm Tr}_{K_1/K}\,\left ( \boldsymbol{\Xi}_{V,1}^{\ep} (\alpha)\right )
&\text{\rm if $n=0.$}
\endcases
$$
An easy computation shows that $ \boldsymbol{\Xi}_{V,0}^\ep \,:\,
\Dc (V) [-1] @>>> \RG_f(K,V)$ is given by the formula
$$
\boldsymbol{\Xi}_{V,0}^\ep
(a)\,=\,\frac{1}{p}\,(-\Ph^{-1}(a),-(p-1)\,a).
$$
In particular, it is homotopic to the map
$
a\mapsto -(0,(1-p^{-1}\Ph^{-1} )\,a).
$
Write
$$
\Xi_{V,n}^\ep \,:\,\Cal D(V)@>>> \bold R^1 \Gamma
(K_n,V)=\frac{t_V(K_n)\oplus \Dc (V)}{\Dc (V)/V^{G_K}}
$$
denote the homomorphism  induced by $\boldsymbol \Xi_{V,n}^\ep.$ Then 
$$
\Xi_{V,0}^\ep (a)=-(0,(1-p^{-1}\Ph^{-1} )\,a) \pmod{\Dc (V)/V^{G_K}}.
$$
If $\Dc (V)^{\Ph=1}=0$ the operator $1-\Ph$ is invertible on $\Dc (V)$ and we can write  
$$
\Xi_{V,0}^\ep (a)= \left (\frac{1-p^{-1}\Ph^{-1}}{1-\Ph} \,a,0 \right ) \pmod{\Dc (V)/V^{G_K}}. \tag{2.1}
$$

\flushpar
For any $i\in \Bbb Z$ let
$\Delta_i\,:\, \Cal D(V)@>>>\dsize\frac{\Dc(V)}{(1-p^i\Ph)\Dc
(V)}\,\otimes \Bbb Q_p(i)
$
be the map given by 
$$\Delta_i(\alpha (X))= \dop^i\alpha(0)\otimes \ep^{\otimes i}\pmod
{(1-p^i\Ph)\Dc(V)}.$$
 Set $\Delta=\oplus_{i\in \Bbb Z}\Delta_i .$
  If $\alpha\in \Cal
D(V)^{\Delta =0},$ then by  \cite{PR1}, Proposition 2.2.1 there exists $F\in \Dc(V)\otimes_{\Bbb Q_p}\Bbb Q_p[[X]]$ 
which converges on the open unit disk 
and such
that $ (1-\Ph)F=\alpha.$ A short computation shows that
$$
\Xi_{V,n}^\ep (\alpha) = p^{-n} ( (\sigma\otimes \Ph)^{-n}
(F)(\zn-1),0)\pmod {\Dc (V)/V^{G_K}},\qquad \text{\rm
if $n\geqslant 1$}
$$
(see \cite{BB}, Lemme 4.9). 
\newline
\newline
{\bf 2.2.3. Construction of the large exponential map.} As $\Bbb Z_p[[X]] \left [1/p\right
] $ is a principal ideal domain and $\Cal H$ is $\Bbb Z_p[[X]]
\left [1/p\right ]$-torsion free,  $\Cal
H$ is flat. Thus
$$
C_{\Iw, \psi}^\dag(V)\otimes^{\bold L}_{\Lambda_{\Bbb
Q_p}}\HG\,=\,C_{\Iw, \psi}^\dag(V)\otimes_{\Lambda_{\Bbb Q_p}
}\HG\,=\,\left [ \HG\otimes_{\Lambda_{\Bbb Q_p}} \bD^\dag
(V)@>\psi-1>>\HG\otimes_{\Lambda_{\Bbb Q_p}}\bD^\dag (V) \right ].
$$
By proposition 1.1.7 on has an isomorphism in $\Cal D (\HG)$
$$
\RG_{\Iw}(K,V)\otimes^{\bold L}_{\Lambda_{\Bbb Q_p}}\HG \simeq
C_{\Iw, \psi}^\dag (V)\otimes_{\Lambda_{\Bbb Q_p}}\HG .
$$
The action of $\HG$ on
$\bD^{\dag}(V)^{\psi=1}$ induces an injection
$
\HG \otimes_{\La_{\Bbb Q_p}} \bD^{\dag}(V)^{\psi=1}
\hookrightarrow \Ddagrig (V)^{\psi=1}.
$
Composing this map with the canonical isomorphism
$\Hi^1(K,V)\simeq \bD^{\dag}(V)^{\psi=1}$ we obtain a map
$
\HG \otimes_{\La_{\Bbb Q_p}} \Hi^1(K,V) \hookrightarrow \Ddagrig
(V)^{\psi=1}.
$
For any $k\in \Bbb Z$ set $\nabla_k=t\dop-k=
t\dsize\frac{d}{dt}-k.$ An easy induction shows that
$
\nabla_{k-1}\circ \nabla_{k-2}\circ \cdots \circ \nabla_0\,=\,t^k\dop^k.
$

Fix $h\geqslant 1$  such that $\F^{-h}\Dc (V)=\Dc (V)$ and $V(-h)^{G_K}=0.$
For any  $\alpha \in \Cal D(V)^{\Delta=0}$ define 
$$
\Omega_{V,h}^\ep(\alpha) \,=\,(-1)^{h-1}\frac{\log \chi (\gamma_1)}{p}\,\nabla_{h-1}\circ
\nabla_{h-2}\circ \cdots \nabla_0 (F (\pi)),
$$
where $F\in \Cal H (V)$ is such that $(1-\Ph)\,F= \alpha.$ 
It is easy to see that $\Omega_{V,h}^\ep(\alpha) \in
\Drig^+(V)^{\psi=1}.$ In \cite{Ber3} Berger shows that
$\Omega_{V,h}^\ep(\alpha)\in \HG \otimes_{\La_{\Bbb
Q_p}}\bD^{\dag} (V)^{\psi=1}$ and therefore  gives rise to a map
$$
\bExp_{V,h}^\ep\,:\,\Cal D(V)^{\Delta=0}[-1] @>>>
\RG_{\Iw}(K,V)\otimes^{\bold L}_{\La_{\Bbb Q_p}}\HG
$$
Let
$$
\Exp_{V,h}^\ep\,:\,\Cal D(V)^{\Delta=0}@>>>\HG\otimes_{\La_{\Bbb
Q_p}}H_{\Iw}^1(K,V)
$$
denote the map induced by $\bExp_{V,h}^\ep$ in degree $1$. The
following theorem is a reformulation of the construction of the large exponential
map given by Berger in \cite{Ber3}.

\proclaim{Theorem 2.2.4} Let
$$
\bExp_{V,h,n}^\ep \,\,:\,\,\Cal D(V)^{\Delta=0}_{\Gam_n}[-1] @>>>
\RG_{\Iw} (K,V)\otimes^{\bold L}_{\La_{\Bbb Q_p}}\Bbb Q_p[G_n].
$$
denote the map induced by $\bExp_{V,h}^\ep.$ Then for any
$n\geqslant 0$ the following diagram in $\Cal D (\Bbb Q_p[G_n])$
is commutative:

$$
\xymatrix{  \Cal D(V)^{\Delta=0}_{\Gamma_n}[-1] \ar[rr]^
{\bExp^{\ep}_{V,h,n}} \ar[d]^{\boldsymbol \Xi^{\ep}_{V,n}} & &
\RG_{\Iw} (K,V)\otimes^{\bold L}_{\La_{\Bbb Q_p}}\Bbb Q_p[G_n] \ar
[d]^{\simeq}
\\
\bold R\Gamma_f (K_n,V) \ar[rr]^{(h-1)!} & &\RG (K_n,V)\,. }
$$

In particular,  $\Exp_{V,h}^\ep$ coincides with the large
exponential map of Perrin-Riou.
\endproclaim
\demo{Proof} Passing to cohomology in the previous diagram one
obtains  the  diagram
$$
\CD \Cal D(V)^{\Delta=0} @> \Exp^{\ep}_{V,h}
 >> \Cal H (\Gam)\otimes_{\La_{\Bbb Q_p}}H^1_{\text{\rm Iw}}(K,V)
\\
@V\varXi^{\ep}_{V,n}VV         @VV\text {\rm pr}_{V,n}V
\\
\bD_{\text{dR}/K_n}(V)\oplus \Dc (V) @>(h-1)!\,\exp_{V,K_n}>>
H^1(K_n, V)
\endCD
$$
which is exactly the definition of the large exponential map. Its
commutativity is proved in \cite{Ber3}, Theorem II.13. Now, the
theorem is an immediate consequence of the following remark. Let
$D$ be a free $A$-module and let $f_1,f_2\,:\,D[-1] @>>>
K^\bullet$ be two maps from $D[-1]$ to a complex of $A$-modules
such that the induced maps $h_1(f_1)$ and $h(f_2)\,:\,
D@>>>H^1(K^\bullet)$ coincide. Then $f_1$ and $f_2$ are homotopic.
\enddemo
{\,}
\flushpar
{\bf Remark.} The large exponential map was first constructed in \cite{PR1}.
See \cite{C1} and \cite{Ben1} for alternative constructions and 
\cite{PR4},  \cite{Na1} and \cite{Ri} for generalizations.

\head {\bf \S3. The $\scr L$-invariant}
\endhead

\flushpar
{\bf 3.1.Definition of the $\scr L$-invariant.} 
\newline
{\bf 3.1.1. Preliminaries.} Let $S$ be a finite set of primes  of $\Bbb Q$
containing $p$  and  $G_S$  the Galois group of the maximal
algebraic extension of $\Bbb Q$ unramified outside $S\cup\{\infty\}$. For each place
$v$ we denote by $G_v$  the decomposition at $v$ group and
by $I_v$ and $f_v$ the inertia subgroup and Frobenius automorphism
respectively.
Let $V$ be a pseudo-geometric $p$-adic  representation of $G_S$.
This means that the restriction of $V$ on the decomposition group
at $p$ is a de Rham representation. Following Greenberg,
for any $v\notin
\{p,\infty\}$  set
 $$
 \RG_f (\Bbb Q_v,V)\,=\,\left [ V^{I_v} @>1-f_v>> V^{I_v} \right
 ],
 $$
where the terms are placed in degrees $0$ and $1$ (see \cite{F3}, \cite{BF}). 
Note that
there is a natural quasi-isomorphism $\RG_f(\Bbb Q_v,V)\simeq
C_c^\bullet (G_v/I_v, V^{I_v}).$ Note that $\bold R^0\Gamma
(\Bbb Q_v,V)=H^0(\Bbb Q_v,V)$ and $\bold R^1\Gamma_f (\Bbb
Q_v,V)=H^1_f(\Bbb Q_v,V)$ where
$$
H^1_f(\Bbb Q_v,V)=\ker (H^1(\Bbb Q_v,V)@>>>H^1(\Bbb Q_v^{\text{\rm
ur}},V)).
$$
For $v=p$ the complex $\RG_f (\Bbb Q_v,V)$ was defined in \S2.
To simplify notation write $H^i_S(V)=H^i(G_S,V)$ for the continuous Galois cohomology
of $G_S$ with coefficients in $V$. The Bloch-Kato's Selmer
group of $V$ is defined as
$$
H^1_f(V)\,=\,\ker \left ( H^1_S(V)@>>>\bigoplus_{v\in S}
\frac{H^1(\Bbb Q_v,V)}{H^1_f(\Bbb Q_v,V)} \right ).
$$
We also set 
$$
H^1_{f,\{p\}}(V)\,=\,\ker \left ( H^1_S(V)@>>>\bigoplus_{v\in S-\{p\}}
\frac{H^1(\Bbb Q_v,V)}{H^1_f(\Bbb Q_v,V)} \right ).
$$
From the Poitou-Tate exact sequence one obtains the following 
exact sequence relating these groups (see for example \cite{PR2}, Lemme 3.3.6)
$$
0@>>>H^1_f(V)@>>>H^1_{f,\{p\}}(V)@>>>\frac{H^1(\Bbb Q_p,V)}{H^1_f(\Bbb Q_p,V)}@>>>H^1_f(V^*(1))\,\,.
$$
We also have the following formula relating dimensions of Selmer groups (see \cite{FP}, II, 2.2.2)
$$
\multline
\dim_{\Bbb Q_p} H^1_f(V) - \dim_{\Bbb Q_p} H^1_f(V^*(1)) -\dim_{\Bbb Q_p} H^0_S(V) + \dim_{\Bbb Q_p} H^0_S(V^*(1)) =\\
\dim_{\Bbb Q_p} t_V(\Bbb Q_p)-\dim_{\Bbb Q_p} H^0(\Bbb R,V).
\endmultline
$$
\flushpar
Set 
$
d_{\pm}(V)=\dim_{\Bbb Q_p}(V^{c=\pm 1}),
$ 
where $c$ denotes the complex conjugation.
\newline
\newline
{\bf 3.1.2. Basic assumptions.} Assume that $V$ satisfies the following conditions

{\bf C1)} $H^1_f(V^*(1))=0$.

{\bf C2)} $H^0_S(V)=H^0_S(V^*(1))=0$.

{\bf C3)} $V$ is crystalline at $p$ and $\varphi \,:\,\Dc (V)@>>>\Dc (V)$ is semisimple
at $1$ and $p^{-1}.$

{\bf C4)} $\Dc (V)^{\varphi=1}=0.$

{\bf C5)} The localisation map
$$
\loc \,\,:\,\ \,:\,H^1_f(V)@>>> H^1_f(\Bbb Q_p,V)
$$
is injective.
\newline
\,

These conditions  appear naturally in the following situation. Let $X$ be a proper smooth variety over $\Bbb Q$. Let $H^i_p(X)$ denote the $p$-adic etale cohomology of $X$.
Consider the Galois representations $V=H^i_p(X) (m).$ 
By Poincar\'e duality together with the hard Lefschetz theorem we have
$$
H^i_p(X)^*\simeq H^i_p(X)\,(i)
$$
and thus $V^*(1) \simeq V (i+1-2m).$ 
The Beilinson conjecture (in the formulation
of Bloch and Kato)  predict that
$$ 
H^1_f(V^*(1))=0 \qquad \text{\rm if} \qquad  w\leqslant -2. 
$$
This corresponds to the hope that there
are no nontrivial extensions of $\Bbb Q(0)$ by motives of weight $\geqslant 0.$
If $X$ has a good reduction at $p$, then $V$ is crystalline [Fa] and  the semisimplicity of $\varphi$
is a well known (and difficult) conjecture. By a result of Katz and  Messing \cite{KM}
$\Dc (V)^{\Ph=1} \ne 0$ can occur only if $i=2m$. Therefore up to eventually replace $V$ by $V^*(1)$ the
conditions {\bf C1, C3-4)} conjecturally hold with except the weight $-1$ case 
$i=2m-1.$ 

The condition $\Dc (V)^{\Ph=1}=0$ imples that the exponential map
$t_V(\Bbb Q_p)@>>>H^1_f(\Bbb Q_p,V)$ is an isomorphism and we denote
by $\log_V$ its inverse.
The composition of the localisation map $\loc$ with the Bloch-Kato logarithm 
$$
r_V\,:\, H^1_f(V)@>>> t_V(\Bbb Q_p)
$$
coincides 
conjecturally  with the $p$-adic (syntomic) regulator. We remark that if $H^0(\Bbb Q_v,V)=0$
for all $v\neq p$ (and therefore $H^1_f(\Bbb Q_v,V)=0$ for all $v\neq p$) then $\loc$ is injective
for all $m\neq i/2, i/2+1$
by a result of Jannsen (\cite{Ja}, Lemma 4 and Theorem 3).

If $H^0_S(V)\neq 0$, then $V$ contains a trivial subextension $V_0=\Bbb Q_p(0)^k$.
For $\Bbb Q_p(0)$ our theory describes  the behavior of the Kubota-Leopoldt
$p$-adic $L$-function and is well known. Therefore we can assume
that $H^0_S(V)=0.$ Applying the same argument to $V^*(1)$ we can also assume that $H^0_S(V^*(1))=0.$
\newline
\newline
From our assumptions we obtain an exact sequence
$$
0@>>>H^1_f(V)@>>>H^1_{f,\{p\}}(V)@>>>\frac{H^1(\Bbb Q_p,V)}{H^1_f(\Bbb Q_p,V)}@>>>0. \tag{3.1}
$$
Moreover 
$$
\aligned
&\dim_{\Bbb Q_p} H^1_f(V)=\dim_{\Bbb Q_p} t_V(\Bbb Q_p)-d_{+}(V),\\
&\dim_{\Bbb Q_p} H^1_{f,\{p\}}(V)=d_{-}(V)+\dim_{\Bbb Q_p} H^0(\Bbb Q_p,V^*(1)).
\endaligned
\tag{3.2}
$$
{\bf 3.1.3. Regular submodules.} In the remainder of this \S we assume that $V$ satisfies {\bf C1-5)}.

\proclaim\nofrills{Definition} {\smc (Perrin-Riou)}. 1) A $\varphi$-submodule $D$ of $\Dc (V)$ is regular 
if $D\cap \F^0\Dc (V)=0$ and the map
$$
r_{V,D}\,:\, H^1_f(V)@>>>\Dc (V)/(\F^0\Dc (V)+D)
$$
induced by $r_V$ is an isomorphism.

2) Dually, a $(\varphi,N)$-submodule $D$ of $\Dc (V^*(1))$ is regular 
if $D+\F^0\Dc (V^*(1))=\Dc (V^*(1))$ and  the map 
$$
D\cap \F^0\Dc (V^*(1)) @>>> H^1(V)^*
$$
induced by the dual map $r_V^*\,:\, \F^0\Dc (V^*(1)) @>>>H^1(V)^*$
is an isomorphism.
\endproclaim

It is easy to see that if $D$ is a regular submodule of $\Dc (V)$, then
$$
D^{\perp}=\text{\rm Hom} (\Dc (V)/D, \Dc (\Bbb Q_p(1))
$$
is a regular submodule of $\Dc (V^*(1)).$ From (3.2)
we also obtain that
$$
\dim D=d_{+}(V), \qquad \dim D^{\perp}=d_{-}(V)=d_{+}(V^*(1)).
$$

Let $D\subset \Dc (V)$ be a regular subspace. As in \cite{Ben2} we 
use the semisimplicity of $\Ph$ to decompose $D$ into the direct sum
$$
D=D_{-1}\oplus D^{\Ph=p^{-1}}.
$$
which gives a four step filtration
$$
\{0\}\subset D_{-1}\subset D \subset \Dc (V).
$$
Let $\bold D$ and $\bold D_{-1}$  denote  the $(\Ph,\Gamma)$-submodules associated  
to $D$ and $D_{-1}$ by Berger's theory, thus
$$
D=\CDcris (\bold D),\qquad  D_{-1}=\CDcris (\bold D_{-1}).
$$ 
Set $W=\text{\rm gr}_0\Ddagrig (V).$ Thus we have two tautological exact sequences
$$
\align
&0@>>> \bold D @>>>\Ddagrig (V) @>>>\bold D'@>>>0,\\
&0@>>>\bold D_{-1} @>>>\bold D@>>>W@>>>0.
\endalign
$$
Note the following properties of cohomology of these modules:

a) The natural maps $H^1(\bold D_{-1}) @>>>H^1(\bold D)$ and
$H^1(\bold D) @>>>H^1(\Ddagrig (V))=H^1(\Bbb Q_p,V)$ are injective. 
This follows from the observation that $\CDcris (\bold D')^{\Ph=1}=0$  by
{\bf C4)}. Since $H^0(\bold D')=\F^0\CDcris (\bold D')^{\Ph=1}$ (\cite{Ben2}, Proposition 1.4.4)
we have $H^0(\bold D')=0$. The same argument works for $W$.

b) $H^1_f(\bold D_{-1})=H^1(\bold D_{-1})$. In particular the exponential map
$\exp_{\bold D_{-1}}\,:\,D_{-1}@>>>H^1(\bold D_{-1})$ is an isomorphism.
This follows from the computation of dimensions of  $H^1(\bold D_{-1})$
and $H^1_f(\bold D_{-1})$. Namely, since $D_{-1}^{\varphi=1}=D_{-1}^{\Ph=p^{-1}}=\{0\}$
the Euler-Poincar\'e characteristic formula \cite{Li} together with Poincar\'e duality give
$$
\dim_{\Bbb Q_p} H^1(\bold D_{-1})=\text{\rm rg} (\bold D_{-1})-\dim_{\Bbb Q_p}H^0(\bold D_{-1})-
\dim_{\Bbb Q_p}H^0(\bold D^*_{-1}(\chi))
=\dim_{\Bbb Q_p} (D_{-1}).
$$
On the other hand since 
$$
\F^0 D_{-1}=D_{-1}\cap \F^0\Dc (V)=\{0\}
$$  
one has 
$\dim_{\Bbb Q_p} H^1_f(\bold D_{-1})=\dim_{\Bbb Q_p} (D_{-1})$
by \cite{Ben2}, Corollary 1.4.5.

c) The exponential map $\exp_{\bold D} \,:\,D@>>>H^1_f (\bold D)$ is an isomorphism.
This follows from $\F^0D=\{0\}$ and $D^{\varphi=1}=\{0\}.$
\newline
\newline
The regularity of $D$ is equivalent to the decomposition
$$
H^1_f(\Bbb Q_p,V) =H^1_f(V) \oplus H^1_f (\bold D). \tag{3.3}
$$

Since  $\loc$ is injective by {\bf C5)}, the localisation map
$
H^1_{f,\{p\}}(V) @>>> H^1(\Bbb Q_p,V)
$
is also injective. Let 
$$
\kappa_D\,\,:\,\,H^1_{f,\{p\}}(V) @>>> \frac{H^1(\Bbb Q_p,V)}{H^1_f(\bold D)}
$$
denote the composition of this map with the canonical projection.

\proclaim{Lemma 3.1.4} i) One has
$$
H^1_f(\Bbb Q_p,V) \cap H^1(\bold D)=H^1_f (\bold D).
$$

ii) $\kappa_D$ is an isomorphism.

\endproclaim
\demo{Proof} i) Since $H^0(\bold D')=0$ we have a commutative diagram
with exact rows and injective colomns
$$
\CD
0 @>>> H^1_f(\bold D) @>>> H^1_f(\Bbb Q_p,V)  @>>> H^1_f(\bold D')\\
@. @VVV @VVV @VVV\\
0 @>>> H^1(\bold D) @>>> H^1(\Bbb Q_p,V)  @>>> H^1(\bold D').
\endCD
$$
This gives i).

ii) Since $H^1_f(\bold D)\subset H^1_f(\Bbb Q_p,V)$ one has
$
\ker (\kappa_D)\subset H^1_f(\Bbb Q_p, V).
$ 
One the other hand (3.3) shows that $\kappa_D$ is injective on $H^1_f(V).$
Thus $\ker (\kappa_D)=\{0\}.$ On the other hand, because
$
\dim_{\Bbb Q_p} H^1_f(\bold D)= \dim_{\Bbb Q_p}(D)
$
we have 
$$
\dim_{\Bbb Q_p} \left (\frac{H^1(\Bbb Q_p,V)}{H^1_f(\bold D)}\right )=d_{-}(V)+
\dim_{\Bbb Q_p}H^0(\Bbb Q_p,V^*(1)).
$$
Comparing this with (3.2) we obtain that $\kappa_D$ is an isomorphism.
% We only need to show that the projection of $H^1(D,V)$
%on the quotient is injective. Let $x\in H^1(\bold D)$. The exact sequence (3.1) 
%shows that there exists
%$\alpha \in H^1_{f,\{p\}}(V)$ such that $\alpha -x\in H^1_f(\Bbb Q_p,V).$
%By (3.3) we can write $\alpha-x=\alpha_1+x_1$ where $\alpha_1 \in H^1_f(V)$
%and $x_1\in H^1_f(\bold D).$ (Here we identify elements of
%different cohomology groups with their images in $H^1(\Bbb Q_p,V)$ under
%canonical injections.) Then $\alpha-\alpha_1=x+x_1\in H^1(\bold D)$ and
%$\kappa_D(\alpha-\alpha_1)= x\pmod{H^1_f(\bold D)}.$ The lemma is proved.
\enddemo
\flushpar
{\bf 3.1.5. The main construction.} Set $e=\dim_{\Bbb Q_p}(D^{\Ph=p^{-1}})$. The $(\Ph,\Gamma)$-module $W$ satisfies 
$$
\F^0\CDcris (W)=0, \qquad \CDcris (W)^{\Ph=p^{-1}}=\CDcris (W).
$$
(Recall that $\CDcris (W)=D^{\Ph=p^{-1}}$.) The cohomology of such  modules was studied
in detail in \cite{Ben2}, Proposition 1.5.9 and section 1.5.10. Namely, $H^0(W)=0,$
$\dim_{\Bbb Q_p}H^1(W)=2e$ and $\dim_{\Bbb Q_p}(W)=e.$ There exists a canonical decomposition
$$
H^1(W)=H^1_f(W)\oplus H^1_c(W)
$$
of $H^1(W)$ into the direct sum of $H^1_f(W)$ and some canonical space $H^1_c(W).$
Moreover there exist canonical isomorphisms
$$
i_{D,f}\,\,:\,\,\CDcris (W)\simeq H^1_f(W), \qquad i_{D,c}\,\,:\,\,\CDcris (W)\simeq H^1_c(W).
$$
These isomorphisms can be described explicitly.
By Proposition 1.5.9 of \cite{Ben2} 
$$
W\simeq \underset{i=1}\to{\overset
e\to\oplus} D_{m_i},  
$$
where $D_{m_i} = \scr R(\vert x\vert x^{m_i}),$ $m_i\geqslant 1.$ 
By Proposition 2.1.5  $H^1_f(D_m)$ is    generated
by $\alpha^*_m $ and  $H^1_c(D_m)$ is the subspace generated by $\beta^*_m$
(see also Proposition 1.1.9). Then 
$$
i_{D_m,f}(x)=x\alpha_m^*,\qquad i_{D_m,c}(x)=x\beta_m^*.
$$

Since $H^0(W)=0$ and $H^2(\bold D_{-1})=0$ we have exact sequences
$$
\align
&0@>>>H^1(\bold D_{-1})@>>>H^1(\bold D)@>>>H^1(W)@>>>0,\\
&0@>>>H^1_f(\bold D_{-1})@>>>H^1_f(\bold D)@>>>H^1_f(W)@>>>0.
\endalign
$$
Since $H^1_f(\bold D_{-1})=H^1(\bold D)$ we obtain that
$$
\dsize\frac{H^1(\bold D)}{H^1_f(\bold D)} \simeq
\dsize\frac{H^1(W)}{H^1_f(W)}.
$$
 Let $ H^1(D,V)$ denote the inverse image of
$H^1(\bold D)/H^1_f(\bold D)$ by $ \kappa_D.$ Then
$\kappa_D$ induces an isomorphism
$$
H^1(D,V) \simeq \frac{H^1(\bold D)}{H^1_f(\bold D)}.
$$
By Lemma 3.1.4  the localisation
map $ H^1(D,V)@>>>H^1(W)$ is well defined and injective.
Hence, we have a diagram
$$
\xymatrix{
\CDcris (W) \ar[r]^{\overset{i_{D,f}}\to \sim} &H^1_f(W)\\
H^1(D,V) \ar[u]^{\rho_{D,f}} \ar[r] \ar[d]_{\rho_{D,c}} 
&H^1(W) \ar[u]_{p_{D,f}}
\ar[d]^{p_{D,c}}
\\
\CDcris (W) \ar[r]^{\overset{i_{D,c}}\to\sim} &H^1_c(W),}
$$
where $\rho_{D,f}$ and $\rho_{D,c}$ are defined as the unique maps making
this diagram commute.
From Lemma 3.1.4 iii)  it follows that  $\rho_{D,c}$ is an isomorphism.
The following definition generalise (in the crystalline case) 
the main construction of \cite{Ben2} where we assumed
in addition that $H^1_f(V)=0.$

\proclaim{Definition} The determinant
$$
\scr L (V,D)= \det \left ( \rho_{D,f} \circ \rho^{-1}_{D,c}\,\mid \,\CDcris (W) \right )
$$
will be called  the $\scr L$-invariant associated to $V$ and $D$.
\endproclaim

{\,}
\flushpar
{\bf 3.2.  $\scr L$-invariant and the large exponential map.}
\newline
{\bf 3.2.1. Derivation of the large exponential map.} In this section we interpret  $\scr L(V,D)$ in terms
of the derivative of the large exponential map.
This interpretation is crucial for the  proof of the main theorem of this paper.
Recall that  $H^1(\Bbb Q_p,\Cal H(\Gamma)\otimes_{\Bbb Q_p}V)= \Cal H(\Gamma)\otimes_{\Lambda (\Gamma)} H_{\Iw}^1(\Bbb Q_p,V)$
injects into $\Ddagrig (V).$  Set 
$$
\align
&F_0H^1(\Bbb Q_p,\Cal H(\Gamma)\otimes_{\Bbb Q_p}V)= \bD
\cap H^1(\Bbb Q_p,\Cal H(\Gamma)\otimes_{\Bbb Q_p}V),\\
&F_{-1}H^1(\Bbb Q_p,\Cal H(\Gamma)\otimes_{\Bbb Q_p}V)= \bD_{-1}
\cap H^1(\Bbb Q_p,\Cal H(\Gamma)\otimes_{\Bbb Q_p}V).
\endalign
$$ 
As in section 2.2 we fix a generator $\gamma \in \Gamma$. 
The following result is a strightforward generalisation of  \cite{Ben3}, Proposition 2.2.2.
For the convenience of the reader we give here  the proof which is 
the same as in {\it op. cit.} modulo  obvious modifications.
 
\proclaim{Proposition 3.2.2} Let $D$ be an admissible subspace of
$\Dc (V).$ For any  $a \in D^{\Ph=p^{-1}}$ let $\alpha  \in \Cal D(V)$
be such that $\alpha (0)=a .$ Then

i) There exists a unique $\beta \in F_0 H^1(\Bbb Q_p,\Cal
H(\Gamma)\otimes V)$ such that
$$
(\gamma-1)\,\beta =\Exp_{V,h}^\ep (\alpha ).
$$

ii) The composition map
$$
\align &\delta_{D,h}\,:\,D^{\Ph=p^{-1}}@>>>F_0 H^1(\Bbb Q_p,\Cal
H(\Gamma)\otimes V)@>>>H^1(W)\\
&\delta_{D,h} (a)\,=\,\beta  \pmod{H^1(\bD_{-1})}
\endalign
$$
is  given
explicitly by  the following formula: 
$$
\delta_{D,h}(\alpha)\,=\,-(h-1)!\,\left (1-\frac{1}{p}\right )^{
 -1}(\log \chi
(\gamma))^{-1}\, i_{D,c} (\alpha).
$$
\endproclaim
\demo{Proof}  
Since $\Dc (V)^{\Ph=1}=0,$ the operator $1-\Ph$ is invertible on
$\Dc (V)$ and we have a diagram
$$
\xymatrix{ \Cal D(V)^{\Delta=0} \ar[rr]^{\Exp_{V,h}^\ep}
\ar[d]^{\Xi_{V,0}^\ep} & & H^1(\Bbb Q_p, \Cal H(\Gamma)\otimes V)
\ar[d]^{\text{\rm pr}_V}\\
\Dc (V) \ar[rr]^{(h-1)!\exp_V} & & H^1(\Bbb Q_p,V). }
$$
where $\Xi_{V,0}^\ep
(\alpha)=\dsize\frac{1-p^{-1}\Ph^{-1}}{1-\Ph}\,\alpha (0)$
(see (2.1)). If
$\alpha \in D^{\Ph=p^{-1}}\otimes \Bbb Z_p[[X]]^{\psi=0},$ then
$\Xi_{V,0}^\ep(\alpha)=0$ and 
$
\text{\rm pr}_V \left (\Exp_{V,h}^\ep (\alpha)\right )\,=\,0.
$
On the other hand, as  $V^{G_K}=0$ the
 map 
$
\left (\Cal H(\Gamma)\otimes_{\Lambda_{\Bbb Q_p}} \Hi^1(\Bbb Q_p,V)\right )_{\Gamma} @>>> H^1(\Bbb Q_p,V)
$
is injective.  Thus there exists a unique 
$\beta \in \Cal H(\Gamma)\otimes_{\Lambda} \Hi^1 (\Bbb Q_p,T)$
such that 
$
\Exp_{V,h}^\ep (\alpha) \,=\,(\gamma-1)\,\beta.
$
Now take $a\in D^{\Ph=p^{-1}}$ and set
$$
f= a \otimes \ell \left (\frac{(1+X)^{\chi (\gamma)}-1}{X}\right ),
$$
where  
$
\ell (g)\,=\,\dsize\frac{1}{p}\log \left(\frac{g^p}{\Ph (g)}\right ). 
$
 An easy computation shows that
$$
\sum_{\zeta^p=1} \ell \left ( \frac{\zeta^{\chi (\gamma)}(1+X)^{\chi (\gamma)}-1}{\zeta (1+X)-1} \right )\,=\,0.
$$
Thus $f \in D^{\Ph=p^{-1}}\otimes \Bbb Z_p[[X]]^{\psi=0}.$
Write $\alpha$ in the form
$\alpha\,=\,(1-\Ph)\,(1-\gamma)\,(a\otimes \log (X)).$ Then 
$$
\Omega_{V,h}(\alpha)\,=\,(-1)^{h-1}\frac{\log \chi (\gamma_1)}{p}\,t^h\partial^h((\gamma-1)\,(a\log (\pi))\,=\,
\frac{\log \chi (\gamma_1)}{p}\,(\gamma-1)\,\beta
$$
where 
$$
\beta \,=\,(-1)^{h-1}t^h \partial^h (a\log (\pi))\,=\,(-1)^{h-1}a t^h \partial^{h-1}
\left (\frac{1+\pi}{\pi}\right ).
$$
This implies immediately that $\beta \in \bD.$ On the other hand
$
D^{\Ph=p^{-1}}\,=\,\Cal D_{\text{\rm cris}}(W)=(W [1/t])^{\Gamma}
$
and we will write  $\tilde a$ for  the image of $a$ in $W\,\left [1/t \right ].$
By \cite{Ben2}, sections 1.5.8-1.5.10 
one has $W\,\simeq \,\underset{i=1}\to{\overset{e}\to \oplus} \bD_{m_i}$
where $\bD_{m}=\scr R(|x|x^m)$
and we denote by $e_m$ the canonical base of $\bD_m.$
Then  without lost of generality  we may assume that $\tilde a\,=\,t^{-m_i}e_{m_i}$ for some $i.$ 
Let $\tilde \beta$ be
denote the image of $\beta$ in $W^{\psi=1}$ and let
$
h^1_0 \,\,:\,\, W^{\psi=1}@>>>H^1(W)
$
be the canonical map furnished by Proposition 1.1.7.  Recall that $h^1_0(\tilde \beta )=\cl (c,\tilde\beta)$ where 
$(1-\gamma)\,c\,=\,(1-\Ph)\,\tilde \beta.$ 
Then $\tilde \beta= (-1)^{h-1}t^{h-m_i}\partial^h \log (\pi)$. By Lemma 1.5.1 of \cite{CC1}
there exists a unique $b_0\in \boB^{\dag, \psi=0}_{\Bbb Q_p}$ 
such that $(\gamma-1)\,b_0=\ell (\pi)$. This implies that 
 $$
(1-\gamma)\,(t^{h-m_i}\partial^h b_0e_{m_i})\,=\,(1-\Ph)\,(t^{h-m_i}\partial^h \log (\pi)e_{m_i})\,=\,(-1)^{h-1}
(1-\Ph)\,\tilde \beta.
$$
Thus $c= (-1)^{h-1}t^{h-m_i}\partial^h b_0e_{m_i}$ and 
$
\res (ct^{m_i-1}dt)=(-1)^{h-1}\res (t^{h-1}\partial^h b_0dt)\,e_{m_i}=0.
$
Next from the congruence
$
\tilde \beta \equiv (h-1)!\,t^{-m_i}e_{m_i}\,
\pmod{\Bbb Q_p[[\pi]]\,e_{m_i}}.
$
it follows  that
$\res (\tilde \beta t^{m_i-1}dt)\,=\,(h-1)!\,e_{m_i}.$ Therefore by \cite{Ben2}, Corollary 1.5.6  we have 
$$
\cl (c,\tilde \beta)\,=\,(h-1)!\,\cl (\beta_m)\,=\,(h-1)! \dsize \frac{p}{\log \chi (\gamma_1)}\,i_{W,c}(a).
$$
 On the other hand 
$$
\alpha (0)\,=\,\left.
a \otimes \ell \left (\frac{(1+X)^{\chi (\gamma)}-1}{X}\right )
 \right \vert_{X=0}\,=\,a\left (1-\frac{1}{p} \right )\,\log (\chi (\gamma)).
$$
These formulas imply that
$$
\delta_{D,h}(\alpha)=(h-1)!\,\left (1-\frac{1}{p}\right )^{
 -1}(\log \chi
(\gamma))^{-1}\, i_{W,c} (\alpha).
$$
and the proposition is proved.
\enddemo
\flushpar
{\bf 3.2.3. Interpretation of the $\scr L$-invariant.}  From the definition of $H^1(D,V)$ and Lemma 3.1.4 we immediately obtain that
$$
\frac{H^1(\Bbb Q_p,V)}{H^1_{f,\{p\}}(V)+ H^1(\bD_{-1})}\simeq 
\frac{H^1(\bD)}{H^1(D,V)+ H^1(\bD_{-1})} \simeq \frac{H^1(W)}{H^1(D,V)}\,.
$$ 
Thus, the map $\delta_{D,h}$ constructed in Proposition 3.2.2 induces a map
$$
D^{\Ph=p^{-1}} @>>> 
\dsize\frac{H^1(\Bbb Q_p,V)}{H^1_{f,\{p\}}(V) +H^1(\bD_{-1})}
$$
which we will denote again by $\delta_{D,h}.$  On the other hand, we have isomorphisms
$$
D^{\Ph=p^{-1}} \overset{\exp_{V}}\to \iso \frac{H^1_f(\Bbb Q_p,V)}{\exp_{V,\Bbb Q_p}(D_{-1})}
\simeq  \frac{H^1_f(\Bbb Q_p,V)}{H^1(\bD_{-1})}
\simeq
\frac{H^1(\Bbb Q_p,V)}{H^1_{f,\{p\}}(V) +H^1(\bD_{-1})}.
$$
\proclaim{Proposition 3.2.4} Let $\lambda_D \,:\,D^{\Ph=p^{-1}}@>>> D^{\Ph=p^{-1}}$ denote the homomorphism
making the diagram
$$\xymatrix{
D^{\Ph=p^{-1}} \ar[dr]^{\delta_{D,h}} \ar[rr]^{\lambda_D}
& & D^{\Ph=p^{-1}} \ar[dl]_{(h-1)!\exp_{V}\,\,\,\,}\\
 & {\dsize\frac{H^1(\Bbb Q_p,V)}{H^1_{f,\{p\}}(V)+H^1(\bD_{-1})}} &
 }
$$
commute. 
Then 
$$
\det \left (\lambda_D\,\vert \, D^{\Ph=p^{-1}} \right )\,=\,(\log \chi (\gamma))^{-e}\left ( 1-\frac{1}{p} \right )^{-e}
 \scr L(V,D).
$$
\endproclaim
\demo{Proof}  The proposition follows from  Proposition 3.2.2 and the following 
elementary fact. Let $U=U_1\oplus U_2$ be the decomposition of a vector space $U$ of dimension
$2e$ into the direct sum of two subspaces of dimension $e$. Let $X\subset U$ be a subspace
of dimension $e$ such that $X\cap U_1=\{0\}.$ Consider the  diagrams 
$$
\xymatrix{
X \ar[r]^{p_1} \ar[d]_{p_2} &U_1  &  U/X  &  U_1 \ar[l]_{i_1} \\
U_2 \ar[ur]_f &  &  U_2 \ar[u]^{i_2} \ar[ur]_g & }
$$
where $p_k$ and $i_k$ are induced by natural projections and inclusions.
Then $f=-g.$ Applying this remark to $U=H^1(W),$
$X=H^1(D,V)$, $U_1=H^1_f(W),$ 
$U_2=H^1_c(W)$ and taking determinants we obtain the proposition. 
\enddemo

\head {\bf \S4. Special values of $p$-adic $L$-functions}
\endhead

\flushpar
{\bf 4.1. The Bloch-Kato conjecture.} 
\newline
{\bf 4.1.1. The  Euler-Poincar\'e line} (see \cite{F3}, \cite{FP},\cite{BF}).
Let $V$ be a $p$-adic pseudo-geometric representation
of $\text{\rm Gal} (\overline{\Bbb Q}/\Bbb Q).$  Thus $V$ is a
finite-dimensional $\Bbb Q_p$-vector space equipped with a
continuous action of the Galois group $G_S$ for a suitable finite
set of places  $S$ containing $p$. 
Write
$
\RG_S (V)= C_c^\bullet (G_S,V)
$
and define
$$
\RG_{S,c}(V)=\text{\rm cone} \left (\RG_S(V)@>>> \underset{v\in S\cup \{\infty\}}
\to \oplus \RG (\Bbb Q_v,V) \right )\,[-1].
$$

Fix a $\Bbb Z_p$-lattice $T$  of $V$ stable under the action of
$G_S$ and set 
$
\Delta_S(V)={\det}^{-1}_{\Bbb Q_p} \RG_{S,c}(V)
$
and $\Delta_S(T)={\det}^{-1}_{\Bbb Z_p} \RG_{S,c}(T).$
Then $\Delta_S(T)$ is a $\Bbb Z_p$-lattice of the one-dimensional
$\Bbb Q_p$-vector space $\Delta_S(V)$ which does not depend on the
choice of $T$. Therefore it defines a
$p$-adic norm on $\Delta_S(V)$ which we denote by $\Vert \,\cdot
\,\Vert_S.$ Moreother, $(\Delta_S(V), \Vert \,\cdot \,\Vert_S)$
does not depend on the choice of $S.$ More precisely, if $\Sigma$
is a finite set of places which contains $S$, then there exists a
natural isomorphism $\Delta_S(V) @>>>\Delta_{\Sigma}(V)$ such that
$\Vert \,\cdot \,\Vert_\Sigma=\Vert \,\cdot \,\Vert_S.$ This allows
to define the Euler-Poincar\'e line $\Delta_{\text{\rm
EP}}(V)$ as $(\Delta_S(V),\Vert \,\cdot \,\Vert_S)$ where
$S$ is sufficiently large. Recall that for any finite place $v\in
S$ we defined
$$
\RG_f (\Bbb Q_v,V)\,=\,\cases \left [ V^{I_v} @>1-f_v>> V^{I_v}
\right ] &{\text{\rm if $v\ne p$}}\\
\left [\Dc (V)@>(\text{\rm pr},1-\Ph)>> t_V(\Bbb Q_p)\oplus \Dc
(V) \right ] &{\text{\rm if $v=p$}}.
\endcases
$$
At $v=\infty$  we set
$
\RG_f(\Bbb R,V)\,=\,\left [V^+@>>>0\right ],
$
where the first term is placed in degree $0.$ Thus $\RG_f(\Bbb R,V) \overset{\sim}\to \rightarrow \RG(\Bbb R,V)$.
For any $v$ we have a canonical morphism $\text{\rm loc}_p\,:\,\RG_f(\Bbb Q_v,V)@>>> \RG (\Bbb Q_v,V)$
which can
be viewed as a local condition in the sense of \cite{Ne2}. Consider
the diagram
$$
\xymatrix{ \RG_S(V) \ar[r] & \underset{v\in S\cup \{\infty\}}\to \oplus
\RG(\Bbb
Q_v,V)\\
 & \underset{v\in S\cup \{\infty\}}\to \oplus \RG_f(\Bbb Q_v,V) \ar[u]
 }
$$
and define
$$
\RG_f(V)=\text{\rm cone} \left (\RG_S(V)\oplus \left
(\underset{v\in S\cup\{\infty\}}\to \oplus \RG_f(\Bbb Q_v,V) \right ) @>>>
\underset{v\in S\cup\{\infty\}}\to\oplus \RG (\Bbb Q_v,V)\right )[-1].
$$
Thus, we have a distinguished triangle
$$
\RG_f(V) @>>>\RG_S(V)\oplus \left
(\underset{v\in S\cup\{\infty\}}\to \oplus \RG_f(\Bbb Q_v,V) \right ) @>>>
\underset{v\in S\cup \{\infty\}}\to\oplus \RG (\Bbb Q_v,V).
\tag{4.1}
$$
Set
$$
\Delta_f (V)={\det}_{\Bbb Q_p}^{-1}\RG_f(V)\otimes
{\det}_{\Bbb Q_p}^{-1}t_V(\Bbb Q_p)\otimes {\det}_{\Bbb Q_p}V^+.
$$
\flushpar
It is easy to see that $\RG_f(V)$ and $\Delta_f (V)$ do
not depend on the choice of $S$. Consider the distinguished triangle
$$
\RG_{S,c}(V)@>>>\RG_f(V)@>>>\underset{v\in S\cup \{\infty\}}\to\oplus \RG_f(\Bbb
Q_v,V).
$$
Since ${\det}_{\Bbb Q_p}\RG_f (\Bbb Q_p,V) \simeq {\det}_{\Bbb Q_p}^{-1}t_V(\Bbb Q_p)$ and
${\det}_{\Bbb Q_p}\RG_f(\Bbb R,V)={\det}_{\Bbb Q_p}V^+$
tautologically, we obtain canonical isomorphisms
$$
\Delta_f(V)\simeq {\det}_{\Bbb
Q_p}^{-1}\RG_{S,c}(V)\simeq \Delta_{\text{\rm EP}}(V).
$$
\flushpar
The cohomology of $\RG_f(V)$ is as follows:
$$
\aligned &\bold R^0\Gamma_f(V)=H_S^0(V),\quad
\bold R^1\Gamma_f(V)=H^1_f(V),\quad
\bold R^2\Gamma_f(V)\simeq H^1_f(V^*(1))^*,\\
&\bold R^3\Gamma_f (V)=\text{\rm coker}\,\left (
H^2_S(V)@>>>\underset{v\in S}\to\oplus H^2(\Bbb Q_v,V)\right )
\simeq H_S^0(V^*(1))^*.
\endaligned
\tag{4.2}
$$
These groups seat in the following exact sequence:
$$
\multline 0@>>>\bold R^1\Gamma(V)@>>>H^1_S(V)@>>>\underset{v\in
S}\to \bigoplus \frac{H^1(\Bbb Q_v,V)}{H^1_f(\Bbb Q_v,V)}
@>>>\bold
R^2\Gamma_f(V)@>>>\\
H^2_S(V)@>>>\underset{v\in S}\to \oplus H^2(\Bbb Q_v,V)@>>>\bold
R^3\Gamma_f(V)@>>>0.
\endmultline
$$
\flushpar
The $L$-function of $V$ is defined as the Euler product
$$
L(V,s)=\underset{v}\to \prod E_v(V,(Nv)^{-s})^{-1}
$$
where
$$
E_v(V,t)=\cases \det \left (1-f_vt\,\vert \,V^{I_v} \right ),
&\text{\rm if $v\ne p$} \\
\det \left (1-\Ph t\,\vert \,\Dc (V)\right ) &\text{\rm if $v=p$}.
\endcases
$$

$\,$
\newline
\flushpar
{\bf 4.1.2. Canonical trivialisations.} In this paper we treat motives in
the formal sense and assume all conjectures about the category of
mixed motives $\Cal M \Cal M$ over $\Bbb Q$ which are  necessary to state the Bloch-Kato
conjecture (see \cite{F3}, \cite{FP}). Let $M$ be a pure motive
over $\Bbb Q$ and let $M_{\text{\rm B}}$ and $M_{\text{\rm dR}}$ denote 
its Betti and de Rham realisations respectively. 
Fix an odd prime $p$  and denote by $V=M_p$ the $p$-adic realisation
of $M$. Then one has comparision isomorphisms
$$
\align
&M_{\text{\rm B}}\otimes_{\Bbb Q}\Bbb C\iso M_{\text{\rm dR}}\otimes_{\Bbb Q}\Bbb C,
\tag{4.3}
\\
&M_{\text{\rm B}}\otimes_{\Bbb Q}\Bbb Q_p\iso V
\tag{4.4}
\endalign
$$
which induce trivialisations
$$
\align
&\Omega_{M}^{(H,\infty)}\,:\, {\det}_{\Bbb Q} M_{\text{\rm B}}\otimes
{\det}_{\Bbb Q}^{-1}M_{\text{\rm dR}} @>>>\Bbb C, \tag{4.5}
\\
&\Omega_{M}^{(\acute et,p)}\,:\,{\det}_{\Bbb Q_p}V \otimes {\det}_{\Bbb Q}^{-1} M_{\text{\rm B}}
 @>>>\Bbb Q_p. \tag{4.6}
\endalign
$$
The complex conjugation acts on $M_{\text{\rm B}}$ and $V$ and 
decomposes the last isomorphism into $\pm$ parts which we denote again 
by $\Omega_{M}^{(\acute et,p)}$ to simplify notation
$$
\Omega_{M}^{(\acute et,p)}\,:\,{\det}_{\Bbb Q_p}V^{\pm}\otimes {\det}^{-1}_{\Bbb Q} M_{\text{\rm B}}^{\pm}
 @>>>\Bbb Q_p. \tag{4.7}
$$
The restriction of  $V$ on the decomposition group at $p$ is a de Rham representation
and $\Dd (V)\simeq M_{\text{\rm dR}}\otimes_{\Bbb Q} \Bbb Q_p.$ 
The comparision isomorphism 
$$
V\otimes \Bd \iso \Dd (V)\otimes \Bd \tag{4.8}
$$
induces a map
$$
\widetilde \Omega_M^{(H,p)}\,:\,{\det}_{\Bbb Q_p}V \otimes {\det}^{-1}_{\Bbb Q_p}
\Dd (V) @>>>\Bd .
$$
It is not difficult to see that there exists a finite extension
$L$ of $\widehat{\Bbb Q_p^{ur}}$ such that $\text{\rm Im}(\tilde \Omega_M^{H,p})\subset  L t^{t_H(V)}$
and we define
$$
\Omega_M^{(H,p)}\,:\,{\det}_{\Bbb Q_p}V \otimes {\det}^{-1}_{\Bbb Q_p}
\Dd (V) @>>>L \tag{4.9}
$$
by $\Omega_M^{(H,p)}=t^{-t_H(V)}\widetilde\Omega_M^{(H,p)}.$
We remark that if $V$ is crystalline at $p$ then one can take
$L=\widehat{\Bbb Q_p^{ur}}$ (see \cite{PR2}, Appendice C.2).  
\newline
\newline
Assume that the groups
$
H^i(M)\,=\,\text{\rm Ext}^i_{\Cal M\Cal M}(\Bbb Q(0),M)
$
are well defined and vanish for $i\ne 0,1.$ It should be possible to  define  a
$\Bbb Q$-subspace $H^1_f(M)$ of $H^1(M)$ consisting of "integral"
classes of extensions which is expected to be finite dimensional.
It is convenient to set $H^0_f(M)=H^0(M).$ 
The conjectures of Tate and Jannsen predict that
the regulator map induces   isomorphisms
$$
H^i_f(M)\otimes_{\Bbb Q}\Bbb Q_p\simeq H^i_f(V), \qquad i=0,1. \tag{4.10}
$$

In this paper $M$ will always denote a motive satisfying the following conditions
\newline
\,

{\bf M1)} $M$ is  pure  of weight $w\leqslant -2$.

{\bf M2)} The $p$-adic realisation $V$ of $M$ is crystalline at $p$.

{\bf M3)} $M$ has no subquotients isomorphic to $\Bbb Q(1).$
\newline
\,
\flushpar
These conditions  imply that $H^0(M)=H^0(M^*(1))=0$ and $H^1(M^*(1))=0$
by the weight argument. and by (4.10) the representation $V$ should 
satisfy the conditions {\bf C1,2,4)} of section 3.1.2. In particular,
from (4.2) it follows that
$${\det}_{\Bbb Q_p} \RG_f(V) \iso {\det}^{-1}_{\Bbb Q_p} H^1_f (V). \tag{4.11}
$$
The semisimplicity of $\Ph$ is a well known conjecture which is actually known
for abelian varieties. Finally {\bf C5)} should follow from the injectivity
of the syntomic regulator.
\newline
\newline
The  comparision isomorphism (4.3) induces an injective  map
$$
\alpha_M\,:\,M_{\text{\rm B}}^+\otimes_{\Bbb Q}\Bbb R@>>>t_M(\Bbb R)
$$
and the  six-term exact sequence of Fontaine and Perrin-Riou (\cite{F3}, section 6.10)  
degenerates into an isomorphism (the regulator map)
$$
r_{M,\infty}\,:\,H^1_f(M)\otimes_{\Bbb Q}\Bbb R \iso \text{\rm coker} (\alpha_M).
$$
The maps $\alpha_M$ and $r_{M,\infty}$ define a map
$$
R_{M,\infty}\,:\,{\det}_{\Bbb Q}^{-1}t_M(\Bbb Q) \otimes
{\det}_{\Bbb Q}M_{\text{\rm B}}^+ \otimes {\det}_{\Bbb Q} H^1_f (M)@>>> \Bbb R
$$
Fix bases $\omega_f\in {\det}_{\Bbb Q}H^1_f(M)$, $\omega_{t_M}\in {\det}_{\Bbb Q}t_M(\Bbb Q)$ and 
$\omega_{M_{\text{\rm
B}}}^+\in {\det}_{\Bbb Q} M_{\text{\rm B}}^+$. Set
$\omega_{M}=(\omega_f,\omega_{t_M},\omega_{M_{\text{\rm
B}}}^+)$ and define
$$
R_{M,\infty}(\omega_{M})=R_{M,\infty}(\omega_{t_M}^{-1}\otimes \omega_{M_{\text{\rm
B}}}^+\otimes \omega_f).
$$
Using (4.11) and the isomorphisms (4.10) define 
$$
i_{\omega_{M},p}\,\,:\,\,\Delta_{\text{\rm EP}}(V)\overset{\sim}\to \rightarrow
{\det}_{\Bbb Q_p}^{-1} t_V(\Bbb Q_p) \otimes {\det}_{\Bbb Q_p}
V^+ \otimes {\det}_{\Bbb Q_p} H^1_f(V) @>>> \Bbb Q_p \tag{4.12}
$$
by
$
x\,=\,i_{\omega_M,p} (x)\,(\omega_{t_M}^{-1}\otimes \omega_{M_{\text{\rm B}}}^+\otimes \omega_f).
$
\newline
\newline
Consider now the case of the dual motive $M^*(1).$ Again one has
$
\Delta_{EP}(V^*(1)) \simeq \Delta_f(V^*(1))
$
where
$$\Delta_f(V^*(1))\simeq {\det}^{-1}_{\Bbb Q_p} t_{V^*(1)}(\Bbb Q_p)\otimes 
{\det}_{\Bbb Q_p} V^*(1)^+\otimes {\det}_{\Bbb Q_p} H^1_f(V).
$$
The map $\alpha_{M^*(1)}\,:\,M^*(1)_{\text{\rm B}}^+\otimes_{\Bbb Q}\Bbb R@>>>
t_{M^*(1)}(\Bbb R)$ is surjective  and it is related to $\alpha_{M}$ by the canonical duality
$\text{\rm coker} (\alpha_M)\times \ker (\alpha_{M^*(1)}) @>>>\Bbb R
$
(see \cite{F3}, section 5.4).
The six-term exact sequence degenerates into an isomorphism
$$
r_{M^*(1),\infty}\,:\, H^1_f(M)^*\otimes_{\Bbb Q}\Bbb R\simeq \ker (\alpha_{M^*(1)}).
$$
This allows to define a map
$$
R_{M^*(1),\infty}\,\,:\,\, 
{\det}_{\Bbb Q}^{-1}t_{M^*(1)}(\Bbb Q) \otimes
{\det}_{\Bbb Q}M^*(1)_{\text{\rm B}}^+ \otimes {\det}_{\Bbb Q} H^1_f (M)@>>> \Bbb R.
$$
Fix bases $\omega_{t_{M^*(1)}}\in {\det}_{\Bbb Q}t_{M^*(1)}(\Bbb Q)$ and 
$\omega_{M^*(1)_{\text{\rm B}}}^+\in {\det}_{\Bbb Q}M^*(1)_B^+.$
Set $\omega_{M^*(1)}=(\omega_{t_{M^*(1)}},\omega_{M^*(1)_{\text{\rm B}}}^+,\omega_f)$
and $R_{M^*(1),\infty}(\omega_{M^*(1)})=
R_{M^*(1),\infty} (\omega_{t_{M^*(1)}}^{-1}\otimes \omega_{M^*(1)_{\text{\rm B}}}^+ 
\otimes \omega_f)$.
 Then again 
this data defines   a trivialisation  
$$
i_{\omega_{M^*(1)},p}\,\,:\,\,\Delta_{\text{\rm EP}}(V^*(1))) @>>> \Bbb Q_p. \tag{4.13}
$$
It is conjectured that the $L$-functions $L(V,s)$ and $L(V^*(1),s)$ 
are well defined complex functions
have meromorphic continuation to the whole $\Bbb C$ and satisfy
some explicit  functional equation (\cite{FP} chapitre III).
One expects that they do not depend on the choice of the prime $p$
and we will denote them by $L(M,s)$ and $L(M^*(1),s)$ respectively.
The conjectures about special values of these functions state as follows.

\proclaim\nofrills{Conjecture } $\,\,${\smc (Beilinson-Deligne)}. 
The $L$-function  $L(M,s)$ does not vanish
at $s=0$ and
$$
\frac{L(V,0)}{R_{M,\infty}(\omega_M)}\in \Bbb Q^*.
$$
The $L$-function $L(M^*(1),s)$ has a zero of order $r=\dim_{\Bbb Q_p} H^1_f(M)$ at $s=0$.
Let $L(M^*(1),0)=\lim_{s\to 0} s^{-r}L(M^*(1),s).$ Then
$$
\frac{L(M^*(1),0)}{R_{M^*(1),\infty}(\omega_{M^*(1)})}\in \Bbb Q^*.
$$
\endproclaim

\proclaim\nofrills{Conjecture} $\,\,${\smc (Bloch-Kato)}. Let $T$  be a $\Bbb Z_p$-lattice of $V$   stable under the action of $G_S.$ Then
$$
\aligned
&i_{\omega_M,p}(\Delta_{\text{\rm
EP}}(T))\,=\,\frac{L(M,0)}{R_{M,\infty}(\omega_M)}\,\Bbb
Z_p,\\
&i_{\omega_{M^*(1)},p}(\Delta_{\text{\rm
EP}}(T^*(1)))\,=\,\frac{L(M^*(1),0)}
{R_{M^*(1),\infty}(\omega_{M^*(1)})}\,\Bbb
Z_p.
\endaligned
$$
\endproclaim 
$\,$
\newline
{\bf 4.1.3. Compatibility with functional equation.} The compatibility of the Bloch-Kato conjecture with the functional equation
follows from the conjecture $C_{EP}(V)$ of Fontaine and Perrin-Riou about
local Tamagawa numbers (\cite{FP}, chapitre III, section 4.5.4). More precisely,  define
$$
\Gamma^* (V)=\prod_{i\in \Bbb Z} \Gamma^*(-i)^{h_i(V)} \tag{4.14}
$$
where 
$h_i(V)=\dim_{\Bbb Q_p}\left (\text{\rm gr}_i(\Dd (V))\right )$ and 
$$
\Gamma (i)=\cases (i-1)! &\text{if $i>0$}\\
\frac{(-1)^i}{(-i)!} &\text{if $i\leqslant 0$}.
\endcases
$$
The exact sequence
$$
0@>>>t_{V^*(1)} (\Bbb Q_p)^* @>>>\Dd (V)@>>>t_V(\Bbb Q_p)@>>>0
$$
%and the isomorphism
%$
%(M^*(1)_{\text{\rm B}}^+)^* \simeq M_{\text{\rm B}}^-
%$
allows to consider $\omega_{M_{\text{\rm dR}}}=\omega_{t_M}\otimes \omega_{t_{M^*(1)}}^{-1}\in {\det}_{\Bbb Q_p}\Dd (V)$.
%and $\omega_{\text{\rm B}}\otimes (\omega_{\text{\rm B}}^*)^{-1}\in {\det}_{\Bbb Q}%M_{\text{\rm B}}.$
Choose bases  $\omega_T^+\in {\det}_{\Bbb Z_p}T^+$ and $\omega_T^-\in {\det}_{\Bbb Z_p}T^-$ 
and set $\omega_T=\omega_T^+\otimes\omega_T^-\in {\det}_{\Bbb Z_p} T$ and
$\omega_{T^*(1)}^+=(\omega_T^-)^*\in {\det}_{\Bbb Z_p}T^*(1)^+.$
Then the conjecture $C_{EP}(V)$ implies that
$$
\frac{{i_{\omega_{M^*(1)},p}(\Delta_{\text{\rm
EP}}(T^*(1)))}}{\Omega_{M^*(1)}^{(\acute et,p)}(\omega_{T^*(1)}^+,\omega_{M^*(1)_{\text{\rm B}}}^+)}\,=\,\Gamma^*(V)\,
\,\Omega_M^{(H,p)}(\omega_T,\omega_{M_{\text{\rm dR}}})\,
\frac{i_{\omega_M,p}(\Delta_{\text{\rm
EP}}(T))}
{\Omega_{M}^{(\acute et,p)}(\omega_T^+,\omega_{M_{\text{\rm B}}}^+)}
\tag{4.15}
$$
(see \cite{PR2}, Appendice C). We remark that for  crystalline representations
$C_{EP}(V)$ is proved in \cite{BB08}.
$$
$$
\flushpar 
{\bf 4.2. $p$-adic $L$-functions.}
\newline
{\bf 4.2.1. $p$-adic Beilinson's conjecture.} We keep previous notation and conventions. Let $M$ be a motive which satisfies the conditions {\bf M1-3)} of section 4.1.2 and let $V$ denote
the $p$-adic realisation of $M$. We fix bases $\omega_{M_{\text{\rm B}}}^+\in {\det}_{\Bbb Q}M_{\text{\rm B}}^+$,
$\omega_{t_M}\in {\det}_{\Bbb Q}t_M(\Bbb Q)$ and $\omega_f\in {\det}_{\Bbb Q}H^1_f(M).$
We also fix a lattice $T$ in $V$ stable under the action of $G_S$ and a base
$\omega_T^+\in {\det}_{\Bbb Z_p}T^+$. To simplify notation we will assume that 
the choices of $\omega_{M_{\text{\rm B}}}^+ $ and $\omega_T^+$ are compatible, namely that 
$\Omega_M^{(\acute et,p)} (\omega_T^+,\omega_{M_{\text{\rm B}}}^+)=1.$
 Let $D$ be a regular subspace of $\Dc (V).$
We fix a $\Bbb Z_p$-lattice  $N$ of $D$ and  a basis $\omega_N\in {\det}_{\Bbb Z_p}N.$
By the analogy with the archimedian case we can consider the $p$-adic
regulator as a map  $r_{V,D}\,:\,H^1_f(V)@>>>\text{\rm coker} (\alpha_{V,D})$ where
$$
\alpha_{V,D}\,:\, D@>>>t_V(\Bbb Q_p)
$$
is the natural projection. Set $\omega_{V,N}=(\omega_{t_M},\omega_N,\omega_f)$
and denote by $R_{V,D}(\omega_{V,N})$ 
the determinant of $r_{V,D}$ 
computed in the bases $\omega_f$ and $\omega_{t_M}\otimes \omega_N^{-1}.$ Namely,
$R_{V,D}(\omega_{V,N})$ is the image of 
$\omega_{t_M}^{-1}\otimes \omega_N\otimes \omega_f$ under the induced isomorphism
$$
R_{V,D}\,\,:\,\,
{\det}^{-1}_{\Bbb Q_p}t_{V}(\Bbb Q_p) \otimes 
{\det}_{\Bbb Q_p} D\otimes {\det}_{\Bbb Q_p}H^1_f(V) @>>>\Bbb Q_p.
$$
Now, consider the projection
$$
\alpha_{V^*(1),D^{\perp}}\,:\,D^{\perp}@>>>t_{V^*(1)}(\Bbb Q_p).
$$
A standard argument from the linear algebra shows that $\alpha_{V^*(1),D^{\perp}}$
is surjective and is related to $\alpha_{V,D}$ by the canonical duality
$\text{\rm coker}(\alpha_{V,D})\times \ker (\alpha_{V^*(1),D})@>>>\Bbb Q_p.$
This  defines  isomorphisms
$$
{\det}^{-1}_{\Bbb Q_p}t_{V^*(1)}(\Bbb Q_p) \otimes 
{\det}_{\Bbb Q_p} D^{\perp} \simeq  {\det}_{\Bbb Q_p}(\text{\rm ker}(\alpha_{V^*(1),D^{\perp}}))
\simeq
{\det}_{\Bbb Q_p}^{-1} (\text{\rm coker}(\alpha_{V,D}))
$$
and composing this map with the determinant of $r_{V,D}$ we have again a
trivialisation
$$
R_{V^*(1),D^{\perp}}\,\,:\,\,{\det}^{-1}_{\Bbb Q_p}t_{V^*(1)}(\Bbb Q_p) \otimes 
{\det}_{\Bbb Q_p} D^{\perp}\otimes {\det}_{\Bbb Q_p}H^1_f(V) @>>>\Bbb Q_p.
$$
Choose a lattice $N^{\perp}\subset D^{\perp},$  fix bases $\omega_{t_{M^*(1)}}$ and $\omega_{N^{\perp}}$ of ${\det}_{\Bbb Q_p}t_{V^*(1)}(\Bbb Q_p)$
and ${\det}_{\Bbb Z_p}N^{\perp}$ respectively and set $\omega_{V,N^{\perp}}=
(\omega_{t_{M^*(1)}}, \omega_{N^{\perp}},\omega_f)$.

 Perrin-Riou conjectured \cite{PR2} 
that there exists an  
analytic  $p$-adic  $L$-function $L_p(T,\omega_N,s)$ which interpolates special 
values of the complex $L$-function $L(M,s).$ In particular one expects that
if  $p^{-1}$ is not an eigenvalue of $\Ph$ acting on $D$ then $L_p(T,\omega_N,s)$
does not vanish at $s=0$ and 
$$
\frac{L_p(T,N,0)}{R_{V,D}(\omega_{V,N})} \,=\,
\Cal E(V,D)\,\frac{L(M,0)}{R_{M,\infty}(\omega_M)}.
$$
where 
$$
\multline
\Cal E(V,D)=\det (1-p^{-1}\Ph^{-1} \mid D)\,\det (1-p^{-1}\Ph^{-1} \mid D^{\perp})=\\
=\det (1-p^{-1}\Ph^{-1} \mid D)\,\det (1-\Ph \mid \Dc (V)/D).
\endmultline
$$
Dually it is conjectured  that there exists a $p$-adic $L$-function $L_p(T^*(1),\omega_{N^{\perp}},s)$ which interpolates special values of $L(M^*(1),s)$.
One expects that if  $1$ is not an eigenvalue of $\Ph$ acting on the quotent $\Dc (V^*(1))/D^{\perp}$ then $L_p(T^*(1),\omega_{N^{\perp}},s)$
has a zero of order $r=\dim_{\Bbb Q}H^1_f(M)$  at $s=0$ and 
$$
\frac{L_p^*(T^*(1),N^{\perp},0)}{R_{V^*(1),D^{\perp}}(\omega_{V^*(1),N^{\perp}})} \,=\,
\Cal E(V^*(1),D^{\perp})\,\frac{L^*(M^*(1),0)}{R_{M^*(1),\infty}(\omega_M)}.
$$
These properties of $p$-adic $L$-functions can be viewed as $p$-adic analogues
of Beilinson's conjectures and we refer the reader to \cite{PR2}, chapitre 4 and \cite{C2}, section 2.8 for more detail.
Note that from the definition it is clear that $\Cal E(V,D)=\Cal E(V^*(1),D^{\perp}).$ One can also write 
$\Cal E(V,D)$ in the form
$$
\Cal E(V,D)=E_p(V,1)\,{\det} \left (\frac{1-p^{-1}\Ph^{-1}}{1-\Ph}  \vert 
D\right ).
$$

$\,$
\newline
{\bf 4.2.2. Trivial zero conjecture.} Assume now that $D^{\Ph=p^{-1}}\ne 0.$ Since $M$ is crystalline at $p$,
this can occur only if $M$ is of weight $-2.$ Set 
$$
e={\dim}_{\Bbb Q_p}D^{\Ph=p^{-1}}={\dim}_{\Bbb Q_p} (D^{\perp}+\Dc(V^*(1))^{\Ph=1})/D^{\perp}).
$$
Assume that the $p$-adic realisation $V$ of $M$ satisfies the conditions {\bf C1-5)}
of section 3.1.2.
Decompose $D$ into the direct sum $D=D_{-1}\oplus D^{\Ph=p^{-1}}$ and
define
$$
\Cal E^+(V,D)=\Cal E^+(V^*(1),D^{\perp})=\det (1-p^{-1}\Ph^{-1} \mid D_{-1})\,
\det (1-p^{-1}\Ph^{-1} \mid D^{\perp}). \tag{4.15}
$$
We propose the following conjecture about the behavior of
$p$-adic $L$-functions at $s=0.$
{\,}
\newline
\newline
{\bf Trivial zero conjecture.} Let $D$ be a regular subspace of $\Dc (V).$ 
Then

1) The $p$-adic $L$-function $L_p(T,N,s)$ has a zero of order $e$ at $s=0$ and
$$
\frac{L^*_p(T,N,0)}{R_{V,D}(\omega_{V,N})} \,=\,
-\scr L(V,D)\,\Cal E^+(V,D)\,\frac{L(M,0)}{R_{M,\infty}(\omega_M)}.
$$
 
2) The $p$-adic $L$-function $L_p(T^*(1),N^{\perp},s)$ has a zero of order
$e+r$ where $r=\dim_{\Bbb Q}H^1_f(M)$ at $s=0$ and 
$$
\frac{L_p^*(T^*(1),N^{\perp},0)}{R_{V^*(1),D^{\perp}}(\omega_{V^*(1),N^{\perp}})} \,=\,
\scr L(V,D)\,
\Cal E^+(V^*(1),D^{\perp})\,\frac{L^*(M^*(1),0)}{R_{M^*(1),\infty}(\omega_{M^*(1)})}.
$$
{\,}
\newline
\newline
{\bf Remarks.} 1) If $H^1_f(M)=0$ the $p$-adic regulator vanishes and  we
recover the conjecture formulated in \cite{Ben2}, section 2.3.2.

2) The regulators $R_{M,\infty}(\omega_M)$  and $R_{V,D}(\omega_{V,N^{\perp}})$
are well defined up to the sign and in order to obtain equalities in the formulation
of our conjecture one should make the same choice of signs in the definitions 
of $R_{M,\infty}(\omega_M)$ and $R_{V,D}(\omega_{V,N^{\perp}})$. See \cite{PR2},
section 4.2 for more detail.

3) Our conjecture is compatible with the expected functional equation for
$p$-adic $L$-functions. See section 2.5 of \cite{PR2} and section 5.2.7 below.

\head {\bf \S5. The module of $p$-adic $L$-functions}
\endhead 

\flushpar
{\bf 5.1. The Selmer complex.}
\newline
{\bf 5.1.1. Iwasawa cohomology.} 
Let $\Gamma$ denote the Galois group of $\Bbb Q(\zeta_{p^\infty})/\Bbb Q$ and 
$\Gamma_n=\Gal (\Bbb Q(\zeta_{p^\infty})/\Bbb Q (\zeta_{p^n})).$ Set
$\Lambda=\Bbb Z_p[[\Gamma_1]]$ and 
$\Lambda (\Gamma)= \Bbb Z_p[\Delta]\otimes_{\Bbb Z_p} \Lambda$.
For any character $\eta \in X(\Delta)$ put
$$
e_{\eta}\,=\,\frac{1}{|\Delta|}\sum_{g\in \Delta} \eta^{-1}(g)g.
$$
Then $\Lambda (\Gamma) =\underset{\eta \in X(\Delta)}\to \oplus \Lambda (\Gamma)^{(\eta)}$ where
$\Lambda (\Gamma)^{(\eta)}=\Lambda e_\eta$ and for any $\Lambda (\Gamma)$-module 
$M$ one has a canonical decomposition 
$$
M\simeq \oplus_{\eta \in X(\Delta)}M^{(\eta)},\qquad M^{(\eta)}=e_\eta (M).
$$ 
We write $\eta_0$ for the trivial character of $\Delta$ and 
identify $\Lambda$ with $\Lambda (\Gamma)e_{\eta_0}.$
\newline
\newline
Let $V$ be a p-adic pseudo-geometric representation  
unramified outside $S.$ Set $d(V)=\dim (V)$ and $d_{\pm}(V)=\dim (V^{c=\pm 1}).$
Fix a $\Bbb Z_p$-lattice $T$ of $V$ stable under the action of $G_{S}.$
Let $\iota \,:\,\Lambda (\Gamma)@>>>\Lambda (\Gamma)$ denote the canonical involution $g\mapsto g^{-1}.$ Recall that
the induced module $\text{\rm Ind}_{\Bbb Q(\zeta_{p^\infty})/\Bbb Q} (T)$
is isomorphic to $(\Lambda (\Gamma)\otimes_{\Bbb Z_p}T)^\iota $ (\cite{Ne2}, section 8.1).
Define
$$
\align
&H^i_{\text{\rm Iw},S}(T)\,=\, 
H^i_S((\Lambda (\Gamma)\otimes_{\Bbb Z_p}T)^\iota),\\
&\Hi^i(\Bbb Q_v,T)\,=\,H^i(\Bbb Q_v,(\Lambda (\Gamma)\otimes_{\Bbb Z_p}T)^\iota)
\qquad \text{\rm for any finite place $v$.}
\endalign
$$
From Shapiro's lemma it follows immediately that
$$
H^i_{\text{\rm Iw},S}(T)\,=\,\varprojlim_{\text{\rm cores}} 
H^i_S(\Bbb Q(\zeta_{p^n}),T),\qquad
\Hi^i(\Bbb Q_p,T)\,=\,\varprojlim_{\text{\rm cores}} H^i(\Bbb Q_p(\zeta_{p^n}),T).
$$
Set $H^i_{\text{\rm Iw},S}(V)\,=\,H^i_{\text{\rm Iw},S}(T)\otimes_{\Bbb Z_p}\Bbb Q_p$
and $H^i_{\text{\rm Iw}}(\Bbb Q_v,V)\,=\,
H^i_{\text{\rm Iw}}(\Bbb Q_v,T)\otimes_{\Bbb Z_p}\Bbb Q_p.$ In \cite{PR2}
Perrin-Riou   proved the following results about the structure of these
modules. 
\newline
\,

i)  $H^i_{\text{\rm Iw},S}(V)=0$ and $H^i_{\text{\rm Iw}}(\Bbb Q_v,T)=0$ if $i\ne 1,2;$
\newline
\,

ii) If $v\ne p$, then for each $\eta\in X(\Delta)$ the $\eta$-component
 $H^i_{\text{\rm Iw}}(\Bbb Q_v,T)^{(\eta)}$ is a finitely generated
torsion $\Lambda$-module. In particular, 
 $H^1_{\text{\rm Iw}}(\Bbb Q_v,T)\simeq 
H^1(\Bbb Q_v^{\text{ur}}/\Bbb Q_v,(\Lambda (\Gamma)\otimes_{\Bbb Z_p}T^{I_v})^\iota)$.
\newline
\,

iii) If $v=p$ then $H^2_{\text{\rm Iw}}(\Bbb Q_p,T)^{(\eta)}$ are finitely generated torsion $\Lambda$-modules. 
Moreover, for each $\eta\in X(\Delta)$

$$
\text{\rm {rg}}_{\Lambda} \left (H^1_{\text{\rm Iw}}(\Bbb Q_p,T)^{(\eta)}\right )\,=\,d,\qquad
H^1_{\text{\rm Iw}}(\Bbb Q_p,T)^{(\eta)}_{\text{\rm tor}}\simeq H^0(\Bbb Q_p(\zeta_{p^\infty})\,,T)^{(\eta)}.
$$
Remark that by local duality $H^2_{\text{\rm Iw}}(\Bbb Q_p,T)\simeq H^0(\Bbb Q_p(\zeta_{p^{\infty}}),V^*(1)/T^*(1))$.
\newline
\,

iv) If  the weak Leopoldt conjecture holds for the pair $(V,\eta)$  i.e. if 
$H^2_S(\Bbb Q(\zeta_{p^\infty}),V/T)^{(\eta)}=0$
then $H^2_{\text{\rm Iw},S}(T)^{(\eta)}$ is $\Lambda$-torsion and  
$$
\text{\rm rank}_{\Lambda}\left ( H^1_{\text{\rm Iw},S}(T)^{(\eta)}\right )\,=\,
\cases d_-(V), &{\text{\rm if $\eta (c)=1$}}\\
d_+(V), &{\text{\rm if $\eta (c)=-1$.}}
\endcases
$$ 

Passing to the projective limit in the Poitou-Tate exact sequence one obtains
an exact sequence
$$
\multline
0@>>> H^2_S(\Bbb Q(\zeta_{p^\infty}), V^*(1)/T^*(1))^{\wedge}@>>>
H^1_{\text{\rm Iw},S}(T)@>>> \underset{v\in S}\to \oplus H^1_{\text{\rm Iw}}(\Bbb Q_v,T)@>>>
H^1_S(\Bbb Q(\zeta_{p^\infty}), V^*(1)/T^*(1))^{\wedge}\\
@>>>H^2_{\text{\rm Iw},S}(T)@>>> \underset{v\in S}\to \oplus H^2_{\text{\rm Iw}}(\Bbb Q_v,T)@>>>
H^0_S(\Bbb Q(\zeta_{p^\infty}), V^*(1)/T^*(1))^{\wedge}@>>>
0.
\endmultline
\tag{5.1}
$$
Define
$$
\align
&\RG_{\text{\rm Iw},S} (T)\,=\,C_c^\bullet (G_S,(\Lambda (\Gamma)\otimes_{\Bbb Z_p} T)^\iota),\\
&\RG_{\text{\rm Iw}}(\Bbb Q_v,T)\,=\,C_c^\bullet (G_v,(\Lambda (\Gamma)\otimes_{\Bbb Z_p} T)^\iota),\\
&\RG_S(\Bbb Q(\zeta_{p^\infty}), V^*(1)/T^*(1))\,=\,
C_c^\bullet (G_{S},\text{\rm Hom}_{\Bbb Z_p}(\Lambda (\Gamma), V^*(1)/T^*(1))).
\endalign
$$
Then 
%$$
%\align
%&\bold R^i\Gamma_{\text{\rm Iw},S} (T)\simeq H^i_{\text{\rm Iw},S}(T),\\
%&\bold R^i\Gamma_{\text{\rm Iw}} (\Bbb Q_v,T)\simeq H^i_{\text{\rm Iw}}(\Bbb Q_v,T),\\
%&\bold R^i\Gamma_S(\Bbb Q(\zeta_{p^\infty}), V^*(1)/T^*(1))\,\simeq\,H_S^i(\Bbb Q(\zeta_{p^\infty}), V^*(1)/T^*(1))
%\endalign
%$$
the sequence (5.1) is induced by the distinguished triangle
$$
\RG_{\text{\rm Iw},S} (T)@>>>\underset{v\in S}\to \oplus \RG_{\text{\rm Iw}}(\Bbb Q_v,T)@>>>
\left (\RG_S(\Bbb Q(\zeta_{p^\infty}), V^*(1)/T^*(1))^\iota \right )^{\wedge}\,[-2]
$$
(\cite{Ne2}, Theorem 8.5.6). Finally, we have  usual descent formulas
$$
\RG_{\text{\rm Iw},S} (T)\otimes_{\Lambda}^{\bold L}\Bbb Z_p \simeq \RG_{S} (T),\qquad 
\RG_{\text{\rm Iw}} (\Bbb Q_v,T)\otimes_{\Lambda}^{\bold L}\Bbb Z_p \simeq \RG (\Bbb Q_v,T)
$$
( \cite{Ne2}, Proposition 8.4.21).
\newline
\newline
{\bf 5.1.2. The complex $\bold R\Gamma^{(\eta_0)}_{\text{\rm Iw},h} (D,V)$. } For the remainder of this chapter we assume that $V$ satisfies  the conditions {\bf C1-5)} of section 3.1.2 and that the weak Leopoldt conjecture holds for $(V,\eta_0)$ and $(V^*(1),\eta_0).$
We remark that these assumptions are not independent. Namely, by \cite{PR2}, Proposition B.5 
{\bf C4)} and {\bf C5)} imply the weak Leopoldt conjecture  for $(V^*(1), \eta_0)$.
From the same result it follows that the vanishing of $H^1_f(V^*(1))$ implies the 
weak Leopoldt conjecture for $(V,\eta_0)$ if in addition we assume that $H^0(\Bbb Q_p,V^*(1))=0.$  
\newline
\newline
To simplify notations we write  $\Cal H$ for $\Cal H(\Gamma_1)$.
Fix a regular subspace $D$ of $\Dc (V)$ and  a $\Bbb Z_p$-lattice $N$ of $D$.
Set $\Cal D_p(N,T)^{(\eta_0)}\,=\,N\otimes_{\Bbb Z_p}\Lambda$, 
$\RG_{\text{\rm Iw},f}^{(\eta_0)}(\Bbb Q_p,N,T) =\Cal D_p(N,T)^{(\eta_0)}[-1]$ and
$\RG_{\text{\rm Iw},f}^{(\eta_0)}(\Bbb Q_p,D,V)=\RG_{\text{\rm Iw},f}^{(\eta_0)}(\Bbb Q_p,N,T)\otimes_{\Bbb Z_p}
\Bbb Q_p.$
Consider the map
$$
\bExp_{V,h}^\ep \,:\,\RG_{\text{\rm Iw},f}^{(\eta_0)}(\Bbb Q_p,T)\otimes_{\La}\Cal H
@>>> \RG_{\text{\rm Iw}}^{(\eta_0)} (\Bbb Q_p,T)\otimes_{\Lambda}^{\bold L}\Cal H
$$
which will be viewed as a local condition at $p$.  If  $v\ne p$ the inertia group $I_v$ acts trivially
on $\Lambda$  set
$$
\RG_{\text{\rm Iw},f}^{(\eta_0)}(\Bbb Q_v,N,T)\,=\,
\left [T^{I_v}\otimes \Lambda^{\iota}  @>1-f_v>>T^{I_v}\otimes \Lambda^{\iota}  \right ] 
$$
where the first term is placed in degree $0$. We have a commutative diagram
 $$
\xymatrix{
\RG_{\text{\rm Iw},S}^{(\eta_0)}(T)\otimes_{\Lambda}\Cal H
\ar[r] & \underset{v\in S} \to \oplus \RG_{\text{\rm Iw}}^{(\eta_0)}(\Bbb Q_v,T)\otimes_{\Lambda} \Cal H\\
 & \left (\underset{v\in S}\to \oplus \RG_{\text{\rm Iw},f}^{(\eta_0)}(\Bbb Q_v,N,T)
\right ) \otimes_{\Lambda} \Cal H\,\,.
 \ar[u]
}\tag{5.2}
$$
Consider the associated Selmer complex
$$
\multline
\RG_{\text{\rm Iw},h}^{(\eta_0)}(D,V)\,=\\
\text{\rm cone}\, \left [
\left ( \RG_{\text{\rm Iw},S}^{(\eta_0)}(T)\oplus 
\left (\underset{v\in S}\to\oplus \RG_{\text{\rm Iw},f}^{(\eta_0)}(\Bbb Q_v,N,T)\right ) \right )\otimes_{\Lambda} \Cal H
@>>>
\underset{v\in S}\to \oplus \RG_{\text{\rm Iw}}^{(\eta_0)}(\Bbb Q_v,T)\otimes_{\Lambda} \Cal H\right ]
[-1]
\endmultline
$$
It is easy to see that it does not depend on the choice of $S$. Our main result 
about this complex is the following theorem.

\proclaim{Theorem 5.1.3} Assume that $V$ satisfies the conditions C1-5) and that the
weak Leopoldt conjecture holds for $(V,\eta_0)$ and $(V^*(1),\eta_0).$
Let $D$ be a regular subspace of $\Dc (V).$ Assume that $\scr L(V,D)\ne 0.$  
Then

i) $\bold R^i\Gamma_{\text{\rm Iw},h}^{(\eta_0)}(D,V)$ are $\Cal H$-torsion modules for all $i$.

ii)  $\bold R^i\Gamma_{\text{\rm Iw},h}^{(\eta_0)}(D,V)=0$ for  $i\ne 2,3$ and
$$
\bold R^3\Gamma_{\text{\rm Iw},h}^{(\eta_0)}(D,V)\simeq \left(H^0(\Bbb Q(\zeta_{p^{\infty}}),V^*(1))^*
\right )^{(\eta_0)}\otimes_{\Lambda}\Cal H.
$$

iii) The complex $ \bold R\Gamma_{\text{\rm Iw},h}^{(\eta_0)}(D,V)$ is semisimple i.e.
for each $i$ the natural map 
$$
\bold R^i\Gamma_{\text{\rm Iw},h}^{(\eta_0)}(D,V)^{\Gamma} @>>>
\bold R^i\Gamma_{\text{\rm Iw},h}^{(\eta_0)}(D,V)_{\Gamma}
$$
is an isomorphism.
\endproclaim
\flushpar
{\bf 5.1.4. Proof of Theorem 5.1.3.} We leave the proof
of the following lemma as an easy exercise. 

\proclaim{Lemma 5.1.4.1} Let $A$ and $B$ be two submodules of a finitely generated 
free $\Cal H$-module $M$.
Assume that the  natural maps $A_{\Gamma_1} @>>> M_{\Gamma_1}$ and $B_{\Gamma_1} @>>> M_{\Gamma_1}$ are both
injective. Then $A_{\Gamma_1}\cap B_{\Gamma_1}=\{0\}$ implies that $A\cap B=\{0\}.$
\endproclaim
\flushpar
{\bf 5.1.4.2.} Since $H_{\text{\rm Iw},S}^0(V)$ and
$\Hi^0(\Bbb Q_v,V)$ are zero, we have   $ \bold R^0\Gamma_{\text{\rm Iw},h}^{(\eta_0)}(D,V)=0.$
Next, by definition $\bold R^1\Gamma_{\text{\rm Iw},h}^{(\eta_0)}(D,V)\,=\,\ker (f)$ where
$$
f\,:\, \left (
H^1_{\text{\rm Iw},S}(T)^{(\eta_0)}\oplus 
\Cal D_p(N,T)^{(\eta_0)}
\underset{v\in S-\{p\}}\to \oplus H^1_{\text{\rm Iw},f}(\Bbb Q_v,T)^{(\eta_0)}
\right )\otimes \Cal H @>>>
\underset{v\in S}\to \oplus \Hi^1(\Bbb Q_v,T)^{(\eta_0)}\otimes \Cal H 
$$ 
is the map induced by (5.2).
If $v\in S-\{p\}$ one has   
$$
H^1_{\text{\rm Iw},f}(\Bbb Q_v,T)^{(\eta_0)}\,=\,\Hi^1(\Bbb Q_v,T)^{(\eta_0)}= 
H^1(\Bbb Q_v^{\text{\rm ur}}/\Bbb Q_v, (\Lambda\otimes T^{I_v})^{\iota}).
$$  
Thus
$$
\bold R^1\Gamma_{\text{\rm Iw},h}^{(\eta_0)}(D,V)\,=\,
\left (H^1_{\text{\rm Iw},S}(T)^{(\eta_0)}\otimes_{\Lambda}\Cal H \right )  \cap \left (\Exp^\ep_{V,h}
\left (\Cal D_p(D,T)^{(\eta_0)}\right  )\otimes_{\Lambda} \Cal H\right ) 
$$ 
in  $\Hi^1(\Bbb Q_p,T)^{(\eta_0)}\otimes_{\Lambda}\Cal H$. Put
$$
A=\Exp_{V,h}^\ep (D_{-1}\otimes \Cal H)\oplus X^{-1} \Exp_{V,h}^\ep (D^{\Ph=p^{-1}}\otimes \Cal H)
\subset \Hi^1(\Bbb Q_p,T)^{(\eta_0)}\otimes_{\Lambda}\Cal H.
$$
By Theorem 2.2.4 and Proposition 3.2.2 $A_{\Gamma_1}$ injects into 
%$\Hi^1(\Bbb Q_p,V)^{(\eta_0)}_{\Gamma_1}
%\hookrightarrow 
$H^1(\Bbb Q_p,V).$ The $\Cal H$-module  $M=\left (\dsize \frac{\Hi^1(\Bbb Q_p,T)}{T^{H_{\Bbb Q_p}}}\right )^{(\eta_0)}\otimes_{\Lambda}\Cal H $ is  free and $A\hookrightarrow M.$ Since  
$T^{G_{\Bbb Q_p}}=0$  one has
$
M_{\Gamma_1}\,=\,\Hi^1(\Bbb Q_p,V)_{\Gamma}\subset  H^1(\Bbb Q_p,V)
$
and we obtain that  $A_{\Gamma_1}$ injects into $M_{\Gamma_1}.$

Set $B=\left (\dsize \frac{H^1_{\text{\rm Iw},S}(T)}{T^{H_{\Bbb Q}}} \right )^{(\eta_0)}\otimes_{\Lambda}\Cal H.$ 
The weak Leopoldt conjecture for $(V^*(1),\eta_0)$ together with
the fact that $\Hi^1(\Bbb Q_v,T)$ are $\Lambda$-torsion for $v\in S-\{p\}$
imply that  $B\hookrightarrow M.$  Since the image of $\Hi^1(\Bbb Q_v,V)_{\Gamma}$ in
$H^1(\Bbb Q_v,V)$ is contained in $H^1_f(\Bbb Q_v,V),$ the image of $H^1_{\text{\rm Iw},S}(V)_{\Gamma}$
in $H^1_S(V)$ is in fact contained in
$H^1_{f,\{p\}}(V).$
From {\bf C5)} it follows that  $H^1_{f,\{p\}}(V)$ injects into $H^1(\Bbb Q_p,V)$ and we have 
$$
H^1_{\text{\rm Iw},S}(V)^{(\eta_0)}_{\Gamma_1}=H^1_{\text{\rm Iw},S}(V)_{\Gamma}\hookrightarrow H^1_{f,\{p\}}(V)\hookrightarrow
H^1(\Bbb Q_p,V).
$$ 
Thus $B_{\Gamma_1}\subset M_{\Gamma_1}.$ We shall prove that $\bold R^1\Gamma_{\text{\rm Iw},h}^{(\eta_0)}(D,V)=0.$
By Lemma 5.1.4.1 it suffices to show that $A_{\Gamma_1} \cap B_{\Gamma_1}=\{0\}.$
Now we  claim that $ A_{\Gamma_1}\cap H^1_{f,\{p\}}(V)=\{0\}.$ 
First note that by Lemma 3.1.4
$$
 H^1_{f,\{p\}}(V)\hookrightarrow \frac{H^1(\Bbb Q_p,V)}{H^1(\bD_{-1})}.
$$ 
On the other hand, from Theorem 2.2.4 
it follows that 
$$
\Exp_{V,h}^\ep(D_{-1}\otimes \Cal H)_{\Gamma_1}=\exp_{V,\Bbb Q_p} (D_{-1})\subset H^1(\bD_{-1}).
$$
Now Proposition 3.2.2 implies that the image of $A_{\Gamma_1}$ in   
$\dsize \frac{H^1(\Bbb Q_p,V)}{H^1(\bD_{-1})}$  coincides with
$H^1_c(W).$ But $\scr L(V,D)\ne 0$
if and only if 
$H^1_D(V)\cap H^1_c(W) =0
$ 
where
$H^1_D(V)$ denotes  the inverse image of $H^1(W)$ in $H^1_{f,\{p\}}(V)$
(see Lemma 3.1.4 iii)).
This proves the claim and implies that $\bold R^1\Gamma_{\text{\rm Iw},h}^{(\eta_0)}(D,V)=0.$
\newline
\newline
{\bf 5.1.4.3.} We shall show that $\bold R^2\Gamma_{\text{\rm Iw},h}^{(\eta_0)}(D,V)$ is $\Cal H$-torsion. 
By definition, we have an exact sequence  
$$
0@>>> \text{\rm coker} (f)@>>>\bold R^2\Gamma_{\text{\rm Iw},h}^{(\eta_0)}(D,V)@>>>\sha^2_{\,\text{\rm Iw},S}(V)^{(\eta_0)}\otimes_{\Lambda_{\Bbb Q_p}}
\Cal H@>>>0,
\tag{5.3}
$$ 
where
$$
\sha^2_{\,\text{\rm Iw},S}(V)\,=\,\ker \left ( H^2_{\text{\rm Iw},S}(V)
@>>>
\underset{v\in S}\to \oplus \Hi^2(\Bbb Q_v,V) \right ).
$$
It follows from the weak Leopoldt conjecture that $\sha^2_{\,\text{\rm Iw},S}(V)$ is $\Lambda_{\Bbb Q_p}$-torsion.
On the other hand, as $\Cal H$ is a Bezout ring \cite{La}, the formulas
$$
\text{\rm rank}_{\Lambda} H^1_{\text{\rm Iw},S} (T)^{(\eta_0)}=d_{-}(V),\quad
\text{\rm rank}_{\Lambda} H^1_{\text{\rm Iw}} (\Bbb Q_p,T)^{(\eta_0)}=d(V),\quad
\text{\rm rank}_{\Lambda} \Cal D_p(N,T)=d_{+}(V)
$$
together with the fact that $\bold R^1\Gamma_{\text{\rm Iw},h}^{(\eta_0)}(D,V)=0$ 
imply that $\text{\rm coker} (f)$ is $\Cal H$-torsion. We have therefore proved that
$\bold R^2\Gamma_{\text{\rm Iw},h}(D,V)$ is $\Cal H$-torsion. Finally,
 the Poitou-Tate exact sequence  gives that 
$$
\bold R^3\Gamma_{\text{\rm Iw},h}^{(\eta_0)}(D,V)= \left (H^0(\Bbb Q(\zeta_{p^\infty}),V^*(1))^*\right )^{(\eta_0)}
\otimes_{\Lambda_{\Bbb Q_p}}\Cal H
$$
is also $\Cal H$-torsion. 
\newline
\newline
{\bf 5.1.4.4.} Now we prove the semisimplicity of 
$\bold R \Gamma_{\text{\rm Iw},h}^{(\eta_0)}(D,V).$
First write $H^1_{\text{\rm Iw},S}(V)^{(\eta_0)}\simeq \Lambda_{\Bbb Q_p}^{d_{-}(V)}\oplus 
H^1_{\text{\rm Iw},S}(V)^{(\eta_0)}_{\text{\rm tor}}. $ 
Since  $H^1_{\text{\rm Iw},S} (V)_{\text{\rm tor}}
\subset H^1_{\text{\rm Iw}} (\Bbb Q_p,V)_{\text{\rm tor}}=V^{H_{\Bbb Q_p}},$ we have  
$(H^1_{\text{\rm Iw},S} (V)_{\text{\rm tor}})_{\Gamma}=0$ by the snake
lemma. Thus $\dim_{\Bbb Q_p} H^1_{\text{\rm Iw},S} (V)^{(\eta_0)}_{\Gamma_1}=d_{-}(V).$
On the other hand 
$\dim_{\Bbb Q_p} H^1_{f,\{p\}}(V)=d_{-}(V)+\dim_{\Bbb Q_p}H^0(\Bbb Q_p,V^*(1))$
by (3.2) and the dimension argument shows that in the commutative diagram
$$
\CD
0@>>>H^1_{\text{\rm Iw},S}(V)^{(\eta_0)}_{\Gamma_1}@>>>H^1_{f,\{p\}}(V)
@>>>H^0(\Bbb Q_p,V^*(1))^*@>>>0
\\
@. @VVV @VVV @VVV @.\\
0@>>>H^1_{\text{\rm Iw}} (\Bbb Q_p,V)^{(\eta_0)}_{\Gamma_1}@>>>H^1 (\Bbb Q_p,V) @>>>
H^0(\Bbb Q_p,V^*(1))^*@>>>0
\endCD
\tag{5.4}
$$
with obviously exact upper line the bottom line is also exact.
This  implies immediately that the natural map
$$
\frac{H^1_{\text{\rm Iw}} (\Bbb Q_p,V)^{(\eta_0)}_{\Gamma_1}}
{H^1_{\text{\rm Iw},S}(V)^{(\eta_0)}_{\Gamma_1}+H^1(\bD_{-1})} 
@>>>
\frac{H^1(\Bbb Q_p,V)}{H^1_{f,\{p\}}(V)+H^1(\bD_{-1})}
$$
is an isomorphism.

Consider the exact sequence
$$
0@>>>
\left (
H^1_{\text{\rm Iw},S}(T)^{(\eta_0)}\oplus 
\Cal D_p(N,T)^{(\eta_0)}
\right )\otimes \Cal H @>>>
\Hi^1(\Bbb Q_p,T)^{(\eta_0)}\otimes \Cal H @>>>\text{\rm coker} (f)@>>>0.
$$ 
Recall that $\Exp_{V,h,0}^\ep\,:\,D @>>>\Hi^1 (\Bbb Q_p,V)_{\Gamma}$ denotes the homomorphism
induced by the large exponential map. 
Applying the snake lemma, and 
taking into account that $\text{\rm Im} (\Exp_{V,h,0}^\ep)=
\exp_{V,\Bbb Q_p} (D_{-1})=H^1(\bD_{-1})$ and  
$\ker (\Exp_{V,h,0}^\ep)= D^{\Ph=p^{-1}}$ (see for example \cite{BB}, Propositions 4.17 and 4.18
or the proof of Proposition 3.3.2) we obtain
$$
\align &\text{\rm coker} (f)^{\Gamma_1}=\ker \left ( 
H^1_{\text{\rm Iw},S}(V)^{(\eta_0)}_{\Gamma_1}\oplus D@>\Exp^\ep_{V,h,0}>> H^1(\Bbb Q_p,V) \right )=D^{\Ph=p^{-1}}, \qquad \text{\rm (by the regularity of $D$)},
\\
&\text{\rm coker} (f)_{\Gamma_1}=\frac{H^1_{\text{\rm Iw}}(\Bbb Q_p,V)^{(\eta_0)}_{\Gamma_1}}{H^1_{\text{\rm Iw},S}(V)^{(\eta_0)}_{\Gamma_1}+H^1(\bD_{-1})} = 
\frac{H^1(\Bbb Q_p,V)}{H^1_{f,\{p\}}(V)+H^1(\bD_{-1})}
.
\endalign
$$
Thus on has a commutative diagram
$$
\CD
\text{\rm coker} (f)^{\Gamma_1} @>>> D^{\Ph=p^{-1}}\\
@VVV   @VV{\delta_{D,h}}V \\
\text{\rm coker} (f)_{\Gamma_1}  @>>> \dsize\frac{H^1(\Bbb Q_p,V)}{H^1_{f,\{p\}}(V)+H^1(\bD_{-1})}.
\endCD
\tag{5.5}
$$
where horizontal arrows are isomorphisms, the left vertical arrow is the natural projection
and the right vertical row is the map defined in section 3.2.3.
From Proposition 3.2.4 it follows that  $\text{\rm coker} (f)^{\Gamma_1}@>>>
\text{\rm coker} (f)_{\Gamma_1}$ is an isomorphism if and only if $\scr L(V,D)\ne 0.$

On the other hand, the arguments \cite{PR2}, section 3.3.4 show that 
$\sha^2_{\,\text{\rm Iw},S}(V)_{\Gamma}=\sha^2_{\,\text{\rm Iw},S}(V)^{\Gamma}=0.$
Remark that Perrin-Riou assumes that $\Dc (V)^{\Ph=1}=\Dc (V)^{\Ph=p^{-1}}=0$, but
her proof  works in our case without modifications and we repeat it for the
commodity of the reader. Consider the commutative diagram (where 
we write $\sha^2_{\,\text{\rm Iw}}(V)$ instead  $\sha^2_{\,\text{\rm Iw},S}(V)$ and
$H^1(\Bbb Q(\zeta_{p^\infty}),V^*(1))^*_{\Gamma}$ instead 
$(H^1_S(\Bbb Q(\zeta_{p^\infty}),V^*(1))^*)_{\Gamma}$ to abbraviate notation)
\eightpoint
$$
\xymatrix{&0 \ar[d]& 0 \ar[d] & & &\\
0 \ar[r] &H^1_{\text{\rm Iw}}(V)_{\Gamma} \ar[r] \ar[d]
&\underset{v\in S}\to \oplus \Hi^1(\Bbb Q_v,V)_{\Gamma} \ar[r] \ar[d] 
&H^1(\Bbb Q(\zeta_{p^\infty}),V^*(1))^*_{\Gamma} \ar[r] \ar[d]^{=}
&\sha^2_{\text{\rm Iw}}(V)_{\Gamma} \ar[r] &0\\
0 \ar[r] &H^1_{f,\{p\}}(V) \ar[r] \ar[d]
&\underset{\Sb v\in S\\v\ne p\endSb}\to \oplus H^1_f (\Bbb Q_v,V)\oplus H^1(\Bbb Q_p,V) \ar[r] \ar[d]
&H^1(V^*(1))^* \ar[r]  &0 &\\
&H^0(\Bbb Q_p,V^*(1))^* \ar[d]\ar[r]^{=}&H^0(\Bbb Q_p,V^*(1))^* \ar[d]  & & &\\
&0&0&&&.
}
$$
%$$
%\xymatrix{&0 \ar[d]& 0 \ar[d] & & &\\
%0 \ar[r] &H^1_{\text{\rm Iw},S}(V)_{\Gamma} \ar[r] \ar[d]
%&\underset{v\in S}\to \oplus \Hi^1(\Bbb Q_v,V)_{\Gamma} \ar[r] \ar[d] 
%&\left (H^1_S(\Bbb Q(\zeta_{p^\infty}),V^*(1))^*\right )_{\Gamma} \ar[r] \ar[d]^{=}
%&\sha^2_{\,\text{\rm Iw},S}(V)_{\Gamma} \ar[r] &0\\
%0 \ar[r] &H^1_{f,\{p\}}(V) \ar[r] \ar[d]
%&\underset{v\in S\setminus\{p\}}\to \oplus H^1_f (\Bbb Q_v,V)\oplus H^1(\Bbb Q_p,V) \ar[r] \ar[d]
%&H^1_S(V^*(1))^* \ar[r]  &0 &\\
%&H^0(\Bbb Q_p,V^*(1))^* \ar[d]\ar[r]^{=}&H^0(\Bbb Q_p,V^*(1))^* \ar[d]  & & &\\
%&0&0&&&.
%}
%$$
\tenpoint
The top row of this diagram is obtained by taking coinvariants in the Poitou-Tate
exact sequence. Thus it is  exact. The middle row is obtained from  the exact sequence
$$
0@>>> H^1_S(V^*(1))@>>>H^1(\Bbb Q_p,V^*(1))\oplus \underset{v\in S-\{p\}}\to \bigoplus
\dsize\frac{H^1(\Bbb Q_v,V^*(1))}{H^1_f(\Bbb Q_v,V^*(1))} @>>>H^1_{f,\{p\}}(V)^*@>>>0
$$
by taking duals. Here we use the condition $H^1_f(V^*(1))=0$.  
The exactness of the left and middle columns  follows from the diagram (5.4).
The isomorphism from  the right column  comes from the exact sequence
$$
0@>>>H^1(\Gamma, H^0_S( \Bbb Q(\zeta_{p^\infty}),V^*(1)))@>>>H^1_S(V^*(1))
@>>>H^1_{S}(\Bbb Q(\zeta_{p^\infty}),V^*(1))^{\Gamma}@>>>0
$$
together with the remark that $H^1(\Gamma, H^0_S( \Bbb Q(\zeta_{p^\infty}),V^*(1)))=0$
because $H^0(\Gamma, H^0_S( \Bbb Q(\zeta_{p^\infty}),V^*(1)))= 
H^0_S( \Bbb Q,V^*(1)))=0$ by {\bf C2)}.
Now an easy diagram search shows  that $\sha^2_{\,\text{\rm Iw},S}(V)_{\Gamma}=0.$
Finally, from  $\text{\rm dim}_{\Bbb Q_p} \sha^2_{\,\text{\rm Iw},S}(V)^{\Gamma}
\leqslant \text{\rm dim}_{\Bbb Q_p} \sha^2_{\,\text{\rm Iw},S}(V)_{\Gamma}$ 
it follows that $\sha^2_{\,\text{\rm Iw},S}(V)^{\Gamma}=0.$
Therefore, applying the snake lemma to (5.3) we obtain a commutative diagram
$$
\xymatrix{
\text{\rm coker} (f)^{\Gamma_1} \ar[r] \ar[d] &\bold R^2\Gamma^{(\eta_0)}_{\text{\rm Iw},h}(D,V)^{\Gamma_1}\ar[d]\\
\text{\rm coker} (f)_{\Gamma_1} \ar[r] &\bold R^2\Gamma^{(\eta_0)}_{\text{\rm Iw},h}(D,V)_{\Gamma_1},}
$$
in which  the horizontal arrows are isomorphisms and the vertical arrows are natural projections.
This proves that $\bold R\Gamma^{(\eta_0)}_{\text{\rm Iw},h}(D,V)$ is semisimple
in degree $2$. Remark that the semisimplicity in degree $3$  is obvious because by ii)
$
\bold R^3\Gamma^{(\eta_0)}_{\text{\rm Iw},h}(D,V)^{\Gamma_1}=
\bold R^3\Gamma^{(\eta_0)}_{\text{\rm Iw},h}(D,V)_{\Gamma_1}=0.
$
 This completes the proof of Theorem 5.1.3.

\proclaim{Corollary 5.1.5} The exponential map induces an isomorphism
of $D^{\Ph=p^{-1}}$ onto $\text{\rm coker} (f)_{\Gamma_1} \simeq \bold R^2\Gamma^{(\eta_0)}_{\text{\rm Iw},h}(D,V)_{\Gamma_1}$
and the diagram 
$$
\xymatrix{
D^{\Ph=p^{-1}} \ar[rr]^{\sim}   \ar[d]^{\lambda_D}& &\bold R^2\Gamma^{(\eta_0)}_{\text{\rm Iw},h}(D,V)^{\Gamma_1}
\ar[d]\\
D^{\Ph=p^{-1}} \ar[rr]^-{(h-1)!\exp_V} & &\bold R^2\Gamma^{(\eta_0)}_{\text{\rm Iw},h}(D,V)_{\Gamma_1}}
$$
in which the map $\lambda_D$ is defined in Proposition 3.2.4, commutes.
\endproclaim

$\,$
\flushpar
{\bf 5.2. The module of $p$-adic $L$-functions.}
\newline
{\bf 5.2.1. The canonical trivialisation.} We conserve the notation and conventions of section 4.2. 
Let $D$ be an admissible subspace of $\Dc (V)$ and assume that $\scr L(V,D)\ne 0.$
We review Perrin-Riou's definition of the module of $p$-adic $L$-functions using the formalism
of Selmer complexes.
Set 
$$
\Delta_{\text{\rm Iw},h} (D,V)\,=\,
{\det}^{-1}_{\Lambda_{\Bbb Q_p}}\left (\RG_{\text{\rm Iw},S}^{(\eta_0)}(V) \oplus \left (\underset{v\in S}\to \oplus
\RG_{\text{\rm Iw},f}^{(\eta_0)}(\Bbb Q_v,D,V)\right )\right ) \otimes 
{\det}_{\Lambda_{\Bbb Q_p}} \left (\underset{v\in S}\to \oplus
\RG_{\text{\rm Iw}}^{(\eta_0)}(\Bbb Q_v,V)\right ).
$$
The exact triangle
$$
\bold R\Gamma_{\text{\rm Iw},S}^{(\eta_0)}(D,V)@>>>
\left (\RG_{\text{\rm Iw},S}^{(\eta_0)}(V) \oplus \left (\underset{v\in S}\to \oplus
\RG_{\text{\rm Iw},f}^{(\eta_0)}(\Bbb Q_v,D,V)\right )\right )\otimes{\Cal H} @>>>
\left (\underset{v\in S}\to \oplus
\RG_{\text{\rm Iw}}^{(\eta_0)}(\Bbb Q_v,V)\right )\otimes \Cal H
$$
gives an isomorphism
$
\Delta_{\text{\rm Iw},h} (D,V)\otimes_{\Lambda_{\Bbb Q_p}}\Cal H \,\simeq\,
{\det}_{\Cal H}^{-1}\bold R\Gamma_{\text{\rm Iw},S}^{(\eta_0)}(D,V).
$
Let $\Cal K$ denote the  field of fractions of $\Cal H.$
By Theorem 5.1.3, all  $\bold R^i\Gamma_{\text{\rm Iw},S}^{(\eta_0)}(D,V)$
are $\Cal H$-torsion and we have a canonical map. 
$$
{\det}_{\Cal H}^{-1}\bold R\Gamma_{\text{\rm Iw},S}^{(\eta_0)}(D,V) \simeq
\underset{i\in\{2,3\}}\to \otimes {\det}^{(-1)^{i+1}}_{\Cal H}\bold R^i\Gamma_{\text{\rm Iw},S}^{(\eta_0)}(D,V)
\hookrightarrow \Cal K.
$$
The composition of these maps  gives a trivialization
$
i_{V,\text{\rm Iw},h}\,\,:\,\, 
\Delta_{\text{\rm Iw},h} (D,V) @>>> \Cal K.
$
\newline
\newline
{\bf 5.2.2. Local conditions.} In this section we compare local conditions
coming from Perrin-Riou's theory to the Bloch-Kato's one.  
Set $\RG_{f} (\Bbb Q_p,D,V)=D[-1]$ and define
$$
S=\text{\rm cone} \left (
\frac{1-p^{-1}\Ph^{-1}}{1-\Ph}\,:\, \RG_{f} (\Bbb Q_v,D,V)@>>>\RG_{f} (\Bbb Q_p,V)\right )\,[-1].
\tag{5.6}
$$
Thus, explicitly
$$
S=\left [ D\oplus \Dc (V)@>>>\Dc (V)\oplus t_V(\Bbb Q_p) \right ]\,[-1]\simeq 
\left [ D\oplus \Dc (V)@>>>\Dc (V)\oplus D) \right ]\,[-1],
$$
where the unique non-trivial map is given by
$$
(x,y)\mapsto \left ( (1-\Ph)\,y,\, \left (\frac{1-p^{-1}\Ph^{-1}}{1-\Ph}\,x+y\right )\,\pmod{\F^0\Dc (V)}
\right ).
$$
Thus
$
H^1(S)\,=\,D^{\Ph=p^{-1}}$ and $H^2(S)\,=\,\dsize\frac{t_V(\Bbb Q_p)}{(1-p^{-1}\Ph^{-1})D}\simeq
\dsize\frac{\Dc (V)}{\F^0\Dc (V)+D_{-1}}$.
From the semi-simplicity of $\dsize\frac{1-p^{-1}\Ph^{-1}}{1-\Ph}$ it follows
that the natural projection  $H^1( S) \oplus H^1_f(V) @>>>H^2( S)$
is an  isomorphism and we have a  canonical trivialization
$$
\alpha_S\,:\,{\det}_{\Bbb Q_p}  S \otimes {\det}_{\Bbb Q_p}\RG_f(V) \simeq {\det}_{\Bbb Q_p}^{-1} H^1( S) \otimes 
{\det}_{\Bbb Q_p} H^2(S)\otimes {\det}_{\Bbb Q_p}^{-1} H^1_f(V)\simeq \Bbb Q_p.
\tag{5.7}
$$
Hence the distingushed triangle 
$$
 S@>>> \RG_{f} (\Bbb Q_p,D,V) @>>>\RG_{f} (\Bbb Q_p,V)@>>>  S[1]
$$
 induces   isomorphisms
$$
\multline
\beta_S\,:\,{\det}_{\Bbb Q_p} t_V(\Bbb Q_p)\otimes {\det}_{\Bbb Q_p}\RG_f(V) \simeq 
{\det}_{\Bbb Q_p}^{-1} \RG_{f} (\Bbb Q_p,V)\otimes {\det}_{\Bbb Q_p}\RG_f(V)\\
\simeq {\det}_{\Bbb Q_p}^{-1} \RG_{f} (\Bbb Q_p,D,V) \otimes
{\det}_{\Bbb Q_p}   S \otimes {\det}_{\Bbb Q_p}\RG_f(V)\\
\simeq
{\det}_{\Bbb Q_p}D \otimes
{\det}_{\Bbb Q_p}   S \otimes {\det}_{\Bbb Q_p}\RG_f(V)
\endmultline
\tag{5.8}
$$ 
and
$$
\vartheta_S\,:\,
{\det}_{\Bbb Q_p} t_V(\Bbb Q_p)\otimes {\det}_{\Bbb Q_p}\RG_f(V) \overset{\beta_S}\to\simeq 
{\det}_{\Bbb Q_p}D \otimes
{\det}_{\Bbb Q_p}   S \otimes {\det}_{\Bbb Q_p}\RG_f(V) \overset{\text{\rm id}\otimes \alpha_S}\to\simeq 
{\det}_{\Bbb Q_p}D.
\tag{5.9}
$$
Fix bases $\omega_{t_V}\in \det_{\Bbb Q_p}t_V(\Bbb Q_p),$ $\omega_D\in 
\det_{\Bbb Q_p}D$ and $\omega_f\in \det_{\Bbb Q_p}H^1_f(V).$
% Write 
%$\omega_D=\omega_{-1}\otimes \omega_{0}$ with $\omega_{-1}\in \det_{\Bbb Q_p}D_{-1}$
%and $\omega_0\in \det_{\Bbb Q_p}D^{\Ph=p^{-1}}.$ 
Let $R_{V,D}(\omega_{V,D})$
denote the determinant of the regulator map
$$
r_{V,D}\,:\,H^1_f(V)@>>>\Dc(V)/(\F^0\Dc (V)+D)
$$
with respect to $\omega_f$ and $\omega_{t_V}\otimes \omega_D^{-1}.$

\proclaim{Lemma 5.2.3} i) Let $f\,:\,W@>>>W$ be a semi-simple endomorphism
of a finitely dimensional $k$-vector space $W$. The canonical
projection $\ker (f)@>>>\text{\rm coker} (f)$ is an isomorphism and 
the tautological
exact sequence 
$$
0@>>>\ker (f)@>>>W@>f>>W@>>>\text{\rm coker} (f)@>>>0
$$
induces an isomorphism
$$
{\det}^*f\,:\, {\det}_k (W) @>>> {\det}_k (W)\otimes  {\det}_k (\ker (f))\otimes
 {\det}^{-1}_k (\text{\rm coker} (f)) @>>>{\det}_k (W).
$$
Then
$
{\det}^*f (x)={\det} (f\,\vert\, \text{\rm coker} (f)).
$

ii) The map $\vartheta_S$ sends $\omega_{t_V}\otimes \omega_f^{-1}$ onto
$$
{\det}^*\left (\frac{1-p^{-1}\Ph^{-1}}{1-\Ph} \vert D \right )^{-1} 
  E_p(V,1)^{-1} R_{V,D}(\omega_{V,D})^{-1} \omega_D
$$
\endproclaim
\demo{Proof} The proof  is straightforward and is omitted here.
\enddemo
{\,}
\flushpar
{\bf 5.2.4. Definition of the module of $p$-adic $L$-functions.}
In this subsection we interpret Perrin-Riou's construction 
of the module of $p$-adic $L$-functions in terms of \cite{Ne2}.
Fix a $\Bbb Z_p$-lattice $N$ of  $D$ and set
$$
\Delta_{\text{\rm Iw},h} (N,T)\,=\,
{\det}^{-1}_{\Lambda}\left (\RG_{\text{\rm Iw},S}^{(\eta_0)}(T) \oplus \left (\underset{v\in S}\to \oplus
\RG_{\text{\rm Iw},f}^{(\eta_0)}(\Bbb Q_v,N,T)\right )\right ) \otimes 
{\det}_{\Lambda} \left (\underset{v\in S}\to \oplus
\RG_{\text{\rm Iw}}^{(\eta_0)}(\Bbb Q_v,T)\right ).
$$
The module of $p$-adic $L$-functions associated to $(N,T)$ is defined  as
$$
\bold L_{\text{\rm Iw},h}^{(\eta_0)}(N,T)  =
i_{V,\text{\rm Iw},h} \left (\Delta_{\text{\rm Iw},h} (N,T)\right ) \subset \Cal K.
$$
Fix a generator $f(\gamma_1-1)$ of $\bold L_{\text{\rm Iw},h}^{(\eta_0)}(N,T)$ and define
a meromorphic $p$-adic function
$$
L_{\text{\rm Iw},h}(T,N,s)=f(\chi (\gamma)^s-1).
$$
{\,}
\newline
\newline
Let now  $V$ be the $p$-adic realisation of a pure motive $M$ over $\Bbb Q$
which satisfies the conditions {\bf M1-3)} of section 4.1.2.  As we saw in section 4.1.2
on expects that $V$ satisfies {\bf C1-5)}.
We fix 
bases $\omega_f\in {\det}_{\Bbb Q}H^1_f(M),$ $\omega_{t_M}\in 
{\det}_{\Bbb Q}t_M(\Bbb Q)$ and use the same notation for their images in
${\det}_{\Bbb Q_p}H^1_f(V)$ and ${\det}_{\Bbb Q}t_V(\Bbb Q_p)$ respectively.
Choose bases  $\omega_{M_{\text{\rm B}}}^+\in {\det}_{\Bbb Q}M_{\text{\rm B}}^+$ and 
$\omega_T^+\in {\det}_{\Bbb Z_p}T^+$  and define the $p$-adic period 
$\Omega_M^{(\acute et,p)} (\omega_T^+,\omega_{M_{\text{\rm B}}}^+)\in \Bbb Q_p$ 
by 
$
\omega_{T}^+=\Omega_M^{(\acute et,p)} (\omega_T^+,\omega_{M_{\text{\rm B}}}^+)
$
using the comparision isomorphism (4.4) and (4.7).
Let $\omega_N$ be a generator  of ${\det}_{\Bbb Z_p}N.$ 

\proclaim{Theorem 5.2.5} Assume that  $V$ satisfies {\bf C1-5)} and that the weak
Leopoldt conjecture holds for $(V,\eta_0)$ and $(V^*(1), \eta_0).$ 
Let $D$ be an admissible subspace of $\Dc (V).$  Assume that $\scr L(V,D)\ne 0.$ Then

i)  $L_{\text{\rm Iw},h}(T,N,s)$ is a meromorphic $p$-adic function which has a zero 
at $s=0$ of order $e=\dim_{\Bbb Q_p}(D^{\Ph=p^{-1}}).$

ii) Let  $L_{\text{\rm Iw},h}^*(T,N,0)= \lim_{s\to 0} s^{-e}L_{\text{\rm Iw},h}(T,N,s)$
be the special value of $L_{\text{\rm Iw},h}(T,N,s)$ at $s=0.$ Then
$$
\frac{L_{\text{\rm Iw},h}^*(T,N,0)}{R_{V,D}(\omega_{V,N})}
\sim_p \Gamma (h)^{d_{+}(V)} \,\scr L (V,D)\, 
\Cal E^+(V,D)
\,\,\frac{i_{\omega_M, p}\,(\Delta_{\text{\rm EP}} (T))}{\Omega_M^{(\acute et,p)} (\omega_T^+,\omega_{M_{\text{\rm B}}}^+)},
$$
where $i_{\omega_M,p}$ and $\Cal E^+(V,D)$ are defined by (4.12) and (4.16) respectively
and $\Gamma (h)=(h-1)!.$
\endproclaim
\flushpar
{\bf 5.2.6. Proof of Theorem 5.2.5.} 
\newline
{\bf 5.2.6.1.} First recall  the formalism of Iwasawa descent which will be used in the proof.
The result we need is proved in \cite{BG}. This is a particular case of Nekov\'a\v r's
descent theory \cite{Ne2}. 
Let $C^\bullet$ be a perfect complex of $\Cal H$-modules and let $C^\bullet_0=C^\bullet\otimes^{\bold L}_{\Cal H}\Bbb Q_p.$
We have a natural distinguished triangle 
$$
C^\bullet @>X>>C^\bullet @>>>C^\bullet_0,
$$
where $X=\gamma_1-1.$ In each degree this triangle gives a short exact sequence
$$
0@>>>H^n(C^\bullet)_{\Gamma_1}@>>>H^n(C^\bullet_0)@>>>H^{n+1}(C^\bullet)^{\Gamma_1}@>>>0.
$$
One says that $C^\bullet$ is semisimple if the natural map
$$
H^n(C^\bullet)^{\Gamma_1}@>>>H^n(C^\bullet)@>>>H^n(C^\bullet)_{\Gamma_1}  \tag{5.10}
$$
is an isomorphism in all degrees. If $C^\bullet$ is semisimple, there exists a
natural trivialisation of ${\det}_{\Bbb Q_p}C^\bullet_0$, namely 
$$
\multline
\vartheta\,\,:\,\,{\det}_{\Bbb Q_p}C^\bullet_0 \simeq \underset{n\in \Bbb Z}\to\otimes
{\det}_{\Bbb Q_p}^{(-1)^n}H^n(C_0) \simeq 
\underset{n\in \Bbb Z}\to\otimes \left ({\det}^{(-1)^n}_{\Bbb Q_p}H^n(C^\bullet)_{\Gamma_1} \otimes 
{\det}^{(-1)^n}_{\Bbb Q_p}H^{n+1}(C^\bullet)^{\Gamma_1} \right ) \\
\simeq\underset{n\in \Bbb Z}\to\otimes \left ({\det}^{(-1)^n}_{\Bbb Q_p}H^n(C^\bullet)_{\Gamma_1} \otimes 
{\det}^{(-1)^{n-1}}_{\Bbb Q_p}H^{n}(C^\bullet)^{\Gamma_1} \right )
\simeq \Bbb Q_p
\endmultline
$$
where the last map is induced by (5.10). 
We now suppose that $C\otimes_{\Cal H}\Cal K$ is acyclic and write 
$
i_\infty \,\,:\,\, {\det}_{\Cal H} C^\bullet @>>> \Cal K
$
for the associated morphism in $\Cal P (\Cal K).$ Then
$
i_\infty ({\det}_{\Cal H} C^\bullet)=f\Cal H,
$
where $f\in \Cal K.$  Let $r$ be the unique integer such that 
$X^{-r}f$ is a unit of the localization $\Cal H_0$ of $\Cal H$ 
with respect to the principal ideal $X\Cal H$. 

\proclaim{Lemma 5.2.6.2} Assume that $C^\bullet$ is semisimple. Then 
$r=\dsize \underset{n\in \Bbb Z}\to \sum (-1)^{n+1} \dim_{\Bbb Q_p}H^n(C^\bullet)^{\Gamma_1}$
and there exists a commutative diagram
$$
\CD
{\det}_{\Cal H} C^\bullet @>X^{-r}i_\infty >>\Cal H_0\\
@V\otimes_{\Bbb Q_p}^{\bold L}VV  @VVV\\
{\det}_{\Bbb Q_p}{C^\bullet_0} @>\vartheta>> \Bbb Q_p
\endCD
$$
in which the right vertical arrow is the augmentation map.
\endproclaim
\demo{Proof} See \cite{BG}, Lemma 8.1. Remark that Burns and Greither consider  complexes
over $\Lambda \otimes_{\Bbb Z_p}  \Bbb Q_p$ but since $\Cal H$ is a B\'ezout ring,  all their
arguments work in our case and are omitted here. 
\enddemo
\flushpar
{\bf 5.2.6.3.}  By Theorem 5.1.3 the complex $\RG_{\text{Iw},h}^{(\eta_0)}(D,V)$
is semisimple and  the first assertion follows from Lemma 5.2.6.2 together with Corollary  5.1.5.
\flushpar
$\,$
\flushpar
{\bf 5.2.6.4.} Now we can prove Theorem 5.2.5. Define

$$
\RG_{f} (\Bbb Q_v,N,T)= \RG_{\text{\rm Iw},f}^{(\eta_0)}(\Bbb Q_v,N,T)\otimes_{\Lambda}^{\bold L} \Bbb Z_p,
\qquad  \RG_{f} (\Bbb Q_v,D,V)=\RG_{f} (\Bbb Q_v,N,T)\otimes_{\Bbb Z_p}\Bbb Q_p.
$$
Remark that for $v=p$ this definition coincides with the definition given in 5.2.2.
Applying $\otimes^{\bold L}_{\Cal H}\Bbb Q_p$ to   the map  
$\RG_{\text{\rm Iw},f}^{(\eta_0)} (\Bbb Q_v,D,V) @>>> \RG_{\text{\rm Iw}}^{(\eta_0)}(\Bbb Q_v,T)
\otimes^{\bold L}_{\Lambda} \Cal H$ we obtain a morphism
$$
\RG_{f} (\Bbb Q_v,D,V)@>>> \RG (\Bbb Q_v, V).
$$ 
If $v\ne p,$ then $\RG_{f} (\Bbb Q_v,D,V)=\RG_f(\Bbb Q_v,V)$ and this morphism coincides
with the natural map $\RG_f(\Bbb Q_v,V)@>>> \RG (\Bbb Q_v,V).$
If $v=p,$ then $\RG_{f} (\Bbb Q_v,D,V)=D[-1]$ and by Theorem 2.2.4 it coincides with
the composition
$$
D@>\frac{1-p^{-1}\Ph^{-1}}{1-\Ph}>> \Dc (V)@>(h-1)! \exp_{V,\Bbb Q_p}>> H^1(\Bbb Q_p,V).
$$

Let $\RG_{f,h} (D,V)$ denote the Selmer complex associated to the diagram
$$
\xymatrix{
\RG_{S} (V)
\ar[r] & \underset{v\in S} \to  \oplus \RG (\Bbb Q_v,V)\\
 & \underset{v\in S}\to \oplus \RG_{f}(\Bbb Q_v,D,V)
 \ar[u]
}
$$

Then we have a distinguished triangle
$$
\RG_{f,h} (D,V) @>>>\RG_{S} (V)\oplus \left (\underset{v\in S}\to \oplus \RG_{f}(\Bbb Q_v,D,V)\right )@>>>\underset{v\in S} \to  \oplus \RG (\Bbb Q_v,V)
\tag{5.11}
$$
which induces isomorphisms
$$
\align
&{\det}^{-1}_{\Bbb Q_p}\RG_S(V) \otimes_{\Bbb Q_p}\left (\underset{v\in S}\to \otimes {\det}_{\Bbb Q_p}\RG (\Bbb Q_v,V)\right )
\otimes {\det}_{\Bbb Q_p}D \iso {\det}^{-1}_{\Bbb Q_p}\RG_{f,h} (D,V),\\
& \xi_{D,h} \,:\, \Delta_{\text{\rm EP}}(V)\otimes_{\Bbb Q_p} \left (  {\det}_{\Bbb Q_p}D\otimes {\det}^{-1}_{\Bbb Q_p}V^+\right )\iso 
{\det}^{-1}_{\Bbb Q_p}\RG_{f,h} (D,V).
\endalign
$$
\flushpar
Next, $\RG_{f,h} (D,V)=\RG_{\text{\rm Iw},h}^{(\eta_0)} (D,V) \otimes_{\Cal H}\Bbb Q_p$
and for any $i$ one has an exact sequence
$$
0@>>>\bold R^i\Gamma_{\text{\rm Iw},h}^{(\eta_0)} (D,V)_{\Gamma_1} @>>>\bold R^i\Gamma_{f,h} (D,V)
@>>>\bold R^{i+1}\Gamma_{\text{\rm Iw},h}^{(\eta_0)} (D,V)^{\Gamma_1}@>>>0.
$$
From Theorem 5.1.3 it follows that 
$$
\bold R^i\Gamma_{f,h} (D,V)=\cases \bold R^{2}\Gamma_{\text{\rm Iw},h}^{(\eta_0)} (D,V)^{\Gamma_1}
&{\text{\rm if $i=1$}}\\
\bold R^{2}\Gamma_{\text{\rm Iw},h}^{(\eta_0)} (D,V)_{\Gamma_1} &{\text{\rm if $i=2$}}\\
 0 &{\text{\rm if $i\ne 1,2$}}.
\endcases
$$
Therefore, the   isomorphism $\bold R^2\Gamma_{\text{\rm Iw},h}(D,V)^{\Gamma_1}@>>>
\bold R^2\Gamma_{\text{\rm Iw},h}(D,V)_{\Gamma_1}$ induces a canonical trivialization
$$
\vartheta_{D,h}\,:\,{\det}_{\Bbb Q_p}\RG_{f,h} (D,V)  \iso  \Bbb Q_p.
$$
 By Lemma 5.2.6.2
we have a commutative diagram
$$
\CD
{\det}_{\Cal H}^{-1}\RG_{\text{\rm Iw},h}^{(\eta_0)} (D,V) @>X^{-e}i_{V,\text{\rm Iw},h}>>
{\Cal H}_0
\\
@V\,^{\bold L}\otimes_{\Cal H}\Bbb Q_p VV    @VVV \\
{\det}_{\Bbb Q_p}^{-1}\RG_{f,h} (D,V)  @>\vartheta_{D,h}^{-1}>>  \Bbb Q_p.
\endCD
$$
Since 
$$
\Delta_{\text{\rm Iw},h}(N,T)\otimes^{\bold L}_{\Lambda}\Bbb Z_p \simeq \Delta_{\text{\rm EP}}(T)\otimes_{\Bbb Z_p}
\omega_N \otimes_{\Bbb Z_p} (\omega^{+}_T)^{-1}
$$ 
it implies that 
$$
\vartheta_{D,h}^{-1}\circ  \xi_{D,h} (\Delta_{\text{\rm EP}}(T)\otimes_{\Bbb Z_p}
\omega_N \otimes_{\Bbb Z_p} (\omega^{+}_T)^{-1} )    \,=\,
\log (\chi (\gamma))^{-e} L^*_{\text{\rm Iw},h}(T,N,0)\,\Bbb Z_p.
\tag{5.12}
$$ 
Consider the  diagram
$$
\xymatrix{
\RG_f(V) \ar[r] &\RG_S(V)\oplus \underset{v\in S\cup\{\infty\}}\to \oplus \RG_f(\Bbb Q_v,V) \ar[r]&
 \underset{v\in S\cup\{\infty\}}\to \oplus \RG(\Bbb Q_v,V)\\
\RG_{f,h}(D,V) \ar[r]\ar[u] &\RG_S(V)\oplus \underset{v\in S}\to \oplus \RG_f(\Bbb Q_v,D,V) \ar[r] \ar[u]
&\underset{v\in S}\to \oplus \RG(\Bbb Q_v,V) \ar[u]\\
L \ar[r]\ar[u] &S\oplus V^+[-1] \ar[r] \ar[u] &V^+[-1] \ar[u]
}
\tag{5.13}
$$ 
where $L=\text{\rm cone} \left (\RG_{f,h}(D,V) @>>> \RG_f(V) \right )\,[-1]$ and the upper
and  middle rows coincide with (4.1) and (5.11) up to the following modification:
the  map ${\text{\rm loc}}_p\,:\, \RG_f(\Bbb Q_p,V)@>>> \RG (\Bbb Q_p,V)$ 
is replaced by
$\Gamma (h)\,{\text{\rm loc}}_p$. Hence $S$ is isomorphic to $L$ in the 
derived category $\Cal D^p(\Bbb Q_p)$  and we have an exact triangle
$$
S@>>> \RG_{f,h}(D,V)@>>>\RG_f(V)@>>>S[1].
$$
An easy diagram search shows that $H^1(S)\simeq \bold R^1 \Gamma_{f,h}(D,V)$ coincides
with ${\text{\rm id}}\,:\,D^{\Ph=p^{-1}} @>>> D^{\Ph=p^{-1}}$ and that 
$$
0@>>>H^1_f(V)@>>>H^2(S)@>>>\bold R^2\Gamma_{f,h}(D,V)@>>>0
$$
coincides with
$$
0@>>>H^1_f(V)@>>>\frac{\Dc (V)}{\F^0\Dc (V)+D_{-1}}@>\Gamma (h)\exp_{V}
>>\frac{H^1(\Bbb Q_p,V)}{H^1_{f,\{p\}}(V)+H^1(\bD_{-1})}@>>>0.
$$
Therefore, we have a commutative diagram 
$$
\CD
{\det}_{\Bbb Q_p} S \otimes {\det}_{\Bbb Q_p}\RG_f(V) @> \alpha>>  {\det}_{\Bbb Q_p}\RG_{f,h} (D,V)\\
@V{\vartheta_S}VV  @V{\vartheta_{D,h}}VV\\
\Bbb Q_p @>\kappa>> \Bbb Q_p
\endCD
$$
where  $\vartheta_S$ was defined in section 5.2.2 and   $\kappa$  is the unique map which makes this diagram commute.
%$$
%\Bbb Q_p \iso {\det}^{-1}_{\Bbb Q_p} H^1(S)\otimes   {\det}_{\Bbb Q_p} H^2(S)
%\iso   {\det}^{-1}_{\Bbb Q_p}\bold R^1 \Gamma_{f,h}(D,V) \otimes
% {\det}_{\Bbb Q_p}\bold R^2 \Gamma_{f,h}(D,V) \iso \Bbb Q_p
%$$
From Proposition 3.2.4  and  Corollary 5.1.5 we obtain immediately that 
$$
\kappa\,=\,(\log \chi (\gamma))^{e} \left (1-\dsize\frac{1}{p}\right )^{e}
\scr L(V,D)^{-1}\, {\text{\rm id}}_{\Bbb Q_p}.
\tag{5.14}
$$
Passing to determinants in the diagram (5.13) we obtain  a commutative diagram
\eightpoint
$$
\xymatrix{
\Delta_{\text{\rm EP}} (V) \otimes 
\left ({\det}_{} (t_V(\Bbb Q_p))\otimes {\det}^{-1}_{}V^+ \right )
\otimes {\det}\,\RG_f (V)
\ar[d]^{\alpha_S} \ar[rr]   & & \Bbb Q_p \ar @{=}[d]\\
\Delta_{\text{\rm EP}} (V) \otimes 
\left ({\det}  D \otimes {\det}_{}S \otimes {\det}\,\RG_f (V)\right )\otimes 
{\det}^{-1}_{}V^+  \ar[r]^-{\xi_{D,h}\otimes 
\alpha}  \ar[d]^{\text{\rm id}\otimes \beta_S}&{\det}_{}^{-1}
\RG_{f,h} (D,V) \otimes {\det}_{}
\RG_{f,h} (D,V)  \ar[r]^-{\text{\rm duality}} \ar[d]^{{\text{\rm id}}\otimes \vartheta_{D,h}}& \Bbb Q_p \ar @{=}[d]\\
\Delta_{\text{\rm EP}}^{} (V) \otimes 
\left ({\det}_{} D \otimes {\det}^{-1}_{\Bbb Q_p}V^+ \right )
\ar[r]^-{\xi_{D,h}\otimes \kappa} &{\det}_{}^{-1}
\RG_{f,h} (D,V) \ar[r]^-{\vartheta_{D,h}^{-1}}    &\Bbb Q_p}
$$
\tenpoint
where the maps $\alpha_S$ and $\beta_S$   were defined in (5.7-5.9). 
The upper row of this diagram sends $\Delta_{EP}(T)\otimes  (\omega_{t_M}\otimes (\omega^{+}_T)^{-1}
\otimes \omega_f)$ onto 
$$\Gamma (h)^{d_{+}(V)}\,\frac{ i_{\omega_M,p} (\Delta_{EP}(T))}
{\Omega_M^{(\acute et,p)} (\omega_T^+,\omega_{\text{\rm B}}^+)}.
\tag{5.15}$$
 From Lemma 5.2.3  it follows that the composition of left vertical maps
$\vartheta_S=(\text{\rm id}\otimes \beta_S)$ sends 
$\Delta_{EP}(T)\otimes  (\omega_{t_M}\otimes (\omega^{+}_T)^{-1}
\otimes \omega_f)$ onto 
$${\det}^* \left (\frac{1-p^{-1}\Ph^{-1}}{1-\Ph}\mid D \right )^{-1} E_p(V,1)^{-1}
R_{V,D}(\omega_{V,N})\,\,
\Delta_{EP}(T)\otimes  (\omega_N\otimes (\omega^{+}_T)^{-1})
 \tag{5.16}
$$
Next,  (5.12) and (5.13) give 
$$
\vartheta_{D,h}^{-1}\circ (\xi_{D,h}\otimes \kappa) (\Delta_{\text{\rm EP}}(T)\otimes \omega_N \otimes (\omega^{+}_T)^{-1})\,=\, 
\left (1-\dsize\frac{1}{p}\right )^{e}
\scr L(V,D)^{-1}      L^*_{\text{\rm Iw},h}(T,N,0)\,\Bbb Z_p.
\tag{5.17}
$$

Putting together (5.15), (5.16) and (5.17) we obtain that
$$
\frac{L_{\text{\rm Iw},h}^*(T,N,0)}{R_{V,D}(\omega_{V,N})}
\sim_p \Gamma (h)^{d_{+}(V)} \,\scr L (V,D)\, 
\,E_p^*(V,1)\,{\det}_{\Bbb Q_p} \left (\frac{1-p^{-1}\Ph^{-1}}{1-\Ph} \,\vert D_{-1} \right ) 
\frac{i_{\omega_M, p}\,(\Delta_{\text{\rm EP}} (T))}{\Omega_M^{(\acute et,p)} (\omega_T^+,\omega_{\text{\rm B}}^+)}
$$
and the theorem is proved.
\newline
\newline
{\bf 5.2.7. Special values of $L_{\text{\rm Iw},h}^*(T,N,s).$}
Let $\widetilde H^1_f(T)$ denote the image of $H^1_f(T)$ in $H^1_f(V)$ and let $\omega_{T,f}$
be a base of $\det_{\Bbb Z_p}\widetilde H^1_f(T).$ Let $R_{V,D}(\omega_{T,N})$ denote
the determinant of $r_{V,D}$ computed in the bases $\omega_{t_{M}}$, $\omega_N$ and
$\omega_{T,f}.$

\proclaim{Corollary 5.2.8} Under the assumptions of Theorem 5.2.5 one has 
$$
\frac{L_{\text{\rm Iw},h}^*(T,N,0)}{R_{V,D}(\omega_{T,N})}
\sim_p \Gamma (h)^{d_{+}(V)} \,\scr L (V,D)\, 
\Cal E^+(V,D)
\,\,\frac{\#\sha(T^*(1))\,\Tam_{\omega_M}^0 (T)}{\#H^0_S(V/T)\,\#H^0_S(V^*(1)/T^*(1))},
$$
where $\sha(T^*(1))$ is the Tate-Shafarevich group of Bloch-Kato \cite{BK} and   $\Tam_{\omega_M}^0 (T)$ is the product of local Tamagawa numbers of $T$ taken over all primes and computed with respect
to a fixed base $\omega_{t_M}$ of ${\det}_{\Bbb Q}t_M(\Bbb Q).$ 
\endproclaim
\demo{Proof} The  computation of the trivialisation of the Euler-Poincar\'e line
(see for example \cite{FP}, chapitre II, Th\'eor\`eme 5.6.3) together with the definition
of $i_{\omega_M,p}$  by (4.12) give
$$
i_{\omega_M,p}(\Delta_{EP}(T))\,=\,\frac{\#\sha(T^*(1))\,\Tam_{\omega_M}^0 (T)}{\#H^0_S(V/T)\,\#H^0_S(V^*(1)/T^*(1))}\,\Omega_M^{(\acute et,p)} (\omega_T^+,\omega_{M_{\text{\rm B}}}^+)\,   [\omega_f:\omega_{T,f}].
$$
Since $R_{V,D}(\omega_{T,N})=R_{V,D} (\omega_{V,N})\,[\omega_f:\omega_{T,f}]$
the corollary follows from Theorem 5.2.5.
\enddemo
{\,}
\flushpar
{\bf 5.2.9. The functional equation.} %To each lattice $T$ of $V$ stable under
%the action of $G_{\Bbb Q_p}$ Berger associated a canonical lattice 
%$\Dc (T) \subset \Dc(V)$ using his theory of Wach modules \cite{Ber2}. 
%In particular, the determinant of the comparision isomorphism
%$$
%\Dc (V)\otimes_{\Bbb Q_p} \Bc \simeq V\otimes_{\Bbb Q_p} \Bc
%$$
%computed in bases of $\Dc (T)$ and $T$ is in $W(\bar\Bbb F_p)^*t^{t_N(V)}$
%where $t_N(V)=\underset_{i\in \Bbb Z}\to\sum i\dim_{\Bbb Q_p}(\text{\rm gr}_i\Dc (V)).$
%This allows us, if $T$ is fixed, to associate to any regular submodule $D\subset \Dc (V)$
%the canonical lattice $N=D\cap \Dc (T)$ of $D$. If $D^{\perp} \subset \Dc (V^*(1))$
%and $N^{\perp}=D^{\perp}\cap \Dc (V^*(1))$
%then
Recall that we set $h_i(V)=\dim_{\Bbb Q_p}(\text{\rm gr}_i\Dd (V))$ and 
$m=\underset_{i\in \Bbb Z}\to\sum i h_i(V).$ Since $V$ is crystalline,
${\det}_{\Bbb Q_p}(V)$ is a one dimensional crystalline representation and 
 ${\det}_{\Bbb Z_p}(T)=  T_0(m)$ where  
$T_0$ is an unramified  $G_{\Bbb Q_p}$-module of rank $1$ over $\Bbb Z_p.$ The module
$(T_0\otimes W(\overline{\Bbb F}_p))^{\Ph=1}\,e_m$ where $e_m=(t^{-1}\otimes\varepsilon )^{\otimes m}$
is a  $\Bbb Z_p$-lattice in ${\det}_{\Bbb Q_p}(\Dc (V))=\Dc (V_0(m))$
which depends only on $T$ and which we denote by $\Dc (T_0(m)).$

Let $D^{\perp}$
be the dual regular module.  The exact sequence
$$
0@>>>D@>>>\Dc (V)@>>>\left (D^{\perp}\right )^*@>>>0
$$ 
gives an isomorphism
$$
{\det}_{\Bbb Q_p}D\otimes {\det}^{-1}_{\Bbb Q_p}D^{\perp}
\simeq {\det}_{\Bbb Q_p}\Dc (V)
$$
and we fix a lattice $N^{\perp}\subset D^{\perp}$ such that
$${\det}_{\Bbb Z_p}N\otimes {\det}^{-1}_{\Bbb Z_p}N^{\perp}
\simeq \Dc (T_0(m)). 
$$ 
Set $\Gamma_{V,h}(s)=\underset{j>-h}\to\prod (j+s)^{\dim\F^{j}\Dd (V)}.$
%We also define 
%$$
%\Gamma^* (V)=\prod_{i\in \Bbb Z} \Gamma^*(-i)^{h_i(V)}
%$$
%where 
%$$
%\Gamma (i)=\cases (i-1)! &\text{if $i>0$}\\
%\frac{(-1)^i}{(-i)!} &\text{if $i\leqslant 0$}.
%\endcases
%$$
The conjecture $\delta_{\Bbb Z_p}(V)$ of \cite{PR1} proved in \cite{BB} implies
that for $h\gg 0$
$$
L_{\text{\rm Iw},h}(T^*(1),N^{\perp},-s) \sim_{\Lambda^*}
\Gamma_{V,h}^{-1}(s) \underset{-h<j<h}\to\prod (j+s)^{d_{+}(V^*(1))}
L_{\text{\rm Iw},h}(T,N,s) 
$$
(see \cite{PR2}, Th\'eor\`eme 2.5.2).
This can seen as the algebraic counterpart of the functional equation
for $p$-adic $L$-functions.
An elementary computation (see \cite{BB}, Lemme 4.7) shows that
$$
\Gamma_{V,h}^{-1}(s) \underset{-h<j<h}\to\prod (j+s)^{d_{+}(V^*(1))}
=\Gamma^*(V) \Gamma (h)^{d_{+}(V^*(1))-d_{+}(V)} s^r +o (s^r). 
$$
where $r={\dim_{\Bbb Q_p} t_V(\Bbb Q_p)-d_+(V)}=\dim_{\Bbb Q_p} H^1_f(V)$
and $\Gamma^*(V)$ is defined by (4.14).
Therefore $L_{\text{\rm Iw},h}(T^*(1),N^{\perp},s)$ has a zero of order
$\dim_{\Bbb Q_p} H^1_f(V)+e$  at $s=0$. Moreover one has
$$
\frac{L^*_{\text{\rm Iw},h}(T^*(1),N^{\perp},0)}{\Gamma (h)^{d_{+}(V^*(1))}}\,\,
\sim_p \,\,
\Gamma^*(V) \,
\frac{L_{\text{\rm Iw},h}(T,N,0)}{\Gamma (h)^{d_+(V)}}. 
$$
From the definition of $R_{V^*(1),D^{\perp}}$ (see section 4.2.1) one obtains easily
that $R_{V^*(1),D^{\perp}}(\omega_{V^*(1),N^{\perp}})=\Omega_M^{(H,p)}(\omega_T,
\omega_{M_{\text{\rm dR}}})^{-1} R_{V,D}(\omega_{V,N})$ where $\Omega_M^{(H,p)}$ denotes the period map defined 
(4.8) and (4.9) and $\omega_{M_{\text{\rm dR}}}=\omega_{t_M}\otimes \omega_{t_{M^*(1)}}^{-1}.$
Taking into account (4.15) we obtain that
$$
\frac{L^*_{\text{\rm Iw},h}(T^*(1),N^{\perp},0)}
{R_{V^*(1),D^{\perp}}(\omega_{V,N^{\perp}})}\,
\sim_p \,\Gamma (h)^{d_{+}(V^*(1))} \scr L(V,D)\,\Cal E^+(V^*(1),D^{\perp})\,
\frac{i_{\omega_M, p}\,(\Delta_{\text{\rm EP}} (T^*(1)))}{\Omega_M^{(\acute et,p)} (\omega_{T^*(1)}^+,\omega_{M^*(1)_{\text{\rm B}}}^+)},
$$
which is the analog of Theorem 5.2.5 for $L_{\text{\rm Iw},h}(T^*(1),N^{\perp},s).$

\head  {\bf Appendix. Galois cohomology of $p$-adic
representations}
\endhead

{\bf A.1.} Let $K$ be a finite extension of $\Bbb Q_p$ and $T$  a
$p$-adic representation of $G_K.$ Fix a topological generator $\gamma$ of $\Gamma$.
Let $\bD (T)=(T\otimes_{\Bbb Z_p} \bold A)^{H_K}$ be the
$(\Ph,\Gamma)$-module associated to $T$ by Fontaine's theory \cite{F2}. Consider the complex
$$
C_{\Ph,\g}(\bD (T))\,=\,\left [\bD(T)@>f>> \bD (T)\oplus \bD (T) @>g>>
\bD (T) \right ]
$$
where the  modules are placed in degrees $0$, $1$ and $2$ and the
maps $f$ and $g$ are given by
$$
f(x)= ((\Ph-1)\,x\,,\, (\g-1)\,x),\qquad g(y,z)= (\g-1)\,y\,-\,
(\Ph-1)\,z.
$$

\proclaim{Proposition A.2} There are canonical and functorial
isomorphisms
$$
h^i\,:\, H^i(C_{\Ph,\g}(\bD (T))) \iso H^i(K,T)
$$
which can be described explicitly by the following formulas:

i) If $i=0,$ then $h^0$ coincides with the natural isomorphism
$$
\bD (T)^{\Ph=1,\g=1}\,=\,H^0(K,T\otimes_{\Bbb Z_p}\A^{\Ph=1})\,=\,
H^0(K,T).
$$
ii)  Let $\alpha,\beta \in \bD (T)$ be such that $(\g-1)\,\alpha
\,=\,(1-\Ph)\,\beta.$ Then $h^1$ sends $\cl (\alpha,\beta)$ to the
class of the cocycle
$$
\mu_1 (g)\,=\, (g-1)\,x\,+\,\dsize \frac{g-1}{\g-1}\,\beta,
$$
where $x \in \bD (T)\otimes_{\A_K} \A$ is a solution of the
equation $(1-\Ph)\,x\,=\,\alpha .$

iii) Let $\widehat\g\in G_K$ be a lifting of $g\in \Gam$ and let
$x$ be a solution of $(\Ph -1)\,x\,=\,\alpha .$ Then  $h^2$ sends
$\alpha$ to the class of the 2-cocycle
$$
\mu_2 (g_1,g_2)\,=\,\widehat \g^{k_1}(h_1-1)\,\frac{\widehat
\g^{k_2}-1}{\widehat \g-1}\,x
$$
where $g_i\,=\,\hat\g^{k_i}h_i,$ $h_i\in H_K.$
\endproclaim
\demo{Proof} The isomorphisms $h^i$ were constructed   in
\cite{H1}, Theorem 2.1. Remark that i) follows directly from this
construction (see \cite{H1}, p.573) and that ii) is proved in
\cite{Ben1}, Proposition 1.3.2 and \cite{CC2}, Proposition I.4.1.
The proof of iii) follows along exactly the same lines. Namely, it
is enough to prove this formula modulo $p^n$ for each $n$. Let $\alpha \in \bD
(T)/p^n\bD (T) .$ By Proposition 2.4 of \cite{H1} there exists
$r\ge 0$ and $y\in \bD(T)/p^n\bD (T)$ such that
$(\Ph-1)\,\alpha\,=\,(\g-1)^r\beta .$ Let
$$
N_x=(\bD (T)/p^n\bD (T))\oplus (\oplus_{i=1}^r
(\A_K/p^n\A_K)\,t_i),
$$
where $\Ph (t_i)\,=\,t_i+(\g-1)^{r-i}(\alpha)$ and $\g
(t_i)\,=\,t_i+t_{i-1}.$ Then $N_x$ is a $(\Ph,\Gamma)$-module and
we have a short exact sequence
$$
0@>>>\bD@>>>N_x@>>>X@>>>0
$$
where $X=N_x/M\simeq \oplus_{i=1}^r \A_K/p^n\A_K \bar t_i.$ An
easy diagram search shows that the connecting homomorphism
$\delta^1_{\bD}\,:\,H^1(C_{\Ph,\g}(\bD (X))) @>>>H^2(C_{\Ph,\g}(\bD(T)))$
sends $\cl (0,\bar t_r)$ to $-\cl (\alpha).$ The functor 
$\bold V(D)=(D\otimes_{\bold A_K} \bold A)^{\Ph=1}$ is a quasi-inverse 
to $\bD $. Thus one has an exact sequence of Galois
modules
$$
0@>>> T/p^nT@>>> T_x @>>> \bold V(X)@>>>0
$$
where $T_x=\bold V(N_x).$ From the definition of $x$ it follows
immediately that $t_r-x\in T_x .$ By ii), $h^1(\cl (0,\bar t_r))$
can be represented by the cocycle $c(g)\,=\,
\dsize\frac{g-1}{\g-1}\,\bar t_r$ and we fix its lifting $\hat
c\,:\,G_K@>>>N_x$ putting $\hat c(g)\,=\, \dsize\frac{g-1}{\g-1}\,
(t_r-x).$ As $g_1\hat c(g_2)-\hat c(g_1g_2)+\hat c(g_1)\,=\, -\mu_2(g_1,g_2), $
the connecting map $\delta_T^1\,:\,H^1(K,\bold
V(X))@>>>H^2(K,T/p^nT)$ sends $\cl (c)$ to $-\cl (\mu_2)$ and iii)
follows from    the commutativity of the diagram
$$
\CD
H^1(C_{\Ph,\g}(X)) @>\delta^1_{\bD}>>H^2(C_{\Ph,\g}(T/p^nT))\\
@Vh^1VV @Vh^2VV\\
H^1(K,\bold V(X))@>\delta^1_T>>H^2(K,T/p^nT).
\endCD
$$
\enddemo

%{\bf A.1.2.}  Let $A^{\bullet}(T)$ denote the complex
%$
%\bD(T) @>\g-1>> \bD (T)
%$
%where the modules are placed in degrees $0$ and $1$. Then
%$
%C_{\Ph,\g}(T) \,=\,\Tot^{\bullet}(A^\bullet (T)@>{\Ph-1}>>
%A^\bullet (T)).
%$

%Consider the ring extension $\A/\Ph (\A).$  The rings $\A$  and
%$\Ph (\A)$ are absolutely unramified, but the residual extension
%is purely inseparable of degree $p$ and we can define a left
%inverse $\psi\,\,:\,\, \A@>>>\A$ for $\Ph $ by
%$$
%\psi (x)=\frac{1}{p}\,\Ph^{-1}\left ( \Tr_{\A/\Ph (\A)}(x) \right
%).
%$$

%Since $\psi$ commutes with the action of $G_K$, it acts on $\bD
%(T)$ and we can consider the complex:
%again in degrees $0,1,2$.
%Recall that we denote by $\RG (K,T)$ the standard complex of
%continuous cochains
%$$
%C^{\bullet}_{\co} (G_K,T)\,=\, \left [
%C^0_{\co}(G_K,T)@>>>C^1_{\co}(G_K,T)@>>>C^2_{\co}(G_K,T)@>>>\hdots
%\right ].
%$$

\proclaim{Proposition A.3} The complexes $\RG (K,T)$ and  $C_{\Ph,\g}(T)$ 
%and$C_{\psi,\g}(T)$ 
are isomorphic in  $\Cal D(\Bbb Z_p)$. 
%Let 
%$$
% C _{\psi,\g}(T)\,=\,\bigl [\,\bD(T)@>f_1>> \bD (T)\oplus \bD
%(T) @>g_1>>\bD (T)\,\bigr ] 
%$$
%be the complex defined by 
%$f_1 (x)= ((\psi-1)\,x,\,(\g-1)\,x)$, 
%$g_1(y,z)= (\g-1)\,y\,-\,(\psi-1)\,z.$
%Then

%i) The diagram
%$$\xymatrix{C_{\Ph,\g}(T)\,\,:  &\bD (T) \ar[r] \ar[d]^{\alpha_0}
%&\bD (T)\oplus \bD (T)
%\ar[r] \ar[d]^{\alpha_1} &\bD (T) \ar[d]^{\alpha_2}\\
%C_{\psi,\g}(T)\,\,:  &\bD (T)\ar [r] &\bD (T) \oplus \bD (T)\ar
%[r] &\bD (T)}
%$$
%where $\alpha_0 (x)= x,$ $\alpha_1(y,z)=(-\psi (y),z),$ and
%$\alpha_2 (z)=-\psi (z),$ commutes and gives rise to a quasi-isomorphism
%$\alpha_* \,\,:\,\, C_{\Ph,\g}(T)@>>>C_{\psi, \g}(T)$.

\endproclaim
\demo{Proof} %i) The commutativity of the diagram follows
%immediately from the relation $\psi \circ \Ph= \text{\rm id}$ and
%we only need to check that in each degree $\alpha_i$ induces an
%isomorphism of $H^i(C_{\Ph,\g}(T))$ onto $H^i(C_{\psi,\g}(T)).$ If
%$i=0$, one has $H^0(C_{\Ph,\g}(T))=\bD (T)^{\Ph=1,\g=1},$
%$H^0(C_{\psi,\g}(T))=\bD (T)^{\psi=1,\g=1}$ and
%$H^0(C_{\Ph,\g}(T))\subset H^0(C_{\psi,\g}(T)).$ Now if $x\in \bD
%(T)^{\psi=1,\g=1},$ then $\psi ((1-\Ph)x)=0$ and
%$(\g-1)\,(1-\Ph)\,x=0.$ Recall that by \cite{CC2} the map
%$\g-1\,:\,\bD (T)^{\psi=0} @>>>\bD (T)^{\psi=0}$ is an
%isomorphism. Thus $(1-\Ph)\,x=0$ and $x\in \bD (T)^{\Ph=1,\g=1}.$
%The isomorphism between $H^1(C_{\Ph,\g}(T))$ and
%$H^1(C_{\psi,\g}(T))$ was established in \cite{CC2}.
%Consider the map $H^2(C_{\Ph,\g}(T))@>>> H^2(C_{\psi,\g}(T))$
%given by $\text{cl} (z)=-\text{cl}(\psi (z)).$ As $-\psi
%(z)=-z+(1-\psi)\,z$ it is surjective. On the other hand, if
%$\text{cl}(z)$ belongs to the kernel of this map, then there
%exists $a,b\in \bD (T)$ such that
%$
%-\psi (z)= (\g-1)\,a - (\psi-1)\,b.
%$
%Put $c=z+(\g-1)\,\Ph (a)\, +\,(\Ph-1)\,b.$ Then
%$
%\psi (c)=\psi (z)+(\g-1)\,a \,+\,(1-\psi )\,b=0.
%$
%Thus there exists $a'\in \bD (T)$ such that $c=(\g-1)\,a'$ and we obtain that
%$
%z=(\g-1)\,(a'-\Ph(a))-(\Ph-1)\,b.
%$
%This implies that $\text{cl} (z)=0$ and we proved that $H^2(C_{\Ph,\g}(T))@>>>
%H^2(C_{\psi,\g}(T))$ is an isomorphism.
The proof  is standard (see for example \cite{BF},
proof of Proposition 1.17). The exact
sequence
$$
0@>>> T@>>>\bD (T)\otimes_{\A_K}\A@>\Ph-1>>\bD
(T)\otimes_{\A_K}\A@>>>0
$$
gives rise to an exact sequence of complexes
$$
0@>>> C^{\bullet}_{\co}(G_K,T)@>>>C^{\bullet}_{\co}(G_K,\bD
(T)\otimes_{\A_K}\A)@>\Ph-1>>C^{\bullet}_{\co} (G_K,\bD
(T)\otimes_{\A_K}\A)@>>>0
$$
Thus $\bold R\Gamma (K,T)$ is quasi-isomorphic to the total
complex
$$
K^\bullet (T)\,=\,\text{Tot}^{\bullet} \left
(C^{\bullet}_{\co}(G_K,\bD
(T)\otimes_{\A_K}\A)@>\Ph-1>>C^{\bullet}_{\co} (G_K,\bD
(T)\otimes_{\A_K}\A)\right ).
$$
On the other hand
$
C_{\Ph,\g}(T)\,=\,\text{Tot}^{\bullet} \left
(A^{\bullet}(T)@>\Ph-1
>>A^{\bullet}(T) \right ),
$
where $A^{\bullet}(T)=[\bD (T) @>\gamma-1>>\bD (T)]$.
Consider the following commutative diagram of complexes
$$
\xymatrix{ \bD (T) \ar[d]^{\beta_0} \ar[r]^{\g-1} &\bD (T) \ar[d]^{\beta_1} \ar[r] &0 \ar[r]\ar[d] & \cdots \\
C^{0}(G_K,\bD (T)\otimes_{\A_K}\A) \ar[r]& C^{1}(G_K,\bD
(T)\otimes_{\A_K}\A) \ar[r]& C^{2}(G_K,\bD (T)\otimes_{\A_K}\A)
\ar[r]& \cdots }
$$
in which  $\beta_0 (x)=x$ viewed as a constant function on $G_K$ and
$\beta_1(x)$ denotes the map $G_K@>>>\bD (T)\otimes_{\A_K}\A)$
defined by
$
(\beta_1 (x))\,(g)=\dsize\frac{g-1}{\g-1}\, x.
$
This diagram induces a map $\Tot^{\bullet}
(A^\bullet(T)@>\Ph-1>>A^{\bullet}(T)) @>>>K^\bullet (T)$ and we
obtain a diagram
$$
C_{\Ph,\g}(T) @>>>K^{\bullet}(T) \leftarrow \RG (K,T)
$$
where the right map is a quasi-isomorphism. Then for each $i$ one
has a map
$$
H^i(C_{\Ph,\g}(T)) @>>>H^i(K^{\bullet}(T)) \simeq H^i(K,T)
$$
and an easy diagram search shows that it coincides with $h^i.$ The
proposition is proved.
\enddemo
\proclaim{Corollary A.4} Let $V$ be a $p$-adic representation of $G_K$.
Then the complexes $\bold R\Gamma (K,V)$, $C_{\Ph,\g}(\bD^{\dag}(V))$ and
$C_{\Ph,\g}(\Ddagrig (V))$ are isomorphic in $\Cal D(\Bbb Q_p).$
\endproclaim
\demo{Proof} This follows from Theorem 1.1 of \cite{Li} together with Proposition A.2.
\enddemo 

{\bf A.5.} Recall that $K_\infty/K$ denotes the cyclotomic
extension obtained by adjoining all $p^n$-th roots of unity. Let
$\Gamma =\G (K_\infty/K)$ and let $\Lambda (\Gamma)\,=\,\Bbb Z_p[[\Gamma]]$
denote the Iwasawa algebra of $\Gam .$ For any  $\Bbb Z_p$-adic
representation  $T$ of $G_K$ the induced representation
$\Ind_{K_{\infty}/K}T$ is isomorphic to $(T\otimes_{\Bbb Z_p}\Lambda
(\Gamma))^\iota $ and we set
$
\R\Gamma_{\text {\rm Iw}} (K,T)\,=\,C^{\bullet}_{\co}
(G_K\,,\Ind_{K_{\infty}/K}T).
$
Consider the complex
$$
C_{\Iw, \psi}(T)\,=\,\left [\bD (T)@>\psi-1>>\bD (T)\right ]
$$
in which the first term is placed in degree $1$.

\proclaim{Proposition A.6} There are canonical and functorial
isomorphisms
$$
h^i_{\Iw}\,\,:\,\,H^i(C_{\Iw,\psi}(T))@>>> H^i_{\Iw}(K,T)
$$
which can be described explicitly by the following formulas:

i) Let $\alpha \in \bD (T)^{\psi=1}.$ Then $(\Ph -1)\,\alpha \in
\bD (T)^{\psi=0}$ and for any $n$ there exists a unique
$\beta_n\in \bD (T)$ such that $(\gn -1)\,\beta_n= (\Ph-1)\,\alpha
.$ The map $h^1_{\Iw}$ sends $\cl (\alpha)$ to
$
(h^1_n (\cl (\beta_n,\alpha)))_{n\in \bold N} \in H^1_{\Iw}
(K_n,T).
$

ii) If $\alpha \in \bD (T),$ then
$
h^2_{\Iw} (\cl (\alpha))\,=\,-(h^2_n(\Ph (\alpha)))_{n\in \bold
N}.
$
\endproclaim
\demo{Proof} The proposition follows from Theorem II.1.3 and
Remark II.3.2 of \cite{CC2} together with Proposition A.2.
\enddemo

\proclaim{Proposition A.7} The complexes $\R\Gamma_{\text {\rm
Iw}} (K,T)$ and $C_{\Iw, \psi}(T)$ are isomorphic in the derived
category $\Cal D (\Lambda (\Gamma)).$
\endproclaim
\demo{Proof} We repeat the arguments used in the proof of
Proposition A.1.2 with some modifications. For any $n\geqslant 1$ one
has an exact sequence
$$
0@>>>\Ind_{K_{n}/K}T @>>> (\bD(T)\otimes_{\Bbb Z_p}\Bbb
Z_p[G_n]^{\iota})\otimes_{\A_K}\A @>\Ph-1>> (\bD(T)\otimes_{\Bbb
Z_p}\Bbb Z_p[G_n]^{\iota})\otimes_{\A_K}\A @>>> 0.
$$
Set $\bD (\Ind_{K_\infty/K}T)=\bD (T)\otimes_{\Bbb Z_p}\Lambda
(\Gamma)^{\iota}$ and
$$
\bD (\Ind_{K_\infty/K}(T)) \hat \otimes_{\A_K}
\A\,=\,\varprojlim_n (\bD (T)\otimes_{\Bbb Z_p}\Bbb
Z_p[G_n]^{\iota})\otimes_{\A_K} \A.
$$
As $\Ind_{K_{n}/K}T$ are compact, taking projective limit one
obtains an exact sequence
$$
0@>>>\Ind_{K_{\infty}/K}T @>>> \bD (\Ind_{K_\infty/K}(T)) \hat
\otimes_{\A_K} \A @>\Ph-1>> \bD (\Ind_{K_\infty/K}(T)) \hat
\otimes_{\A_K} \A @>>>0 .
$$

Thus $\R_{\text {\rm Iw}} (K,T)$ is quasi-isomorphic to
$$
K^{\bullet}_{\Iw}(T)\,=\, \text{Tot}^{\bullet} \left
(C^{\bullet}_{\co}(G_K,\bD (\Ind_{K_\infty/K}T) \hat
\otimes_{\A_K} \A) @>\Ph-1>> C^{\bullet}_{\co}(G_K,\bD
(\Ind_{K_\infty/K}T) \hat \otimes_{\A_K} \A) \right ).
$$
We construct a quasi-isomorphism
$f_{\bullet}\,:\,C_{\Iw,\psi}(T)@>>> K^{\bullet}_{\Iw}(T).$ Any
$x\in \bD (T)$  can be written in the form $x=(1-\Ph
\psi)\,x\,+\,\Ph \psi (x)$ where $\psi (1-\Ph \psi)\,x\,=\,0.$
Then for each $n \geqslant  0$ the equation
$
(\gn -1)\,y_n\,=\,(\Ph \psi -1)\,x
$
has a unique solution  $y_n\in \bD (T)^{\psi=0}$ (\cite{CC2},
Proposition I.5.1). In particular,
$
y_n\,=\,\dsize \frac{\gamma_{n+1}-1}{\gn-1}\,y_{n+1}
$
and we have a compatible system of elements
\linebreak
$
Y_n=\dsize\sum_{k=0}^{|G_n|-1} \g^k \otimes \g^k(y_n)\,\in \bD (T)
\otimes_{\Bbb Z_p} \Bbb Z_p[G_n]^{\iota}.
$
Put $Y=(Y_n)_{n\ge 0}  \in \bD (\Ind_{K_\infty/K}T).$ Then
$$
(\gn-1)\,Y_n\,=\,(\g-1)\,Y \pmod{\bD (\Ind_{K_n/K}T)}.
$$
Let $\eta_x \in C^{1}_{\co}(G_K,\bD (\Ind_{K_\infty/K}T)\hat
\otimes_{\A_K} \A)$ be the map defined by
$
\eta_x (g)\,=\,\dsize \frac{g-1}{\g-1}\,(1\otimes x).
$
Define
$
f_1\,:\,\bD (T)@>>> K^1_{\Iw}(T)\,=\, C^{0}_{\co}(G_K,\bD
(\Ind_{K_\infty/K}T)\hat \otimes_{\A_K} \A) \oplus
C^{1}_{\co}(G_K,\bD (\Ind_{K_\infty/K}T)\hat \otimes_{\A_K} \A)
$
by $f_1(x)=(Y,\eta_x)$ and
$
f_2\,:\,\bD (T) @>>> C^{1}_{\co}(G_K,\bD (\Ind_{K_\infty/K}T)\hat
\otimes_{\A_K} \A)\subset K^2_{\Iw}(T)
$
by $f_2(z)= -\eta_{\Ph (z)}.$ 
It is easy to check that $f_{\bullet}$ is a morphism of complexes. This gives   a diagram
$$
C_{\Iw,\psi}(T) @>>>K^{\bullet}_{\Iw}(T) \leftarrow \RG_{\Iw}
(K,T)
$$
in which the right map is a quasi-isomorphism. Using Proposition A.1.4 it is not difficult to check 
that for each  $i$ the induced  map
$$
H^i(C_{\Iw,\psi}(T)) @>>>H^i(K^{\bullet}_{\Iw}(T)) \simeq
H_{\Iw}^i(K,T)
$$
coincides with $h^i_{\Iw}.$ The proposition is proved.
\enddemo

\proclaim{Corollary A.8} The complexes $\bold R\Gamma_{\Iw} (K,T)$ and
$C_{\Iw,\psi}^{\dag}(T)$ are isomorphic in $\Cal D(\Lambda (\Gamma)).$
\endproclaim
\demo{Proof} One has $\bD^{\dag} (T)^{\psi=1}=\bD (T)^{\psi=1}$ (\cite{CC1}, Proposition 3.3.2)
and $\bD^{\dag} (T)/(\psi-1) =\bD (T)/(\psi-1)$ (\cite{Li}, Lemma 3.6).
This shows that the inclusion $C_{\Iw,\psi}^{\dag}(T) @>>> \bD (T)^{\psi=1}$ is a quasi-isomorphism.
\enddemo

\flushpar 
{\bf Remark A.9.} These results can be slightly improved.
Namely, set $r_n=(p-1)p^{n-1}.$ The method used in the proof of Proposition III.2.1 
\cite{CC2} allows to show that $\psi (\bD^{\dag,r_n}(T))\subset \bD^{\dag,r_{n-1}}(T)$
for $n\gg 0.$ Moreover, for any $a\in  \bD^{\dag,r_n}(T)$ the solutions of 
the equation $(\psi-1)\,x=a$ are in $\bD^{\dag,r_n}(T).$ Thus
$C_{\Iw}^{\dagger,r_n} (T)\,=\, \left [\bD^{\dag,r_n}
(T)@>\psi-1>>\bD^{\dag, r_n} (T)\right ]$, $n\gg 0$ is a well-defined complex 
which is quasi-isomorphic to $C_{\Iw,\psi}^{\dag}(T).$  
Further, as $\Ph (\A^{\dagger,r/p})\,=\,\A^{\dagger,r}$ 
we can consider the complex 
$$
C^{\dagger ,r_n}_{\Ph,\g}(T)\,=\,\bigl [\bD^{\dagger
,r_{n-1}}(T)@>f>> \bD^{\dag,r_n} (T)\oplus \bD^{\dagger , r_{n-1}}
(T) @>g>> \bD^{\dag,r_n}
(T) \bigr ], \qquad n\gg 0
$$
in which $f$ and $g$ are defined by the same formulas as before.
Then the inclusion 
$
C^{\dagger, r_n}_{\Ph,\g} (T) @>>>C_{\Ph,\g}(T)
$
is a quasi-isomorphism.

\enddocument
\bye